\documentclass[11pt]{article}
\usepackage{fancyhdr}
\usepackage{bigints}
\usepackage{isomath}
\usepackage{mathtools} 
\usepackage{amsbsy}
\usepackage{amssymb}
\usepackage{amscd,amsfonts}
\usepackage{bigints}
\usepackage{graphicx}
\usepackage{verbatim}
\usepackage{euscript}
\usepackage{alltt}
\usepackage{stmaryrd}
\usepackage{relsize}
\usepackage{enumerate}
\usepackage{url}
\usepackage[stable]{footmisc}
\usepackage{breakurl}
\usepackage{hyperref}
\usepackage{comment}
\usepackage[style=ieee,dashed=false]{biblatex}
\usepackage[font=small,labelfont=bf]{caption}
\usepackage[font=small,labelfont=bf]{subcaption}
\usepackage{float}
\usepackage{mwe}
\usepackage{algorithm}
\usepackage{algpseudocode}
\usepackage{xcolor}
\usepackage[super]{nth}
\usepackage{soul}
\usepackage{centernot}
\DeclareGraphicsExtensions{.pdf}
\usepackage{amsmath}
\usepackage{amsbsy}
\usepackage{amssymb}
\usepackage{amscd}
\usepackage{amsfonts}

\newcommand{\R}{\mathbb R}

\newcommand{\beq}{\begin{equation}}
\newcommand{\eeq}{\end{equation}}
\newcommand{\beqs}{\begin{eqnarray}}
\newcommand{\eeqs}{\end{eqnarray}}
\newcommand{\beql}{\begin{equation} \label}
\newcommand{\half}{\frac{1}{2}}

\newcommand{\calB}{{\cal B}}

\newcommand{\calN}{{\cal N}}

\newcommand{\calS}{{\cal S}}

\newcommand{\calU}{{\cal U}}
\newcommand{\calV}{{\cal V}}

\newcommand{\parderiv}[2]{\frac{\partial #1}{\partial #2}}

\newcommand{\veps}{\varepsilon}

\newcommand{\p}{\partial}

\usepackage[margin=1in]{geometry}

\newcommand{\dee}{\mathcal{D}}
\newcommand{\scl}{\mathcal{L}}

\newcommand{\udkk}[1]{{\color{black} #1}}

\date{}
\begin{document}
\title{Inviscid Burgers as a degenerate elliptic problem}

\author{Uditnarayan Kouskiya\thanks{Department of Civil \& Environmental Engineering, Carnegie Mellon University, Pittsburgh, PA 15213, email: udk@andrew.cmu.edu.} $\qquad$ Amit Acharya\thanks{Department of Civil \& Environmental Engineering, and Center for Nonlinear Analysis, Carnegie Mellon University, Pittsburgh, PA 15213, email: acharyaamit@cmu.edu.}}

\maketitle
\begin{abstract}
\noindent We demonstrate the feasibility of a scheme to obtain approximate weak solutions to the (inviscid) Burgers equation in conservation and Hamilton-Jacobi form, treated as degenerate elliptic problems. We show different variants recover non-unique weak solutions as appropriate, and also specific constructive approaches to recover the corresponding entropy solutions.

\end{abstract}

\section{Introduction}
In this paper we continue, following \cite{KA1}, the assessment of a recently introduced duality based approach to solving differential equations involving evolution in time. A nonlinear ODE system was considered in \cite{KA1}, along with some linear PDE; here we consider the (inviscid) Burgers equation in conservation and Hamilton-Jacobi form. 

Our essential idea (see App.~\ref{app:rev_duality}) is to treat the primal PDE under consideration as constraints and invoke a more-or-less arbitrarily designable strictly convex, auxiliary potential to be optimized. Then, a dual variational principle for the Lagrange multiplier (dual) fields can be designed involving a dual-to-primal (DtP) mapping which is an adapted change of variables.  The dual variational principle has the special property that its Euler-Lagrange equations are exactly the primal PDE system, interpreted as equations for the dual fields using the DtP mapping, and this, even though the primal system may not have the required symmetries necessary to be the Euler-Lagrange equations of any objective functional of the primal fields alone. We use a simple Galerkin Finite element discretization of the dual Euler-Lagrange system (which is the primal system of interest using a change of variables) to look for a computational approximation of weak solutions to (inviscid) Burgers equation. The system is degenerate-elliptic in space-time domains. 

To our knowledge, our approach for generating approximate solutions to a nonlinear hyperbolic problem (albeit scalar here, but seamlessly generalizable to systems, formally at least)  is new. Brenier \cite[Sec.~4]{brenier2018initial} solves the inviscid Burgers problem without approximation through a maximization of a functional related to ours, but without the use of `base states,' (see Sec.~\ref{sec:Burgers_formulation} and App.~\ref{app:base_state_example}) which we find crucial in making our ideas work; in fact, we have to employ a sequence of dual functionals parametrized by evolving base states, the latter self-consistently prescribed by the scheme. As very nicely explained by a `safe mountain climbing' analogy, Brenier's exact scheme \cite[Sec.~4]{brenier2018initial} is different from our approach, and computational results are shown in \cite[pp.101-105]{brenier_book}; he also established connections to an Optimal Transport based method of attack, which can be utilized for computational approximations via the celebrated Benamou-Brenier formulation of optimal transport theory, as mentioned in \cite[Sec.~4.1, p.~597]{brenier2018initial}.

An outline of this paper is as follows: Sec.~\ref{sec:Burgers_formulation} and Sec.~\ref{sec:HJ_formulation}  comprise the development of the weak formulation for the dual inviscid Burgers equation in conservation and Hamilton-Jacobi forms, respectively. Sec.~\ref{sec:algorithm} describes the algorithm for the results computed in the paper. In Sec.~\ref{sec:results}, five selected problems are solved to illustrate and evaluate the features of the formulation and algorithm developed in  Sections \ref{sec:Burgers_formulation} and \ref{sec:HJ_formulation}, as they relate to generating approximate solutions to the (inviscid) Burgers equation. The paper also contains five appendices supporting various sections of the main narrative. A word on notation: we always use the Einstein summation convention for indices unless otherwise mentioned and except when the indices are $x,t$. \udkk{Our work is mathematically formal, and it is not within our terms of reference to speak definitively about the regularity classes of the various functions involved, particularly because we do not do existence proofs and because of the non-standard nature of our undertaking. Nevertheless, we have endeavored to state, to the best of our ability, the classes of functions to which we expect the various functions we employ to belong.}

\section{Dual formulation of the Burgers Equation}\label{sec:Burgers_formulation}
Burgers equation \cite{Burger,Bateman} (cf.~App.~\ref{app:new}) is a (convection-dominated) convection-diffusion partial differential equation given by
\begin{equation}\label{eq:viscous_Burgers}
\p_t u + u\,\p_x u = \nu \p_{xx} u\qquad \text{in } \Omega,
\end{equation}
where $\Omega = (0,L) \times (0,T) \subset \R \times \R^+$ is a domain in (1-d)space-time, $u: \Omega \to \R$ and $\nu \in \R^+$ are the flow velocity and viscosity, respectively, in the context of fluid mechanics.
\begin{figure}[h!]
    \centering
    \includegraphics[width=0.3\textwidth]{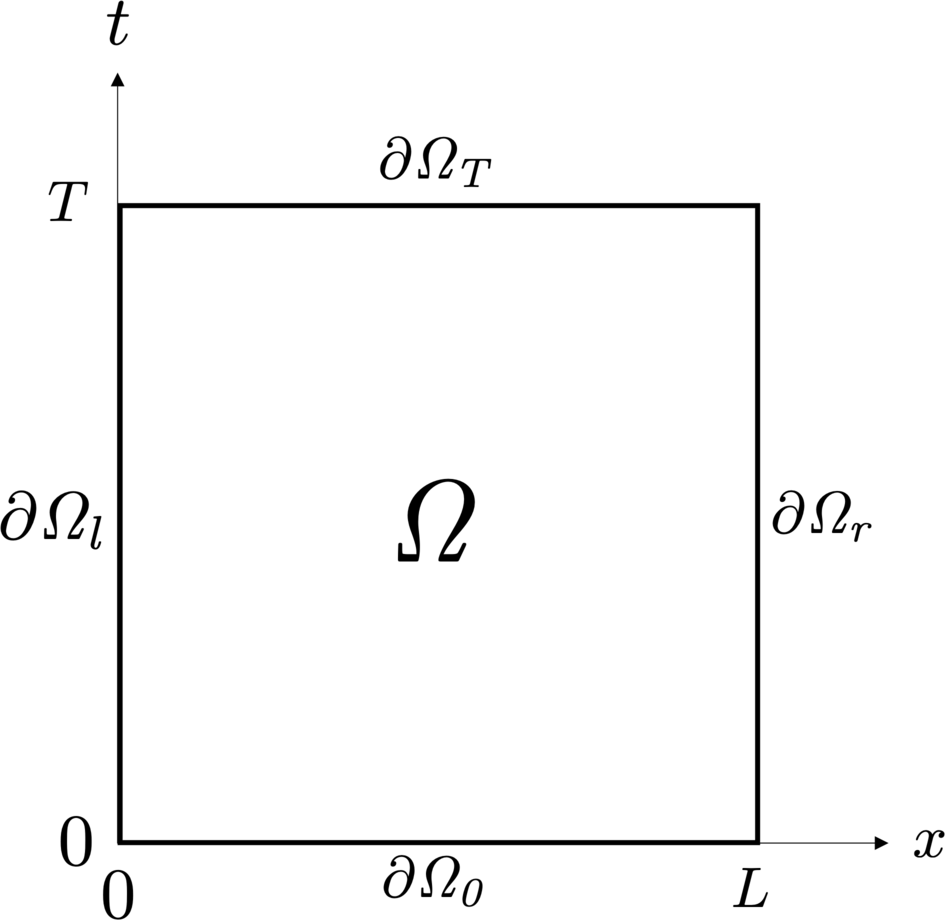}
    \caption{Domain of interest}
    \label{fig:domain}
\end{figure}
A schematic of the domain is shown in Fig.~\ref{fig:domain}. Henceforth, $L = 1$. The boundaries $\p \Omega_l$, $\p \Omega_r$, $\p \Omega_0$ and $\p \Omega_T$ will be referred to as the left, right, bottom, and top boundaries, respectively.

\udkk{In the inviscid case, $\nu = 0$, we think of solutions of Burgers equation in the class of bounded, piecewise-continuously-differentiable functions on $\Omega$ except possibly on a finite number of curves in $\Omega$. In the viscous case, we think of solutions $u$ belonging to the space $C^{0}(\Omega)$ (and possibly with higher regularity).}

We will be mainly interested in the inviscid Burgers equation, obtained by setting $\nu=0$ in \eqref{eq:viscous_Burgers}:
\begin{equation*}
    \p_t u + u\,\p_x u = 0 \qquad \text{in }\Omega,
\end{equation*} 
and, as is standard, explore weak solutions of this quasilinear equation expressed in conservation form:
\begin{equation}
\label{eq:Burgers_conservation_primal}
    \p_t u + \p_x \left( \half u^2  \right) = 0 \qquad \text{in }\Omega,
\end{equation}
with initial and the boundary conditions given by 
\udkk{
\begin{equation}
    u(x,0) = u_0(x); \qquad u(0,t) = \pm\, u_l(t).
    \label{eq:Burgers_primal_bc}
\end{equation}
where $u_0:(0,L) \to \R$ and $u_l:(0,T) \to \R$ are prescribed functions.}

Following the ideas in Appendix \ref{app:rev_duality}, we form the product of \eqref{eq:Burgers_conservation_primal} with a dual field $\lambda: \Omega \to \R$ (that acts like a Lagrange multiplier) and integrate by parts in the space-time domain to define a pre-dual functional
\begin{equation} \label{eq:Burgers_prefunc}
    \widehat{S}_H\,[u,\lambda]
    =\int_\Omega \,
    \left(-u \, \p_t \lambda  - \frac{u^2 \, \p_x \lambda}{2} + H(u,x,t)\right)\,dt\,dx - \int_0^L \,u_0(x)\,\lambda(x,0)\,dx
    - \int_0^T \,\frac{(u_l(t))^2}{2}\,\lambda(0,t)\,dt,
    \end{equation}
where $H$ is a `free' auxiliary potential of the arguments shown. \udkk{While not strictly necessary, we write many of the required expressions in terms of our general formalism to show the uniformity of our scheme over a range of problems. 

For our considerations related to the dual formulation of Burgers equation, the array $U = (u)$ and $\dee = (\lambda, \nabla \lambda) = (\lambda, \p_x \lambda, \p_t \lambda)$. We will index the primal array $U$ by lowercase Latin indices and the array $\dee$ by uppercase Greek indices. If the primal variable $U$ takes values in the set $\calU$, then $\hat{H}: \calU \times \calU \to \R$, and $H(U,x,t) = \hat{H}(U, \bar{U}(x,t))$, where both functions $\hat{H}$ and $\bar{U}:\Omega \to \calU$ are, in principle, open to design, and may be judiciously chosen.

Defining the `Lagrangian' identified from \eqref{eq:Burgers_prefunc}} as
$$ 
\mathcal{L}_H(U,\mathcal{D},x,t) := -u \, \p_t \lambda  - \frac{u^2 \, \p_x \lambda}{2} + H(u,x,t),$$
we can rewrite \eqref{eq:Burgers_prefunc} as follows:
\begin{equation*}
    \widehat{S}_H\,[U,D]= \int_\Omega \,\mathcal{L}_H(U,\mathcal{D},x,t)\, dt\, dx- \int_0^L \,u_0(x)\,\lambda(x,0)\,dx
    - \int_0^T \,\frac{(u_l(t))^2}{2}\,\lambda(0,t)\,dt.
\end{equation*}
\udkk{Here, we think of $\lambda$ to belong to the space of continuous functions on $\Omega$ with square-integrable distributional derivatives in $\Omega$. Additionally, $\lambda$ needs to have enough regularity such that space-time natural boundary conditions (e.g.~\eqref{eq:Burgers_dual_int_weak_form}) in which $\hat{u}$, defined in \eqref{eq:DtP_Burgers} below, appears need to make sense, as well as for the prescription of Dirichlet boundary conditions such as \eqref{eq:Burgers_dual_dbc}.}

The design of the functional involves the inclusion of only those space-time boundary terms that contain information available from the primal problem \cite{action_3, KA1}. Additionally, we impose the following Dirichlet space-time boundary conditions on $\lambda$:
\begin{subequations}\label{eq:Burgers_dual_dbc}
\begin{align}
        \lambda(L,t) &= \lambda_r(t); \label{eq:Burgers_dual_dbc1}\\ 
    \lambda(x,T) &= \lambda_T(x), \label{eq:Burgers_dual_dbc2}
\end{align}
\end{subequations}
\udkk{where $\lambda_r: (0,T) \to \R$ and $\lambda_T: (0,L) \to \R$ are arbitrarily specified continuous functions}. 
The choice of the auxiliary potential $H$ needs to be such that it renders the following equation `solvable' for $u$ (Appendix \ref{app:rev_duality}):
\udkk{\begin{equation}
    \frac{\p \mathcal{L}_H}{\p U} : -\p_t\lambda - u\,\p_x \lambda + \frac{\p}{\p U}H(U,x,t)=0.
    \label{eq:DtP-Burgers-eqn}
\end{equation}}
 We employ the following shifted quadratic form for $H(U,x,t)$:
\begin{equation}
\label{eq:Burgers_H_fun}
    H(U,x,t) = \frac{\beta_u}{2}\big(U-\bar{U}(x,t)\big)^2 = \frac{\beta_u}{2}\big(u-\bar{u}(x,t)\big)^2,
\end{equation}
 where $\beta_u \gg 1$. The function $\bar{U} = \bar{u}$ is referred to as a \textit{base state} (see \cite[Sec.~5]{action_3}-\cite[Sec.~4]{KA1}). The choice of this base state is often crucial for the success of the scheme (see Appendix \ref{app:base_state_example} for a simple example)  and their algorithmic use is shown in \cite[Sec.~5.3]{KA1} and in Section \ref{sec:results} of this work, among others. Solving for $U$ in \eqref{eq:DtP-Burgers-eqn}, one obtains the following \textit{dual-to-primal mapping} (DtP):
\begin{equation}\label{eq:DtP_Burgers}
\begin{aligned}
    U^{(H)}(\mathcal{D},x,t) & = \bar{U} + \frac{\bar{U}\,\p_x \lambda  + \p_t \lambda}{\beta_u-\p_x \lambda}\\
    \hat{u}(x,t) & := U^{(H)}\big(\dee(x,t), x,t\big).
\end{aligned}
\end{equation} 
The \textit{dual} functional is now defined as
\begin{equation}\label{eq:Burgers_functional}
\begin{aligned}
S_H[\lambda]:=\widehat{S}_H\left[U^{(H)}(\mathcal{D},x,t),\lambda\right] & = 
    \int_\Omega \,\mathcal{L}_H\left(U^{(H)}(\mathcal{D},x,t),\mathcal{D},x,t\right)\, dt\, dx \\
    &  - \int_0^L \,u_0(x)\,\lambda(x,0)\,dx
    - \int_0^T \,\frac{(u_l(t))^2}{2}\,\lambda(0,t)\,dt.
\end{aligned}
\end{equation}
The above functional can be explicitly written in terms of $\lambda$ as
\begin{multline}\label{eq:Burgers_functional_explicit}
    S_H[\lambda] = \int_\Omega \,\mathbb{K}\big|_{\dee}\,\left(\p_t\lambda + \bar{u}\, \p_x\lambda\right)^2\,dt\,dx - \int_\Omega \,\left(\bar{u} \, \p_t \lambda +\bar{u}^2\,\frac{\p_x \lambda }{2} \right)\,dt\,dx\\
    - \int_0^L \,u_0(x)\,\lambda(x,0)\,dx
    - \int_0^T \,\frac{(u_l(t))^2}{2}\,\lambda(0,t)\,dt,
\end{multline}
{
where 
$$\mathbb{K}\big|_{\dee}=\frac{-1}{2\beta_u\left(1 - \frac{\p_x \lambda}{\beta_u}\right)}$$
and the subscript $\dee$ represents the dependence of this coefficient on $\dee$. }Hence, for $\beta_u\gg \left|\p_x \lambda\right|$, the nonlinear part of the bulk integrand is non-positive.
\udkk{We now consider variations (test functions) $\delta \lambda: \Omega \to \R$ which are continuous on $\Omega$ and piecewise continuously differentiable, with $\delta \lambda = 0$ on $\p\Omega_T$ and $\p \Omega_r$ and define
\[
\delta^{(1)}S_H[\lambda;\delta\lambda] := \frac{d}{d\veps} S_H [\lambda + \veps \delta \lambda]\bigg|_{\veps = 0},
\]
which is the first variation of the above functional $S_H$ in a direction $\delta\lambda$. Similarly, we define $\delta \dee [\lambda; \delta \lambda] = \frac{d}{d\veps} \dee [\lambda + \veps \delta \lambda]\bigg|_{\veps = 0}$ (cf.~ App.~\ref{app:rev_duality}). Noting that in this case
\[
\parderiv{\scl_H}{U_i}
\parderiv{U^{(H)}_i}{\dee_\Gamma}
\delta \dee_\Gamma  +
\parderiv{\mathcal{L}_H}{\mathcal{D}_\Gamma}
\delta \mathcal{D}_\Gamma = \parderiv{\mathcal{L}_H}{\mathcal{D}_\Gamma}
\delta \mathcal{D}_\Gamma
\]
by \eqref{eq:DtP-Burgers-eqn} 
}, and recognizing that $\scl_H$ is necessarily affine in $\dee$ (with $H$ independent of $\dee$), one has
\udkk{
\begin{equation}
\begin{aligned}\label{eq:Burgers_dual_int_weak_form}
\delta^{(1)}S_H[\lambda;\delta\lambda]  &= \int_\Omega \left(-\hat{u} \, \p_t \delta \lambda  - \frac{\hat{u}^2 \, \p_x \delta \lambda}{2}\right) dt\, dx  \\ &- \int_0^L \,u_0(x)\, \delta \lambda(x,0)\,dx
    - \int_0^T \,\frac{(u_l(t))^2}{2}\, \delta \lambda(0,t)\,dt.  
    \end{aligned}
    \end{equation}
    }
Using $\delta \lambda$ consistent with \eqref{eq:Burgers_dual_dbc}, the E-L equations and side-conditions are given by 
\begin{subequations}
\begin{gather}
\p_t \hat{u} + \frac{\p}{\p x}\left(\frac{\hat{u}^2}{2}\right) = 0 \qquad \text{in }\Omega;
\label{eq:Burgers_primal_in_dual_1}\\
\hat{u}(x,0) = u_0(x); \qquad \hat{u}(0,t) = \pm\,u_l(t), \label{eq:Burgers_primal_in_dual_2}
\end{gather} 
\label{eq:Burgers_primal_in_dual}%
\end{subequations}
which are simply \eqref{eq:Burgers_conservation_primal}-\eqref{eq:Burgers_primal_bc} with the replacement $u \to \hat{u}$. \udkk{Equations \eqref{eq:Burgers_dual_int_weak_form} results under the assumption that $\hat{u}$ is $C^1(\Omega)$. When this is not true, e.g.~along a finite number of curves in $\Omega$, additional contributions arise along such curves, as in App.~\ref{sec:weak_form_implications}. 
}

Thus, solutions of the problem defined by \eqref{eq:Burgers_conservation_primal}, \eqref{eq:Burgers_primal_bc} can be generated by solving \eqref{eq:Burgers_primal_in_dual} along with \eqref{eq:Burgers_dual_dbc}, where the DtP mapping \eqref{eq:DtP_Burgers} is utilized to bridge the primal variable $u$ with the dual variable $\lambda$. In short, our scheme may be interpreted as designing an adapated \textit{change of variables} for solving Burgers equation.

In \cite{dual_cont_mech_plas}, the degenerate ellipticity of the dual formulation of a significant class of equations from continuum mechanics is analyzed. We apply those ideas next to the dual formulation of Burgers equation given by \eqref{eq:Burgers_primal_in_dual} written as
\begin{equation*}
    \p_t \big(\mathcal{F}_1(\hat{u})\big) + \p_x \big(\mathcal{F}_2(\hat{u})\big) = 0,
\end{equation*}
some of whose analytical properties are governed by the terms
\begin{equation*}
    A_{ij} = \frac{\p \,\mathcal{F}_i}{\p (\nabla D)_j} = \left. \frac{\p \,\mathcal{F}_i}{\p u} \right|_{\hat{u}}\, \frac{\p \,u}{\p (\nabla D)_j}.
\end{equation*}
\udkk{ The expressions for $A_{ij}$ are given by}
\begin{equation*}
    \begin{gathered}
        A_{11} = \frac{1}{\beta_u - \p_x \lambda}; \quad 
        A_{12} = \frac{\hat{u}}{\beta_u - \p_x \lambda}; \quad
        A_{21} = \frac{\hat{u}}{\beta_u - \p_x \lambda}; \quad
        A_{22} = \frac{\hat{u}^2}{\beta_u - \p_x \lambda}.
    \end{gathered}
\end{equation*}
Consequently, for $\nabla D = (0,0)$
\begin{equation*}
    c_i \, A_{ij} \, c_j = \frac{(c_1 + \bar{u} \,c_2)^2}{\beta_u} \geq 0 \qquad \forall (c_1,c_2) \in \mathbb{R}^2
\end{equation*}
which establishes that $A_{ij}|_{\nabla D =0}$ is positive semi-definite. \udkk{Let $\mathcal{N}$ represent a neighborhood around $\mathcal{D} = 0$ given by
$\mathcal{N} = \left\{a \in \R^3: \left \lVert a \right \rVert_3 < \beta_u \right\},$
(which is conservative). Then, for $\mathcal{D} \in\mathcal{N}$ the following expression holds:}
\begin{equation*}
c_i \, A_{ij} \, c_j = \frac{(c_1 + \hat{u} \,c_2)^2}{\beta_u - \partial_x \lambda} = \frac{(c_1 + \hat{u} \,c_2)^2}{\beta_u\left(1-\frac{\p_x \lambda}{\beta_u}\right)} \geq 0\qquad \forall (c_1,c_2) \in \mathbb{R}^2.
\end{equation*}
This relation establishes the positive semi-definiteness of $A$ within this specific neighborhood. Consequently, the equation \eqref{eq:Burgers_primal_in_dual} is locally degenerate elliptic and for any surface with normal $(n_1,n_2)$ in the space-time domain ($n_1$ and $n_2$ represents the projection of normal in the time and space directions, respectively), ellipticity fails when
$$(n_1,n_2) = \kappa(-\hat{u},1) \quad \text{for } \kappa\in\mathbb{R},$$
where $\hat{u} $ is given by \eqref{eq:DtP_Burgers}.

\subsection{Weak formulation of the dual Burgers equation}
A weak formulation of \eqref{eq:Burgers_primal_in_dual} with the side conditions \eqref{eq:Burgers_dual_dbc} can be generated in the usual way through integration by parts or by considering the first variation of the functional \eqref{eq:Burgers_functional} or \eqref{eq:Burgers_functional_explicit}  and setting it equal to zero. For any $\delta \lambda$ \udkk{ in the class discussed in Sec.~\ref{sec:Burgers_formulation}} satisfying the conditions stated below, we intend to find the dual field $\lambda$ which satisfies the following equations:
\begin{equation}
\begin{gathered}
 R[\lambda;\delta \lambda]:=\int_\Omega \,
    \left(-\hat{u} \, \p_t \delta\lambda - \frac{\hat{u}^2 \,\p_x\delta\lambda}{2}  \right)\,dt\,dx - \int_0^L \,u_0(x)\,\lambda(x,0)\,dx
    - \int_0^T \,\frac{(u_l(t))^2}{2}\,\lambda(0,t)\,dt=0; \\
     \delta \lambda(x,T) = 0; \qquad \delta \lambda(L,t) = 0; \\
      \udkk{ \lambda(x,T) = \lambda_T(x);} \qquad
      \lambda(L,t) = \lambda_r(t),
\end{gathered} \label{eq:weak_dual_Burgers}
\end{equation}
\udkk{where the functions $\lambda_T(\cdot)$ and $\lambda_r(\cdot)$ are discussed in Sec.~\ref{sec:Burgers_formulation}} satisfying $\lambda_T(L) = \lambda_r(T)$.
The dual scheme formally guarantees that the solution to \eqref{eq:weak_dual_Burgers} implies the solution to the set of equations \eqref{eq:Burgers_primal_in_dual}. We subsequently make use of the above weak form to numerically compute an approximate solution for the dual field and utilize the DtP mapping \eqref{eq:DtP_Burgers} to obtain the corresponding field for the primal problem i.e.~the inviscid Burgers equation. 

\section{Dual formulation of the (inviscid) Burgers Equation in Hamilton-Jacobi form}\label{sec:HJ_formulation}

Here we consider the Burgers equation \eqref{eq:Burgers_conservation_primal} in Hamilton-Jacobi form:
\begin{equation}\label{eq:HJ_original}
    \p_t Y + \frac{(\p_x Y)^2}{2} = \nu \, \p_{xx} Y \qquad \text{in }\Omega.
\end{equation}
Differentiating \eqref{eq:HJ_original} w.r.t. $x$, and defining $\p_x Y =: u$ gives \eqref{eq:viscous_Burgers}. Throughout the following text, we will refer to this form of the Burgers equation as \textit{Burgers-HJ} and when $\nu=0$, we will refer to it as the \textit{inviscid Burgers-HJ} equation. 
 As already stated, our primary interest is in the inviscid case $\nu=0$. 
We will consider the $\nu\neq 0$ case to shed light on the inviscid case, the motivation for which will be provided in the relevant examples. 

It is clear that \eqref{eq:HJ_original} requires $Y$ to be specified at one point of the domain, above and beyond information available from the corresponding Burgers equation. If the function $u$ satisfies Burgers equation and $\p_x Y = u$ then
$\p_t Y (x,t) + \half u^2 (x,t) - \nu \p_x u(x,t) = f(x^*,t)$ for an arbitrarily fixed $x^*$, where $f(x^*,t) = \p_t Y (x^*,t) + \half u^2 (x^*,t) - \nu \p_x u(x^*,t)$. In all problems considered in the text and Appendix, we will always work with $f(x^*,t) = 0$ with knowledge of the Burgers solution of interest known at the point $x^*$ (mostly a boundary point) from which $Y(x^*,t)$ is determined, with an initial condition $Y(x^*,0)$ arbitrarily specified.
The dual formulation of \eqref{eq:HJ_original} for $\nu=0$ is obtained by writing it in first-order form \cite[Sec.~6.2]{action_2}:
\begin{subequations}
\begin{gather}
    \p_t Y = -\frac{u^2}{2}  \quad \text{in }\Omega; \label{eq:HJ_primal_1}\\
    \p_x Y = u \qquad\text{in }\Omega\label{eq:HJ_primal_2}.
\end{gather}    \label{eq:HJ_primal}%
\end{subequations}
\udkk{Here, we think of $Y \in C^0(\Omega)$.} The initial and the boundary conditions can be given as
\begin{equation}
    Y(x,0) = Y_0(x); \qquad Y(0,t) = Y_l(t).
    \label{eq:HJ_primal_bc}
\end{equation}
\udkk{where $Y_0:(0,L) \to \R$ and $Y_l:(0,T) \to \R$ are prescribed functions.} Corresponding to the equations \eqref{eq:HJ_primal_1} and \eqref{eq:HJ_primal_2} we introduce the dual fields $\lambda$ and $\gamma$, respectively, \udkk{both with regularity as specified for $\lambda$ for the discussion of Burgers in Sec.~\ref{sec:Burgers_formulation}.} Similar to the setup explained in Sec.~\ref{sec:Burgers_formulation}, we employ a shifted quadratic form for the auxiliary potential $H$ in both $Y$ and $u$:
\begin{equation}
\label{eq:HJ_H_fun}
    H(Y,u,x,t) = \frac{\beta_Y}{2}(Y-\bar{Y})^2 + \frac{\beta_u}{2}(u-\bar{u})^2,
\end{equation}
where $\bar{Y}$ and $\bar{u}$ correspond to the base states for $Y$ and $u$, respectively. With this choice of $H$, the DtP mapping for $Y$ and $u$, given by $Y^{(H)}$ and $u^{(H)}$, can be expressed in terms of the dual fields and their derivatives as
\begin{subequations}
\begin{gather}
Y^{(H)}(\mathcal{D},x,t) = \bar{Y}+\frac{\p_t \lambda + \p_x \gamma}{\beta_Y};\label{eq:DtP_HJ1}
\\
u^{(H)}(\mathcal{D},x,t) = \bar{u}+\frac{\gamma - \lambda \bar{u}}{\beta_u + \lambda},
\label{eq:DtP_HJ2}
\end{gather}
\label{eq:DtP_HJ}%
\end{subequations}
\udkk{where $\mathcal{D}$ for this Burgers H-J problem, following the general formalism, is given by the array $\mathcal{D}= (D, \nabla D) = (\lambda,\gamma,\p_t \lambda, \p_x \lambda, \p_t \gamma, \p_x \gamma)$ (with the array $U = (u, Y)$ and $D = (\lambda, \gamma)$).}
Additionally, the following boundary conditions are imposed:
\begin{subequations}
\begin{align}
\gamma(L,t) &= \gamma_r(t);\label{eq:HJ_dual_BC1} \\ \lambda(x,T) &= \lambda_T(x),\label{eq:HJ_dual_BC2}
\end{align}\label{eq:HJ_dual_BC}%
\end{subequations}
\udkk{where $\gamma_r: (0,T) \to \R$ and $\lambda_T: (0,L) \to \R$ are arbitrarily specified continuous functions.} The \textit{dual functional} obtained using our scheme is
\begin{multline}
S_H[\lambda,\gamma] 
    =
    \bigintsss_\Omega \,
    \Biggl(-\hat{Y}\,\p_t \lambda \,  + \frac{\hat{u}\,^2 \lambda}{2} 
    -\hat{Y}\,\p_x\gamma  - \hat{u}\,\gamma +H\left(\hat{Y},\hat{u},x,t\right)\Bigg)\,dt\,dx\\ - \int_0^L \,Y_0(x)\,\lambda(x,0)\,dx
    - \int_0^T \,Y_l(t)\,\gamma(0,t)\,dt , \label{eq:HJ_functional}    
\end{multline}
where
\[
\hat{Y}(x,t) := Y^{(H)}\big(\mathcal{D}(x,t),x,t\big); \qquad\hat{u}(x,t) := u^{(H)}\big(\mathcal{D}(x,t),x,t\big).
\] 
The functional above can be written explicitly in terms of the dual variables as
\begin{multline} \label{eq:HJ_functional_explicit}S_H[\lambda,\gamma] = \int_\Omega \,\left( -\frac{\mathbb{K}_1}{\beta_Y}\,(\p_t \lambda \,+\, \p_x \gamma)^2 \, -  \frac{\mathbb{K}_2}{\beta_u}(\bar{u}\lambda - \gamma)^2 \right) \,dt\,dx\\
+ \int_\Omega \,\left(-\bar{Y}(\p_t \lambda + \p_x \gamma) + \bar{u}\left(\frac{\lambda\bar{u}}{2} - \gamma\right)\right)\,dt\,dx \\- \int_0^L \,Y_0(x)\,\lambda(x,0)\,dx
    - \int_0^T \,Y_l(t)\,\gamma(0,t)\,dt,
\end{multline}
where 
$$\mathbb{K}_1 = \frac{1}{2}; \qquad \mathbb{K}_2 = \frac{1}{2\left(1+\frac{\lambda}{\beta_u}\right)},$$
and for $\beta_u\gg |\lambda|$, $\mathbb{K}_2 > 0$. Hence, the integrand containing nonlinear terms in the dual fields is negative semi-definite for this range of $|\lambda|$.
The Euler-Lagrange equations of the dual functional \eqref{eq:HJ_functional}, when extracted using the conditions \eqref{eq:HJ_dual_BC}, recover the primal equations \eqref{eq:HJ_primal} and \eqref{eq:HJ_primal_bc}, expressed in terms of their dual counterparts: 
\begin{subequations}
\begin{equation}
\begin{gathered}
    \p_t\hat{Y} = -\frac{\bigl(\hat{u}\bigl)^2}{2}  \quad \text{in }\Omega; \\
    \p_x\hat{Y} = \hat{u} \quad \text{in 
 } \Omega;
 \end{gathered}
 \end{equation}
 \begin{equation}
    \hat{Y}(x,0) = Y_0(x) \qquad \hat{Y}(0,t) = Y_l(t).
    \end{equation}\label{eq:HJ_primal_in_dual}%
\end{subequations}
Thus, a solution to the the above set of equations in dual variables implies a solution the primal problem using the DtP map \eqref{eq:DtP_HJ}. \udkk{As in the case of Burgers equation, the above equations are derived based on the assumption that $\hat{Y}$ and $\hat{u}$ are continuously differentiable. When this is not true, e.g.~along a finite number of curves in $\Omega$, additional contributions arise along such curves,} see \ref{sec:weak_form_implications}.

To examine the ellipticity of the dual equations \eqref{eq:HJ_primal_in_dual}, we rewrite them as
\begin{equation*}
    \begin{gathered}
        \p_t \big(\mathcal{F}_{11}(\hat{Y},\hat{u})\big) + \p_x \big(\mathcal{F}_{12}(\hat{Y},\hat{u})\big) + \mathcal{G}_1\big(\hat{Y},\hat{u}\big) = 0; \\
        \p_t \big(\mathcal{F}_{21}(\hat{Y},\hat{u})\big) + \p_x \big(\mathcal{F}_{22}(\hat{Y},\hat{u})\big) + \mathcal{G}_2\big(\hat{Y},\hat{u}\big) = 0,
    \end{gathered}
\end{equation*}
where $\mathcal{F}_{12} = \mathcal{F}_{21} = 0$. The behavior  of the above set of equations is now governed  by matrix
\begin{equation*}
\begin{gathered}
    \mathbb{A}_{ijk\ell} = \frac{\p \mathcal{F}_{ij}}{\p (\nabla D)_{k\ell}}; \qquad \nabla D= 
\begin{bmatrix}
    \p_t \lambda & \p_x \lambda \\
    \p_t \gamma & \p_x \gamma \\
\end{bmatrix},  
\end{gathered}
\end{equation*}
\udkk{and each of the indices $i,j,k,\ell$ ranging from $1$ to $2$}. Also, $\mathbb{A}_{ijk\ell}=0$ except 

\begin{equation*}
    \mathbb{A}_{1111} = \mathbb{A}_{1122} = \mathbb{A}_{2211} = \mathbb{A}_{2222} = \frac{1}{\beta_Y}.
\end{equation*}
Equivalently, we represent $\mathbb{A}$ as a $2 \times 2$ matrix $A$:
\begin{equation}\label{eq:HJ_ellipticity_matrix}
    A = 
\begin{bmatrix}
    \mathbb{A}_{1111} & \mathbb{A}_{1112} & 
    \mathbb{A}_{1121} & 
    \mathbb{A}_{1122}  \\
    \mathbb{A}_{1211} & \mathbb{A}_{1212} & 
    \mathbb{A}_{1221} & 
    \mathbb{A}_{1222}  \\
    \mathbb{A}_{2111} & \mathbb{A}_{2112} & 
    \mathbb{A}_{2121} & 
    \mathbb{A}_{2122}  \\
    \mathbb{A}_{2211} & \mathbb{A}_{2212} & 
    \mathbb{A}_{2221} & 
    \mathbb{A}_{2222}  \\
\end{bmatrix}
\end{equation}
which then satisfies
\begin{equation*}
    c_i A_{ij} c_j = \frac{(c_1 + c_4)^2}{\beta_Y}\geq 0 \qquad \forall c\in \mathbb{R}^4 ,
\end{equation*}
 establishing that $\mathbb{A}_{ijk\ell}$ is positive semi-definite on the space of $2 \times 2$ matrices and 
$$C_{ij} \mathbb{A}_{ijk\ell} C_{k\ell} = 0 \quad \mbox{for} \quad \left\{ C \in \mathbb{R}^{2\times2} \big| \ C_{11} + C_{22} = 0 \right\}.$$
Hence the equation set \eqref{eq:HJ_primal_in_dual} is degenerate elliptic. Furthermore, for any rank-one matrix $C$ of the form $a\otimes n$ satisfying $C : \mathbb{A} C = 0$, for some $a\in\mathbb{R}^2$ and some $n \in \R^2, |n| = 1$, the latter representing a unit normal in the space-time domain,  the degenerate ellipticity (positive semi-definiteness of $\mathbb{A}$) implies a loss of ellipticity along the direction $n$.
 Clearly, for any arbitrarily fixed unit normal $n$, there always exists a vector $a$ such that $a_1 n_1 = -a_2 n_2$. As a result, ellipticity, i.e.~$\det(\mathbb{A}_{ijk\ell} n_j n_{\ell}) >0$ fails along every unit normal in the domain (which is also apparent from the direct calculation  $\mathbb{A}_{ijk\ell} n_j n_\ell = \frac{1}{\beta_Y} n_i n_j$ which is a rank-one matrix). 
 

\subsection{Weak formulation for the dual Hamilton-Jacobi form of inviscid Burgers}
We construct a weak form for the system of equations \eqref{eq:HJ_primal_in_dual} with the side condition \eqref{eq:HJ_dual_BC} in the usual way through integration by parts. As in the case of Burgers equation, the weak form can also be obtained from the first variation of the dual functional \eqref{eq:HJ_functional} or \eqref{eq:HJ_functional_explicit}. For any field $\delta \lambda$ and $\delta \gamma$ \udkk{in the class discussed in Sec.~\ref{sec:HJ_formulation} satisfying the conditions stated below, we intend to find the dual fields $\lambda$ and $\gamma$ which satisfy the following equations:}
\begin{equation}
\begin{gathered}
\begin{aligned}\label{eq:weak_dual_HJ}
   R\,[\lambda, \gamma;\delta \lambda, \delta \gamma]:= \int_\Omega \,
    \bigg( -\hat{Y}\,\p_t \delta \lambda \,  + &\frac{\hat{u}\,^2 \,\delta \lambda}{2} 
    -\hat{Y}\,\p_x\delta \gamma  - \hat{u}\,\delta \gamma \bigg) \,dt\,dx\\ & - \int_0^L \,Y_0(x)\,\big(\delta \lambda(x,0)\big)\,dx
     - \int_0^T \,Y_l(t)\,\big(\delta \gamma(0,t)\big)\,dt=0; \end{aligned} \\
     \delta \lambda(x,T) = 0; \qquad \delta \gamma(L,t) = 0 ; \\
       \udkk{\lambda(x,T) = \lambda_T(x);} \qquad
      \gamma(L,t) = \gamma_r(t),
\end{gathered}  
\end{equation}
\udkk{where the functions $\lambda_T(x)$ and $\gamma_r(t)$ are discussed in Sec.~\ref{sec:HJ_formulation}}.



\section{Algorithm}
\label{sec:algorithm}
As introduced in \cite{KA1}, we compute approximate solutions to the problems of interest by solving them in a time-concatenated series of space-time subdomains whose (closed) union forms (the closure of) the entire domain $\Omega$ of interest. We solve a distinct, dual (space-time) bvp in each of these stages, step-by-step, marching in time, where the initial condition and base state in any stage depends on the output produced from the previous stage, details of which will be explained in the description of example problems to follow. 
Each of these subdomains, $\Omega^{(s)}$, is referred to as a \textit{stage}, indexed by $s = 1,2,3,\ldots, N$, $\overline{\Omega} = \cup_{s=1}^n \overline{\Omega^{(s)}}$.

\begin{figure}
\centering
\begin{subfigure}{.29\textwidth}
  \centering
  \includegraphics[width=.93\linewidth]{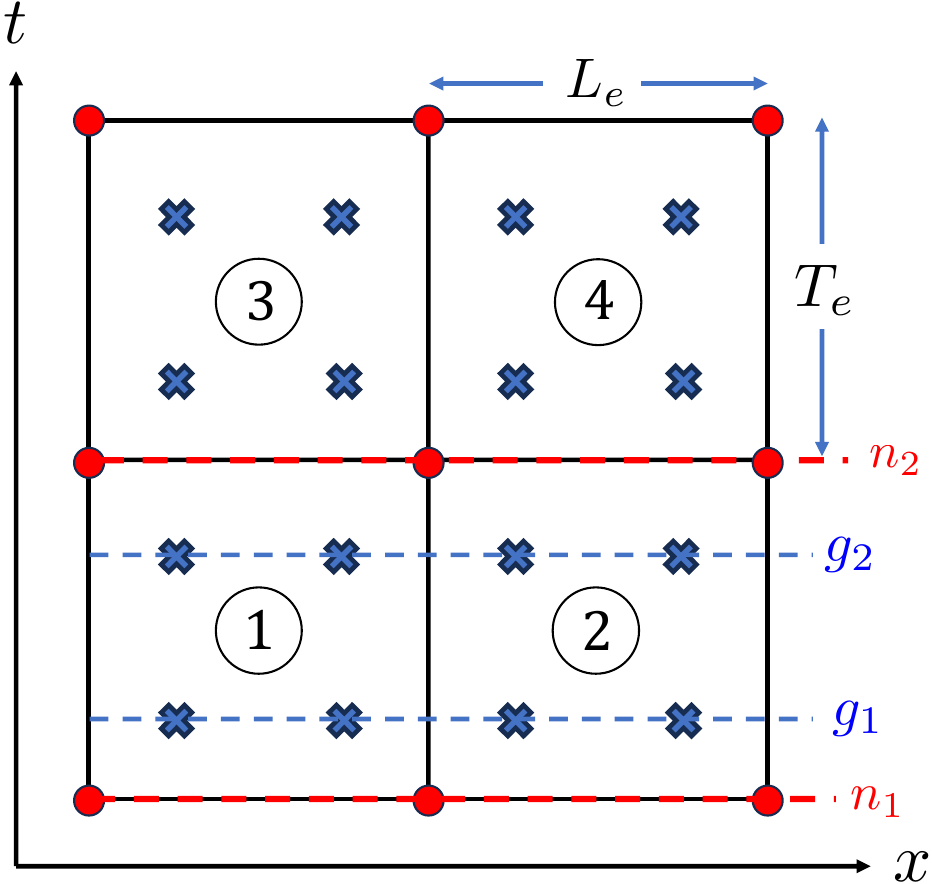}  
  \caption{FEM element geometry}
  \label{fig:Gauss_times}
\end{subfigure}
\begin{subfigure}{.29\textwidth}
  \centering
  \includegraphics[width=.85\linewidth]{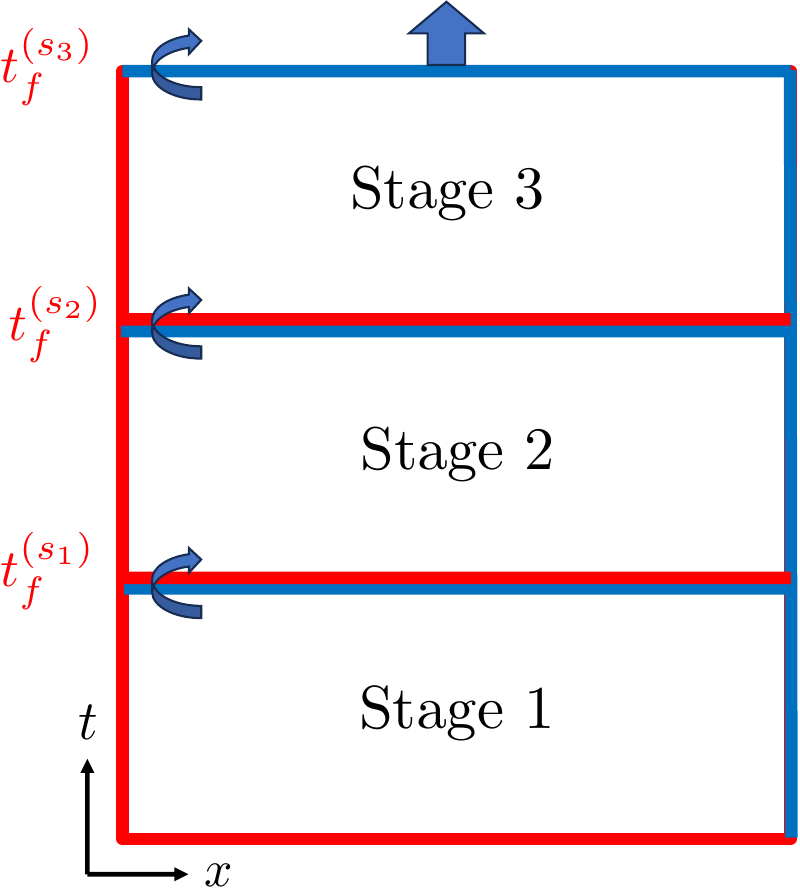}  
  \caption{Time marching}
  \label{fig:time_march}
\end{subfigure}
\begin{subfigure}{.4\textwidth}
  \centering
  \includegraphics[width=.8\linewidth]{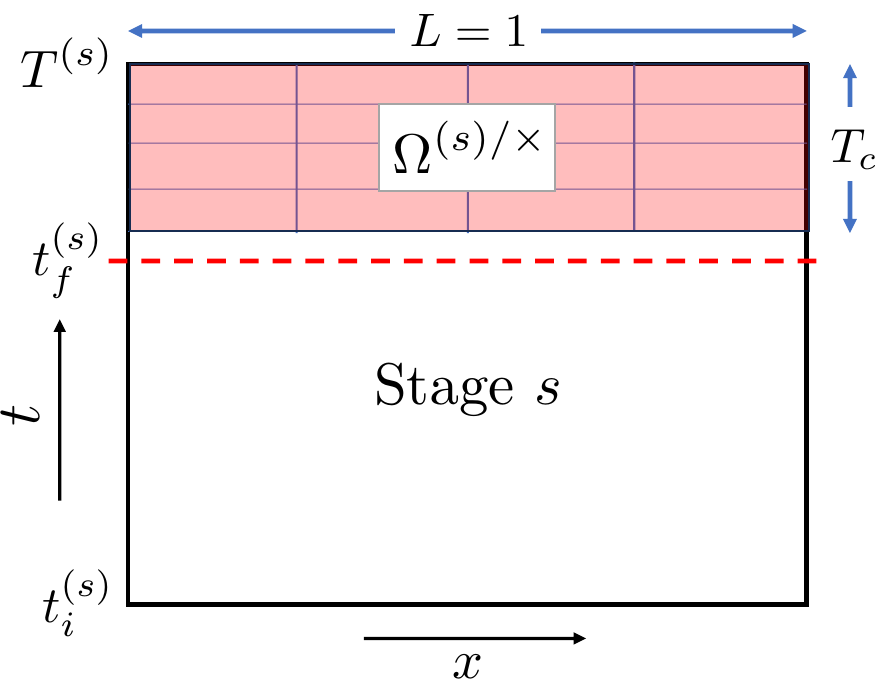}  
  \caption{$\Omega^{(s)}$}
  \label{fig:stage_2}
\end{subfigure}
\caption{(a) A sample $2\times2$ mesh with element numbers at center. Red dashed lines, passing through the nodal points (red dots), represent nodal timelines. $g_{1}$ and $g_{2}$ represent the first and the second Gauss
    timelines passing through the Gauss points of the elements $1$ and $2$. (b) represents the concatenated domains after truncation. Additionally, $t_f^{(s)} = t_i^{(s+1)}$.}
\label{fig:stage}
\end{figure}

A two-point Gauss quadrature scheme, in each of the $x$ and $t$ directions, has been utilized to approximate all the domain and boundary integrals appearing in the work. We will refer to the collection of nodes corresponding to a specific time (falling on the nodes of discretization) as a \textit{nodal timeline}. Similarly, we will refer to the collection of all Gauss points corresponding to a particular time as a \textit{Gauss timeline}. The setup for a sample mesh of $2\times2$ has been shown in Fig.~\ref{fig:Gauss_times}.

Due to the nature of the arbitrarily chosen (without loss of generality) $0$ b.c.~utilized on the dual fields (\eqref{eq:Burgers_dual_dbc}  and \eqref{eq:HJ_dual_BC}) at the final-time boundary of a stage,  strong gradients/boundary-layers can arise in the dual solution that can require excessive refinement to resolve.
Since the time-like extent of any stage is  an arbitrary choice subject only to computational expediency and accuracy, in every stage, solutions obtained for a specified layer of elements in the time-like direction near the final time boundary of the stage are ignored. We refer to this operation as `truncation'. The following stage initiates from a suitable time near the end of the current stage on the retained mesh (either at a Gauss timeline or at a nodal timeline, see Fig.~\ref{fig:Gauss_times}), referred to as the `cutoff' line. 
A visual representation illustrating this concept is presented in the Fig.~\ref{fig:stage_2}.

We employ the Newton Raphson method (NR) within each stage to approximate the solutions. A linear span of globally continuous, piecewise smooth finite element shape functions corresponding to an FE mesh for $\Omega^{(s)}$ is used to achieve this discretization.
These shape functions are represented by $N^{(\cdot)}$, where $(\cdot)$ denotes the index of any node on the space-time mesh. 
The discretized dual fields and test functions are expressed as
\begin{equation}
\label{eq:discrete_duals}
D_i(x,t) = D^A_i\, N^A(x,t); \qquad   \delta D_i(x,t) = \delta D^A_i\, N^A(x,t); \qquad dD_i(x,t) = dD^A_i \,N^A(x,t),
\end{equation}
where $D^A_i$ denotes the finite element nodal degrees of freedom. The discretized version of the appropriate weak form  \eqref{eq:weak_dual_Burgers} or \eqref{eq:weak_dual_HJ} generates a discrete residual $R^A_i(\dee)$ (derived in Sec.~\ref{sec:Burgers_examples} for the Burgers equation and Sec.~\ref{sec:HJ_examples} for the Burgers-HJ)  given by 
\begin{equation*}
    R[D; \delta D] = \delta D^A_i \, R^A_i(\dee)
\end{equation*}
and its variation in the direction $dD$, the Jacobian
\begin{equation*}
    J[D; \delta D, d D] = \delta D^A_i \, J^{AB}_{ij}(\dee) \, dD^B_j.
\end{equation*}
Any stage begins with an initial guess of $D_i^{A(0)} = 0$ followed by solving the following matrix equation for the $k^{th}$ correction:
\begin{equation}
    \begin{aligned}
     - R^A_i\left(\dee^{(k-1)}\right) & =  J^{AB}_{ij}\left(\dee^{(k-1)}\right)\, dD^B_j\\
         D^{B(k)}_j & = D^{B(k-1)}_j + dD^B_j,
    \end{aligned}
    \label{eq:update_duals}
\end{equation}
where $\mathcal{D} = (\p_t \lambda,\p_x\lambda)$ for the Burgers equation and $\mathcal{D}=(\lambda,\gamma,\p_t \lambda, \p_x \gamma)$ for Burgers-HJ. The convergence criteria for NR is set as follows:  $$\underset{A,i}{\max} \,\left|R^A_i\right| < tol,$$
where $tol$ represents a user-defined threshold tolerance. 

The following notation is used in the algorithm:
\begin{center}
 \renewcommand{\arraystretch}{1.3} 
  \begin{tabular}{|c|c|}
    \hline
    $(\cdot)^{(s)(k)}$ & value of $(\cdot)$ at $k^{th}$ NR iterate for stage $s$\\
    $N_x$ & number of elements in $x$ direction\\
    $N_t^{(s)}$ & number of elements in $t$ direction per stage\\ 
    $t_i^{(s)}$ & time at the start of stage $s$\\
    $T^{(s)}$ & time over which the stage $s$ is solved\\
    $t_f^{(s)}$ & time at the cutoff line for stage $s$\\
    $\Omega^{(s)}$ & $(0,L)\times\left(t_i^{(s)},t_i^{(s)} + T^{(s)} \right)$; domain for stage $s$\\
    $\bar{U}_k^{(s)}$ & base state for $U_k$ in stage $s$\\
    $\beta_k$ & coefficient(s) defining $H$ in \eqref{eq:Burgers_H_fun} or \eqref{eq:HJ_H_fun}\\
    $U_0^{(s)}$ & primal initial condition for stage $s$\\
    $\Omega^{(s) \times}$ & discarded domain for the stage $s$\\
    $N_c \times N_x $ & number of discarded elements per stage\\
    \hline
  \end{tabular}
\end{center}

\begin{table}[!ht]
{\begin{algorithm}[H]
\caption*{\textbf{Algorithm}} 
\begin{algorithmic}
\State \textbf{Initialization}:
\begin{enumerate}
\item Set $s=1$, $t_f^{(0)} = 0$ and $U_0^{(1)}(x)$ from \eqref{eq:Burgers_primal_bc} or \eqref{eq:HJ_primal_bc} depending on the problem under consideration.
\item Choose the values for $T^{(s)}$, $N_x$, $N_t^{(s)}$, $\beta_k$, $N_c$ and $tol$.
\end{enumerate}

\\\hrulefill
\State \textbf{$\boldsymbol{s_{th}}$ stage}: 
\begin{enumerate}
    \item \label{start-stage} Set $t_i^{(s)} = t_f^{(s-1)}$ and  generate $U_0^{(s)}$ ($s>1$, see Sec.~\ref{sec:algo1} and Sec.~\ref{sec:algo2} for details). Over the domain $\Omega^{(s)}$, set appropriate $\bar{U}_k^{(s)}$.
    \item Initiate NR: Set $D_k^{A(s)(0)} = 0$.
    \item[] \textbf{For} $j \geq 0$:
    \begin{enumerate}[i]
        \item \label{step-i} Evaluate $R_i^{A(s)(j)}$.
        \item Set $d^{(s)(j)} = \underset{A,i}{\max} \,\left|R^{A(s)(j)}_i\right|.$ 
        \item[]\textbf{if} $d^{(s)(j)}<tol$ \textbf{then} go to step \ref{stepout} (exit loop).
        \item[]\textbf{else} continue. 
        \item Evaluate $J_{ij}^{AB(s)(j)}$.
        \item Evaluate $D_k^{A(s)(j+1)}$ using \eqref{eq:update_duals}.  
        \item \textbf{do } $j=j+1$ and \textbf{go to} step \ref{step-i} 
    \end{enumerate}
    \item \label{stepout} $D_k^{(s)(j)}$ serves as the dual solution for the current stage, while the corresponding $U_k$ serves as the solution for the primal problem.
    Discard the results obtained in the region $\Omega^{(s)\times}$.
    \item Set an appropriate $t_f^{(s)}$ near the final time obtained on $\Omega^{(s)}\backslash\Omega^{(s)\times}$.
    \item Set $s=s+1$ and repeat steps 1-4 until $t_f^{(s)} \geq T$
\end{enumerate}
\end{algorithmic}
\end{algorithm}}
\caption{Algorithm to solve Burgers equation or Burgers-HJ. The subscript index $k$ takes the value of 1 for Burgers and ranges over $(1,2)$ for Burgers-HJ.}
\label{algo:euler_algorithm}
\end{table}

The algorithm has been presented in table \ref{algo:euler_algorithm}. Additional details for the algorithm are as follows:
\begin{itemize}
    \item Corresponding to the discretized dual fields $D$, the primal field $U$ can be evaluated  using the appropriate DtP mapping \eqref{eq:DtP_Burgers} or \eqref{eq:DtP_HJ} on the Gauss points within any element, which is further used in the evaluation of $R_i^A$ and $J_{ij}^{AB}$.
    \item For each stage, we discard the results over the final $T_c$ duration of time in that stage (see Fig.~\ref{fig:stage_2}). Denoting by $T_e$ the edge-length of elements in the time-like direction (see Fig.~\ref{fig:Gauss_times}), the truncation operation is carried out over $N_c := \frac{T_c}{T_e}$ space-like layers of elements.
    
    \item For the Burgers equation, $U_0(x):= u_0(x)$, whereas for the Burgers-HJ, $U_0(x):= Y_0(x)$. For the first stage, initial condition is used from the problem definition \eqref{eq:Burgers_primal_bc}or \eqref{eq:HJ_primal_bc}, and for all subsequent stages, the initial condition is generated from the output of the previous stage (see Sec.~\ref{sec:algo1} and Sec.~\ref{sec:algo2}). The primal boundary conditions obviously need no adjustment in the course of dual time evolution.  These initial and boundary conditions enter into the dual formulation as natural boundary conditions in each stage.
    \item To evaluate the quality of the results produced via the algorithm, the corresponding exact entropy solution to each example is presented in Appendix \ref{app:exact_sol}.
    \end{itemize}

The figures presented in this work are produced as follows: 
For each Gauss timeline, we compute the average of the primal data obtained from the two Gauss points within each element through which the Gauss timeline traverses. This average is then assigned to the $x$-coordinate of the element center and $t$-coordinate of the Gauss timeline under consideration.

We demonstrate and discuss the results of the computation of five selected examples that are used to evaluate and understand our dual scheme applied to the Burgers equation and its corresponding Hamilton-Jacobi form. 
Sec.~\ref{sec:Burgers_examples} and Sec.~\ref{sec:HJ_examples} deal with the inviscid case of the Burgers equation and Burgers-HJ, respectively, and we often omit the adjective `inviscid.' Specific details related to each of these equations are presented below.

\subsection{Algorithmic details for the Burgers equation}\label{sec:algo1}
 In order to evaluate solutions to the dual form of  the Burgers equation \eqref{eq:Burgers_primal_in_dual}, we utilize the residual \eqref{eq:weak_dual_Burgers} and  generate the corresponding Jacobian by considering its variation in a direction $d\lambda$:
\begin{equation*}
 J[\delta \lambda;d\lambda] = \bigintsss_\Omega \,
    \left(-\p_t \delta\lambda \, \left(\frac{\partial \hat{u}}{\partial \mathcal{D}_i} d\mathcal{D}_i\right) - \hat{u}\,\p_x \delta\lambda \,\left(\frac{\partial \hat{u}}{\partial \mathcal{D}_i} d\mathcal{D}_i\right) \right)\,dt\,dx, 
\end{equation*}
where $\mathcal{D} := (\p_t \lambda,\p_x\lambda)$. 

We employ the approximate dual field \eqref{eq:discrete_duals} in the residual \eqref{eq:weak_dual_Burgers} to generate a discrete residual $R^A$ 
at each node $A$ given by
\begin{equation*}
 R^A = \bigintsss_\Omega \,
    \left(-\hat{u} \, \p_t N^A  - \frac{\hat{u}^2\,\p_x N^A}{2} \right)\,dt\,dx - \int_0^L u_0(x)\,N^A(x,0)\,dx
    - \int_0^T \,\frac{(u_l(t))^2}{2}\,N^A(0,t)\,dt. 
\end{equation*}
However, $R^A = 0$ is not imposed for all the nodes corresponding to the right and top boundaries as a consequence of the Dirichlet b.cs \eqref{eq:Burgers_dual_dbc}. The discrete Jacobian matrix associated with the pair of degrees of freedom {$(A, B)$} (the coefficient of the term $(\delta \lambda^A d\lambda^B)$) is given by
\begin{equation*}
\begin{gathered}
 J^{AB} = \int_\Omega \,
    \Big(-\p_t N^A - \hat{u}\,\p_xN^A\Big) \, \left(\frac{\partial \hat{u}}{\partial \mathcal{D}_i} \p_i  N^B\right)\,dt\,dx;  \\ 
    \frac{\p \hat{u}}{\p(\p_t \lambda)} = \frac{1}{\beta_u-\p_x \lambda}; \qquad
    \frac{\p \hat{u}}{\p(\p_x \lambda)} = \frac{\hat{u}}{\beta_u-\p_x \lambda},
    \end{gathered}
\end{equation*}
where $\p_1 (\cdot)= \p_t (\cdot)$ and $\p_2 (\cdot)= \p_x (\cdot)$.
 
Based on the algorithm presented in Table.~\ref{algo:euler_algorithm}, after truncation, the cutoff line for each stage (from which the next stage commences) is the $g_2$ Gauss timeline of the last layer of retained elements with $t_f^{(s)}$ the corresponding physical time. The data for $\hat{u}$  obtained at this time (required only at Gauss points) is denoted by $\hat{u}^{(s)}\left(x,t_f^{(s)} \right)$. The initial condition and the base state for the next stage is set as
\begin{equation*}
\begin{aligned}
 u_0^{(s+1)}(x) & = \hat{u}^{(s)}\left(x,t_f^{(s)}\right); \\
    \bar{u}^{(s+1)}(x,t) & = \mathcal{S}\left[\hat{u}^{(s)}\left(x,t_f^{(s)}\right)\right],
\end{aligned}
\end{equation*}
where $\mathcal{S}$ denotes a smoothing operator, and the specific details can be found in Appendix~\ref{app:Smoothing}. In essence, the primal fields generated in each stage, relying on the dual fields and their derivatives, may exhibit oscillations, particularly in the vicinity of any shock; the smoothing operation mitigates these high-wave number oscillations in the base state for the subsequent stage that remains constant throughout the subsequent cycle. \textit{We note that such smoothing is not used to define the initial condition for the next stage}.
\subsection{Algorithmic details for the Burgers-HJ}\label{sec:algo2}

To assess solutions to \eqref{eq:HJ_primal_in_dual}, we employ the residual \eqref{eq:weak_dual_HJ} and generate a corresponding Jacobian by considering its variation in the direction $dD:=(d\lambda,d\gamma)$, expressed as
\begin{equation}
 J[\delta \lambda, \delta \gamma;d\lambda,d\gamma] = \bigintsss_\Omega \,
    \left(\bigl(-\p_t (\delta\lambda) -\p_x(\delta \gamma)\bigl) \, \left(\frac{\partial \hat{Y}}{\partial \mathcal{D}_i} d\mathcal{D}_i\right) + \bigl(\hat{u}\,\delta\lambda - \delta \gamma\bigl) \,\left(\frac{\partial \hat{u}}{\partial \mathcal{D}_i} d\mathcal{D}_i\right) \right)\,dt\,dx, 
\label{eq:jaco_HJ}
\end{equation}
where $\mathcal{D} := (\p_t \lambda,\p_x\gamma,\lambda,\gamma)$. 
The discrete version for the residual \eqref{eq:weak_dual_HJ}, computed using the approximate dual fields $D$ as defined in \eqref{eq:discrete_duals}, for each node $A$ is given by
\begin{equation*}
 R_1^A = \bigintsss_\Omega \,
    \Biggl(-\hat{Y}\,\p_t N^A \,  + \frac{\hat{u}\,^2 N^A}{2} 
     \Biggl)\,dt\,dx - \int_0^L \,Y_0(x)\,(N^A(x,0))\,dx   ; \end{equation*}
     \begin{equation*}
     R_2^A = \int_\Omega \,
    -\bigl( 
    \hat{Y}\,\p_x N^A  + \hat{u}\,N^A \bigl) \,dt\,dx-\int_0^T \,Y_l(t)\,(N^A(0,t))\,dt.
\label{eq:resi_HJ_discrete}
\end{equation*}
$R_1^A=0$ for all the nodes corresponding to the top boundary  and  $R_2^A=0$ for all the nodes corresponding to the right boundary are not imposed as a consequence of the Dirichlet b.cs \eqref{eq:HJ_dual_BC}. The discrete version of the Jacobian corresponding to the degree of freedom pair {$(A,i),(B,j)$} is given as follows:
\begin{equation*}
\begin{gathered}
 J_{11}^{AB} = \bigintsss_\Omega \,\left(
    -\frac{\partial \hat{Y}}{\partial (\p_t \lambda)} \,\p_t N^A \,\p_t N^B +  \hat{u}\,\frac{\p \hat{u}}{\p \lambda}\, N^A  \, N^B \right)\,dt\,dx
      ;\\
    J_{12}^{AB} = \bigintsss_\Omega \,\left(
    -\frac{\partial \hat{Y}}{\partial (\p_x \gamma)}\,\p_t N^A  \,\p_x N^B +  \hat{u}\, \frac{\p \hat{u}}{\p \gamma}\,N^A  \, N^B\right)\,dt\,dx;\\
    J_{21}^{AB} = \bigintsss_\Omega \,\left(
    -\frac{\partial \hat{Y}}{\partial (\p_t \lambda)}\,\p_x N^A\,  \p_t N^B - \frac{\p \hat{u}}{\p \lambda}\, N^A \, N^B \right)\,dt\,dx
    ;\\
    J_{22}^{AB} = \bigintsss_\Omega \,\left(
    -\frac{\partial \hat{Y}}{\partial (\p_x \gamma)} \,\p_x N^A  \,\p_x N^B - \frac{\p \hat{u}}{\p \gamma}\,N^A \, N^B\right)\,dt\,dx;\\
\frac{\partial \hat{Y}}{\partial (\p_t \lambda)} = \frac{1}{\beta_Y}; \quad \frac{\partial \hat{Y}}{\partial (\p_x \gamma)} = \frac{1}{\beta_Y};
\\
\frac{\p \hat{u}}{\p \lambda} = -\frac{\hat{u}}{\beta_u + \lambda}; \quad \frac{\p \hat{u}}{\p \gamma} = \frac{1}{\beta_u + \lambda}.
    \end{gathered}
\end{equation*}

 In contrast to the examples associated with the Burgers equation, the oscillations observed in the solution (for $\hat{Y}$ here) obtained on the Gauss points are significantly more pronounced.
To address this issue, in each stage, the truncation is followed by an $L^2$ projection of $\hat{Y}$  (Appendix.~\ref{app:L2}) along the $n_2$ nodal timeline of the last layer of elements obtained on the retained mesh. Such an operation provides us with a continuous $C^0$ approximant along the considered nodal timeline. This timeline also serves as the cutoff line for the current stage.
The initial condition and the base states for the following stage are set as follows:
\begin{equation*}
\begin{gathered}
    Y_0^{(s+1)}(x) = L^2\left[Y^{(s)}\left(x,t_f^{(s)}\right)\right]; \\
 \bar{Y}^{(s+1)}(x,t)=Y^{(s)}\left(x,t_f^{(s)}\right); \quad \bar{u}^{(s+1)} = \frac{d}{dx}Y^{(s)}\left(x,t_f^{(s)}\right), 
\end{gathered}
\end{equation*}
where $L^2[\cdot]$ represents the $L^2$ projection operator. Unlike the scheme for the Burgers equation, no smoothing for base states is used for the inviscid Burgers-HJ algorithm.

\section{Results}\label{sec:results}
To evaluate the dual scheme, we  demonstrate five examples for the two types of equations, each presenting distinct levels of complexity. For each of the examples introduced in Sec.~\ref{sec:Burgers_examples}, a related example is generated in Sec.~\ref{sec:HJ_examples}. This associated example is supplied with the initial and boundary conditions in a manner that the exact entropy solution produced by the Burgers-HJ equation should have its derivative equal to the exact entropy solution obtained for the corresponding problem from the Burgers equation.
\subsection{Burgers Equation}\label{sec:Burgers_examples}
For each of the examples presented in this section, the following parameters were set for each stage in the algorithm: $T^{(s)}=5\times10^{-3}$, $N_x=100$, $N_t^{(s)}=100$, $\beta_k=10^6$, $N _c = 5$ and $tol=10^{-16}$. This amounts to a fairly fine discretization of the problem, but our main goal here is not to probe the computational efficiency of the proposed elementary scheme, but to simply use it to evaluate and understand the capabilities of our dual variational formulation for nonlinear PDE in this specific context.
\subsubsection{Expansion Fan}\label{sec:Burgers_fan}
We consider the Riemann problem corresponding to an expansion fan solution with initial condition given by:
\begin{equation}
    u_0(x) = \begin{cases}
U_L = 0 &  \mbox{for } x< 0.5 \\ 
U_R = 1 &  \mbox{for } x>0.5.
\end{cases} \label{eq:ini_Burgers_fan}
\end{equation}
The problem has non-unique weak solutions (consisting of the continuous rarefaction wave and an entire family of discontinuous solutions with shocks), and a unique `entropy' solution is obtained by imposing an entropy condition \cite{Lax_SIAM, GILBERT_STRANG} (also see App.~\ref{app:exact_sol}). This entropy solution is a rarefaction wave, which for $t>0$ in an infinite domain is given by \eqref{eq:Fan_ex_u}.

\begin{figure}
\centering
\begin{subfigure}{.29\textwidth}
  \centering
  \includegraphics[width=.93\linewidth]{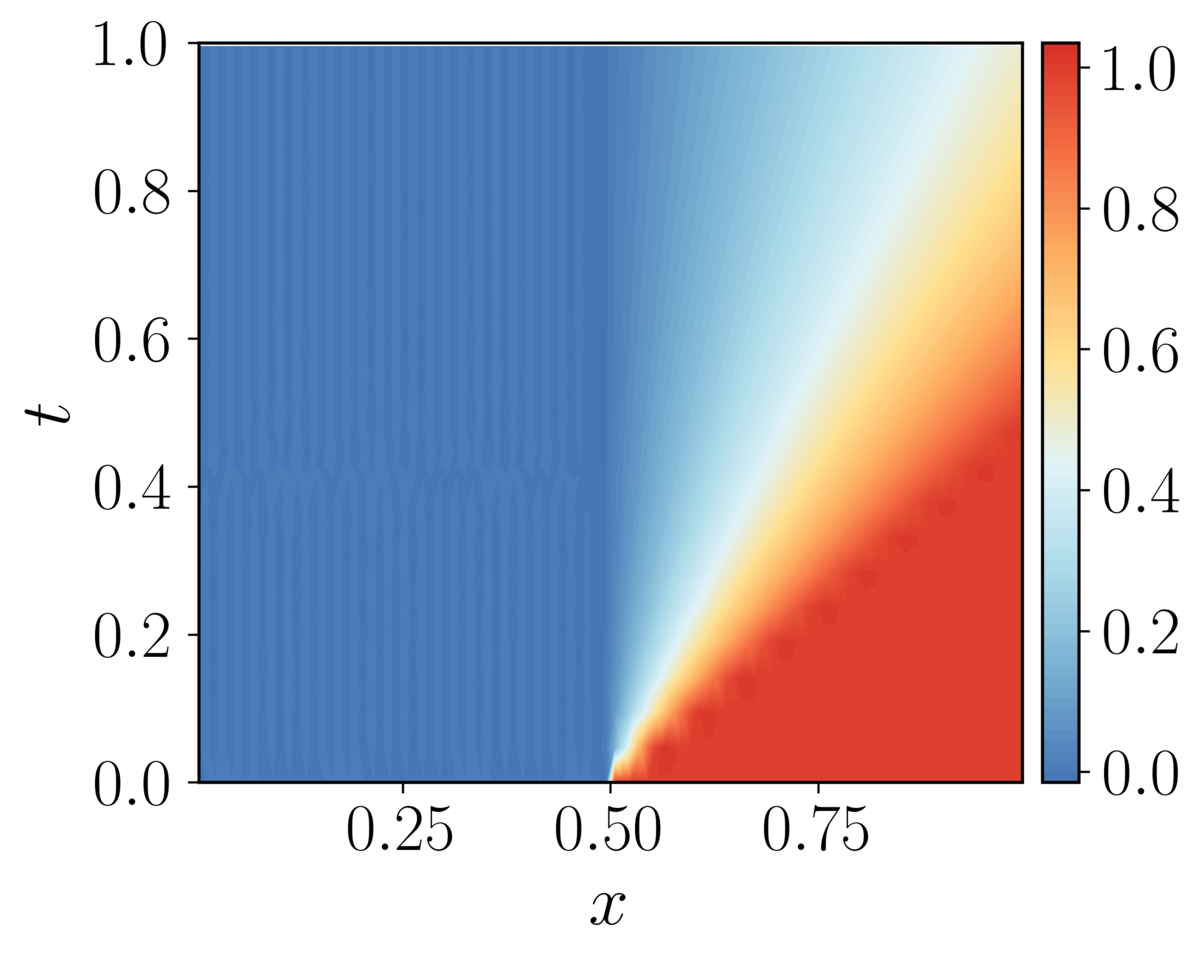}  
  \caption{$u(x,t)$}
  \label{fig:sub-bur_p1_heatmap}
\end{subfigure}
\begin{subfigure}{.29\textwidth}
  \centering
  \includegraphics[width=.93\linewidth]{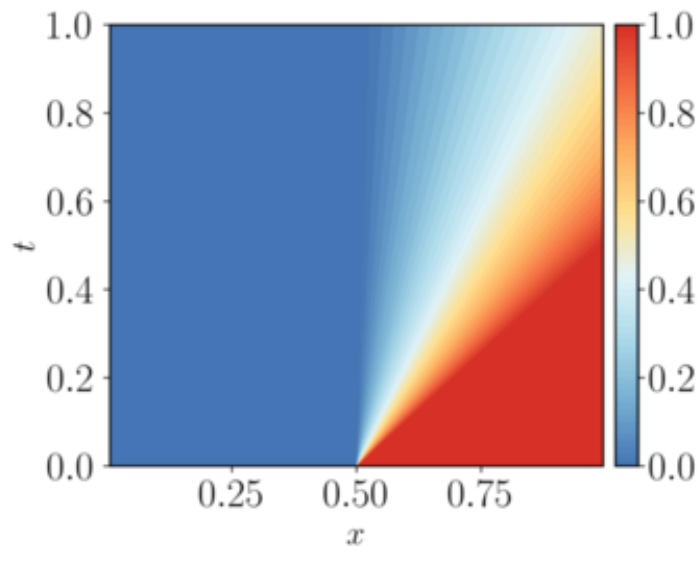}  
  \caption{Exact $u(x,t)$ based on \eqref{eq:Fan_ex_u}} 
  \label{fig:sub-bur_p1_ex}
\end{subfigure}
\begin{subfigure}{.4\textwidth}
  \centering
  \includegraphics[width=.9\linewidth]{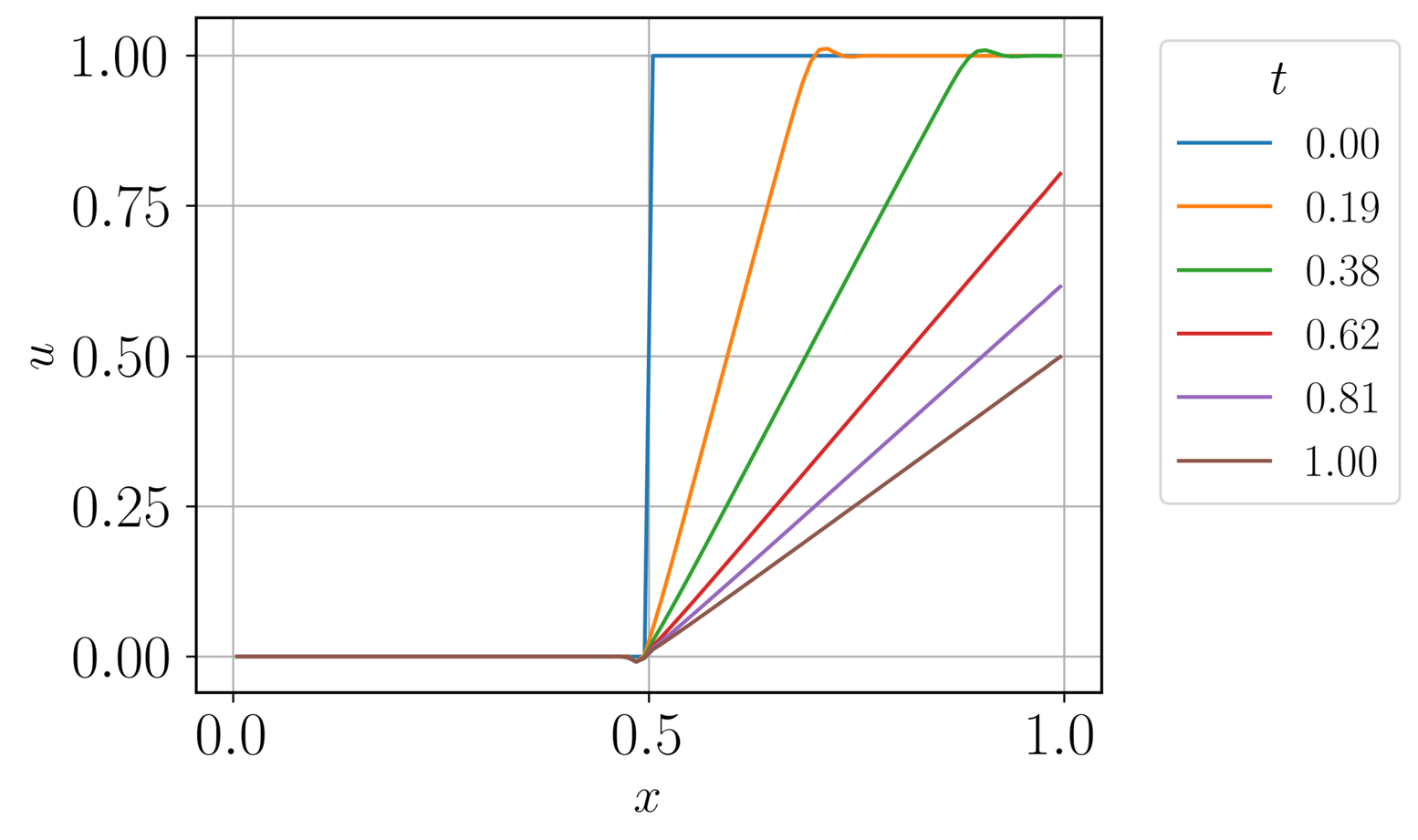}  
  \caption{Line plots for $u(x,t)$}
  \label{fig:sub-bur_p1_line}
\end{subfigure}
\caption{DtP mapping generated primal field $u$ for the expansion fan. Fig.~(c) shows minor overshoots as the fan opens up which may be attributed to the $C^0$ approximation of the dual fields.} 
\label{fig:bur_p1}
\end{figure}

We  use the primal boundary condition $u_l(t) = 0$ and apply the Dirichlet b.cs $\lambda_T(x) = \lambda_r(t) = 0$.  The results for this setup are shown in Fig.~\ref{fig:bur_p1}. 

Comparing Fig.~\ref{fig:sub-bur_p1_heatmap} against Fig.~\ref{fig:sub-bur_p1_ex}, it is surprising to see that the dual scheme automatically picks up only the entropy solution without enforcing any further conditions, and seems to be incapable of recovering the other possible weak solutions for the equation.

\subsubsection{Shock}\label{sec:Burgers_shock}
We next consider the Riemann problem \eqref{eq:ini_Burgers_fan} with $$U_L = 1; \qquad U_R=0.$$ The exact entropy solution to this problem on an infinite domain is given by \eqref{eq:Shock_ex_u}.

\begin{figure}
\begin{subfigure}[t]{.49\textwidth}
  \centering
  \includegraphics[width=0.7\linewidth]{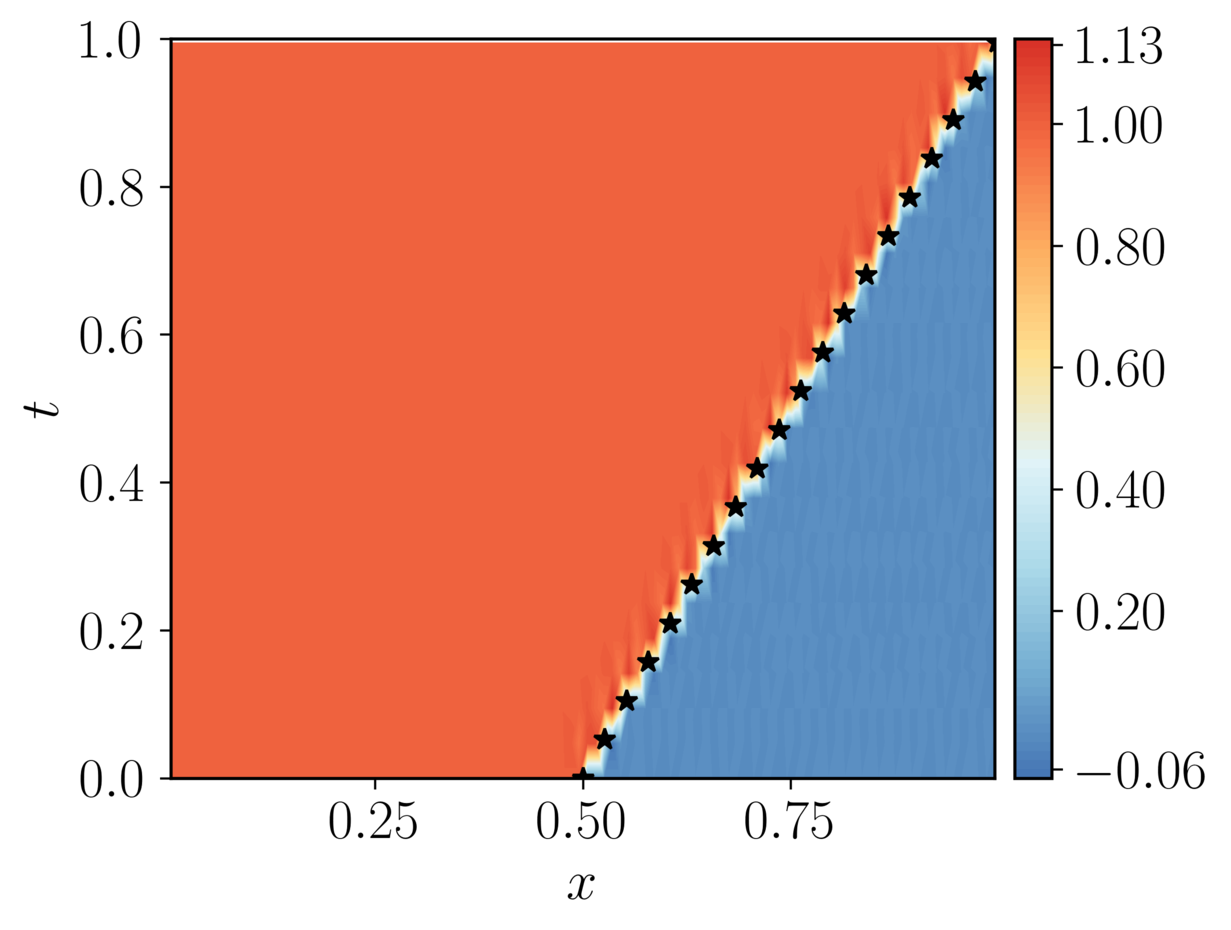}  
  \caption{$u(x,t)$ }
  \label{fig:sub-bur_p2_1}
\end{subfigure}
\begin{subfigure}[t]{.49\textwidth}
  \centering
  \includegraphics[width=0.85\linewidth]{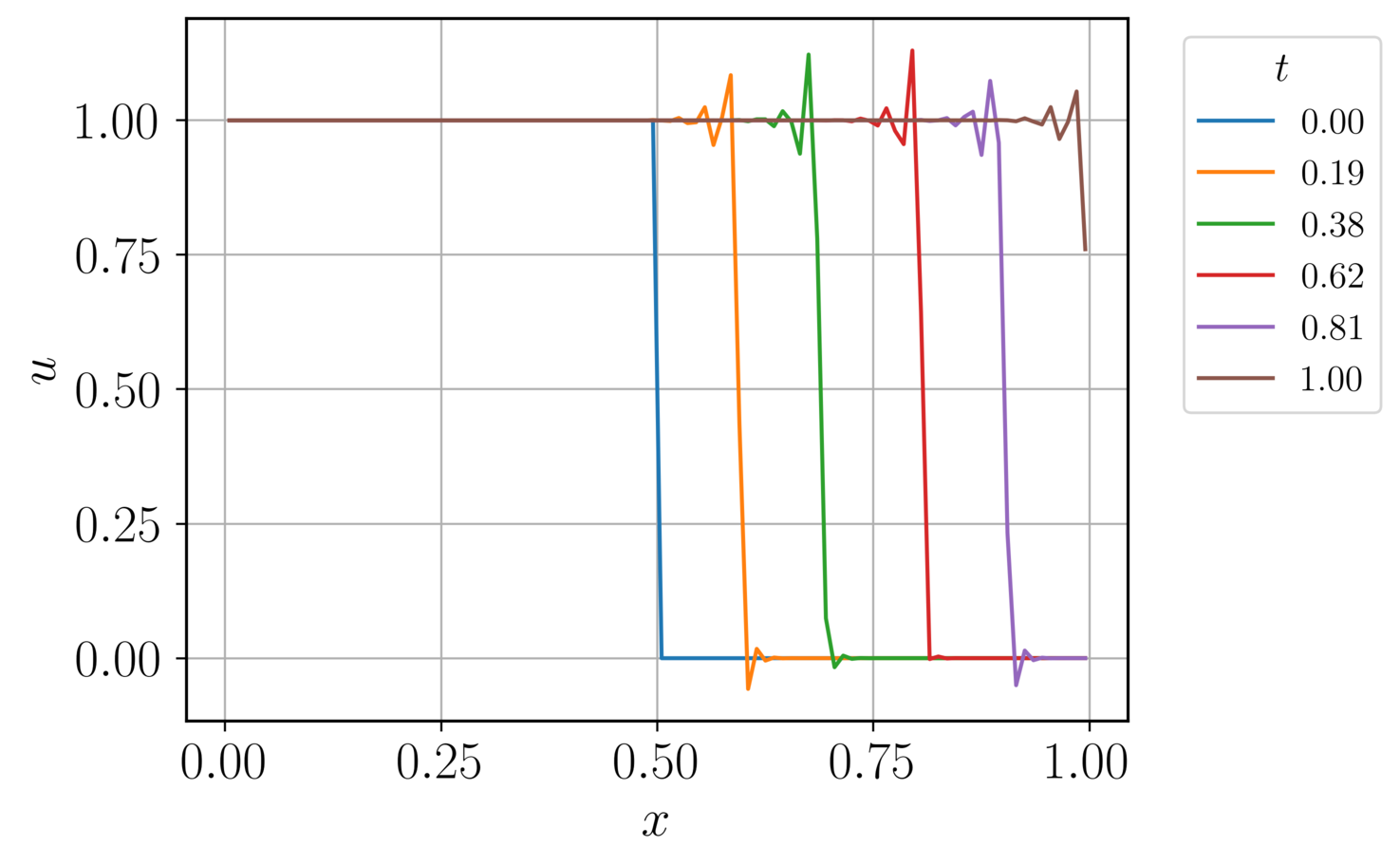}  
  \caption{Line plots for $u(x,t)$}
  \label{fig:sub-bur_p2_2}
\end{subfigure}
\caption{DtP mapping generated primal field $u$ for the shock problem . In Figure (a), the black asterisks represent the shock trajectory based on the exact entropy solution, which is being superimposed on the original plot.}
\label{fig:bur_p2}
\end{figure}

We  use the primal boundary condition $u_l(t) = 1$ and apply the Dirichlet b.cs $\lambda_T(x) = \lambda_r(t) = 0$. The results for this setup has been shown in Fig.~\ref{fig:bur_p2}. The over(under)shoot seen around the shock profile at any time in Fig.~\ref{fig:sub-bur_p2_2} arise due to the $C^0$ FE interpolation of the dual fields. Nevertheless, the shock profile effectively captures the accurate height and speed, aligning closely with the exact entropy solution.

\subsubsection{Double Shock}\label{sec:Burgers_dshock}
Here we assess the dual formulation in the setting of two traveling shocks interacting with each other. Consider the following initial condition
\begin{equation*}
u_0(x) =
\begin{cases}
    \begin{aligned}
        1 & \quad \text{for } x < 0.25 \\
        0.5 & \quad \text{for } 0.25 < x < 0.5 \\
        0 & \quad \text{for } x > 0.5.
    \end{aligned}
\end{cases}
\end{equation*}
The solution to this problem on an infinite domain is given by \eqref{eq:DShock_ex_u}. 

\begin{figure}
\begin{subfigure}[t]{.49\textwidth}
  \centering
  \includegraphics[width=0.7\linewidth]{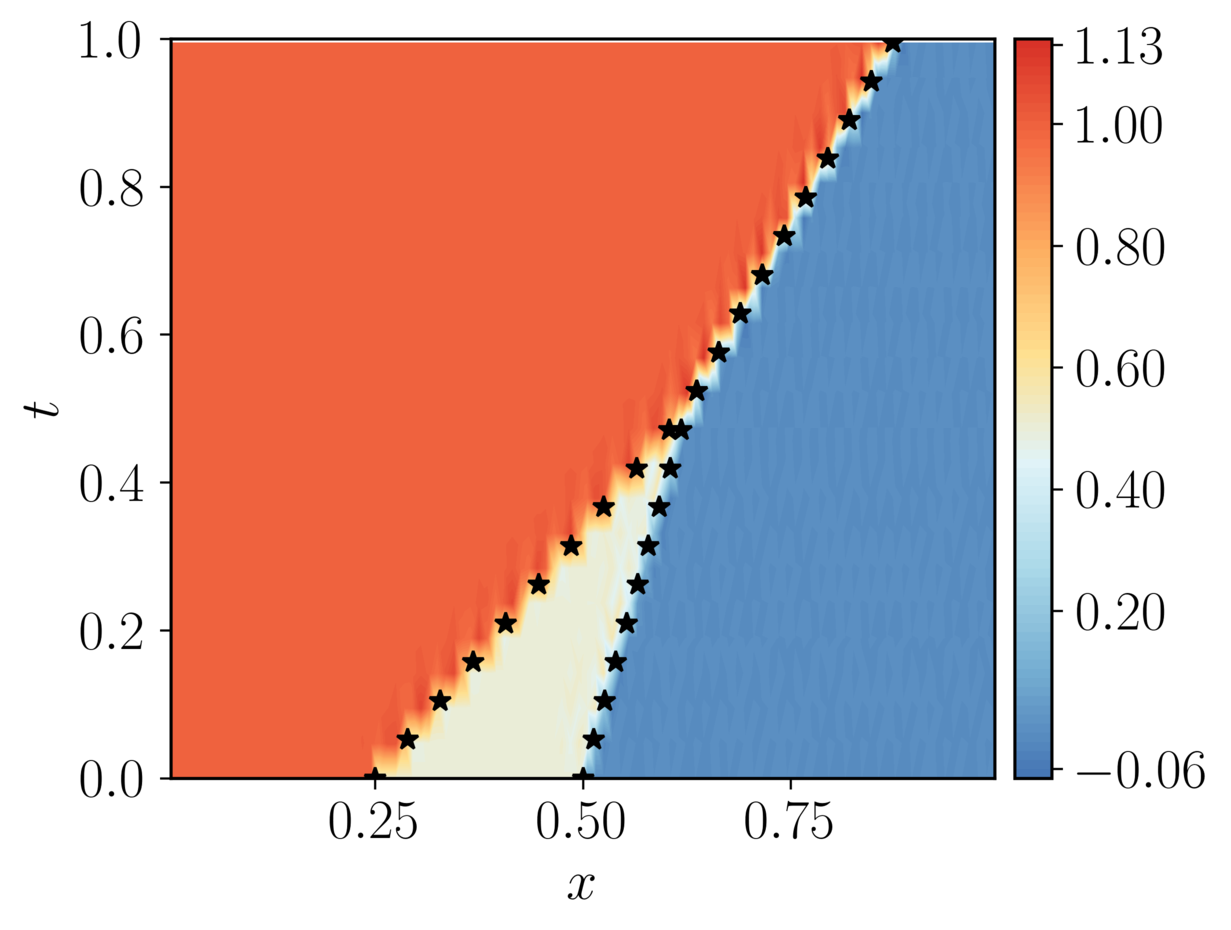}  
  \caption{$u(x,t)$ }
  \label{fig:sub-bur_p3_1}
\end{subfigure}
\begin{subfigure}[t]{.49\textwidth}
  \centering
  \includegraphics[width=0.85\linewidth]{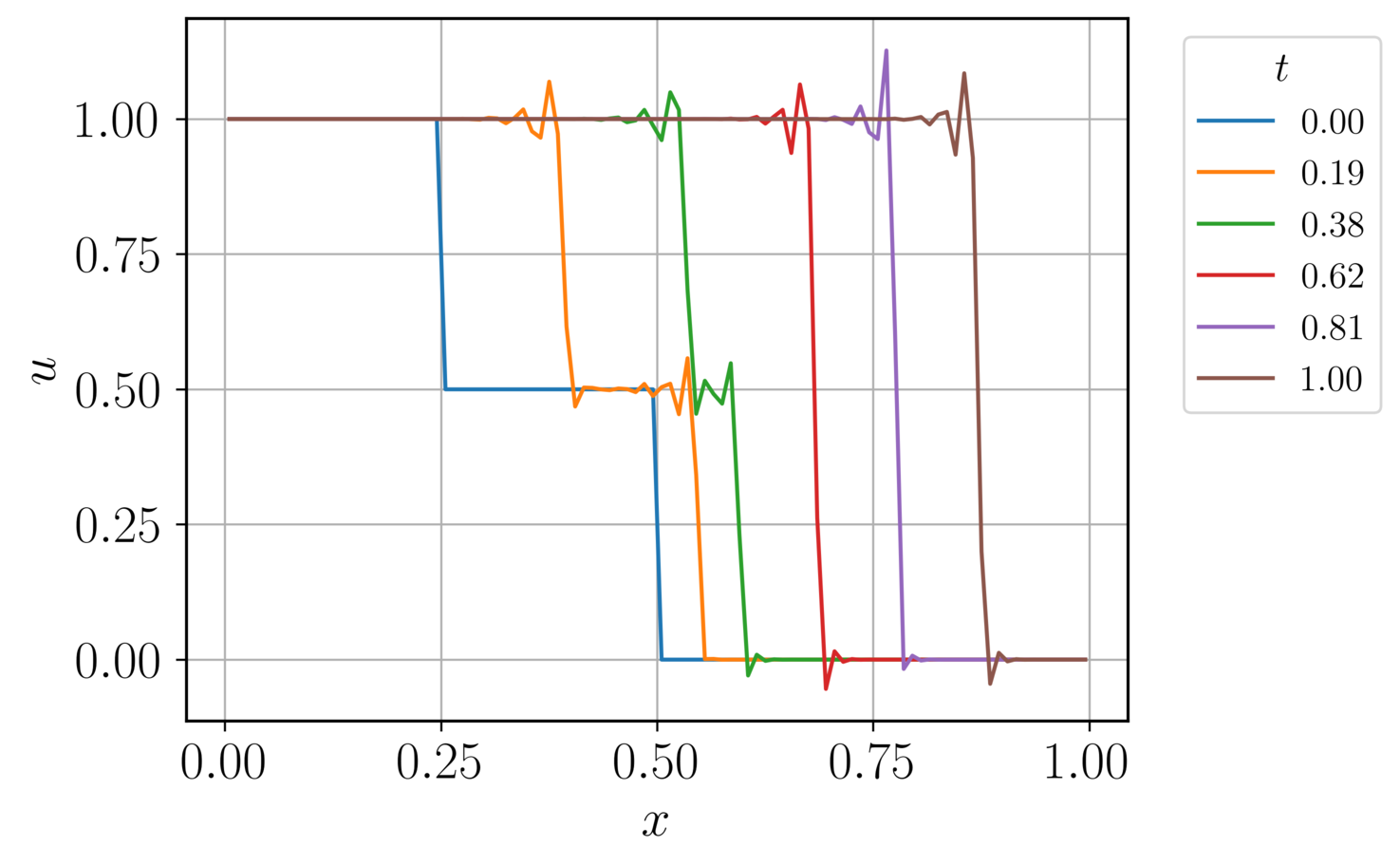}  
  \caption{Line plots for $u(x,t)$}
  \label{fig:sub-bur_p3_2}
\end{subfigure}
\caption{DtP mapping generated primal field $u$ for a double shock. In Figure (a), the black asterisks represent the shock trajectories for the two shocks based on the exact entropy solution, which is being superimposed on the original plot. The two shocks merge at $(x,t)=(0.625,0.5)$}
\label{fig:bur_p3}
\end{figure}

We  use the primal boundary condition $u_l(t) = 1$ and apply the Dirichlet b.cs $\lambda_T(x) = \lambda_r(t) = 0$.  The results for this setup are shown in Fig.~\ref{fig:bur_p3}.

 Evident form Fig.~\ref{fig:sub-bur_p3_2}, the instabilities near the two shocks, caused by the continuous interpolation, do not accumulate as the shocks merge. Instead, they vanish, leaving no indication of a double shock presence. Based on Fig.~\ref{fig:sub-bur_p3_1}, it is also evident that the dual problem correctly captures the accurate shock speeds for each individual shock.

\subsubsection{Half N-wave}\label{sec:Burgers_halfwave}

\begin{figure}
\begin{subfigure}[t]{.49\textwidth}
  \centering
  \includegraphics[width=0.7\linewidth]{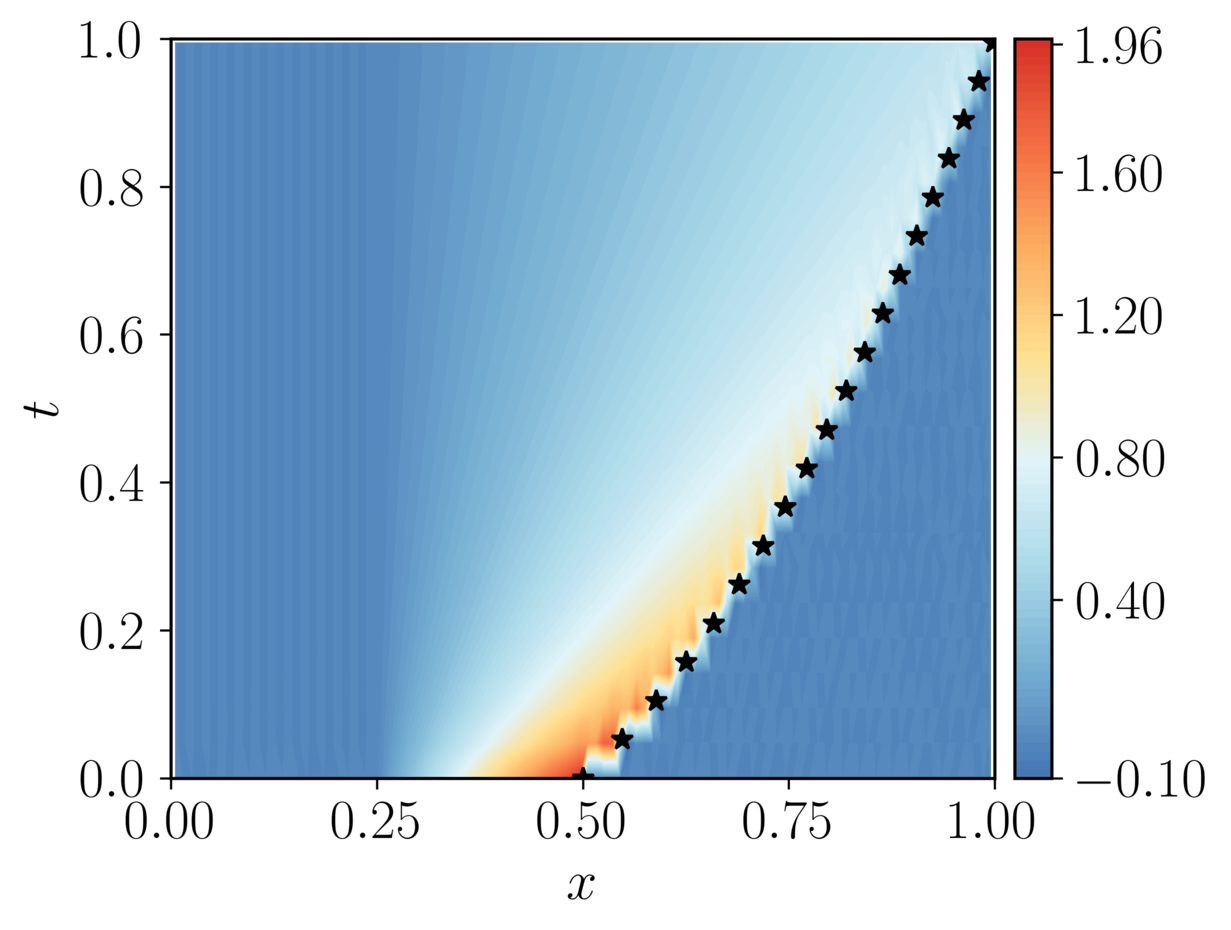}  
  \caption{$u(x,t)$ }
  \label{fig:sub-bur_p4_1}
\end{subfigure}
\begin{subfigure}[t]{.49\textwidth}
  \centering
  \includegraphics[width=0.85\linewidth]{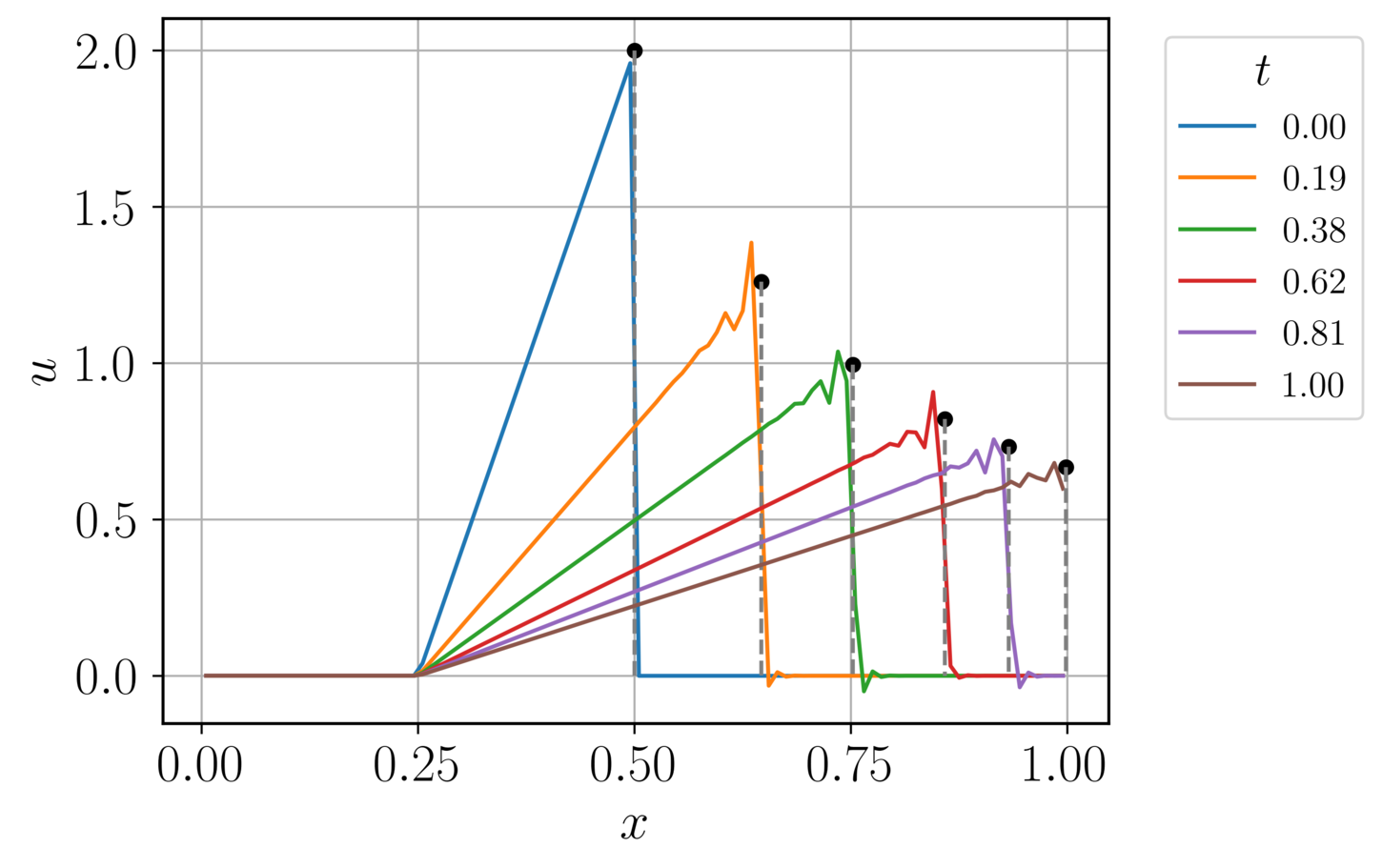}  
  \caption{Line plots for $u(x,t)$}
  \label{fig:sub-bur_p4_2}
\end{subfigure}
\caption{DtP mapping generated primal field $u$ for the half N-wave problem. In Fig.(a), the black asterisks represent the shock trajectory based on the exact entropy solution, which is being superimposed on the original plot. In Figure (b), the adjacent black lines, aligned with the original plot lines at various time points, depict the position and height of the shock based on the exact entropy solution for those respective times.} 
\label{fig:bur_p4}
\end{figure}

The objective of this problem is to evaluate the dual formulation when a shock wave decreases in size and experiences a change in speed during its propagation Fig.~\ref{app:halfwave}. We call this problem as a \textit{half N-wave}, adapting standard terminology in the literature (even though it looks like a `reflected' half $N$). Consider the following initial condition:
\begin{equation*}
\begin{aligned}
u_0(x) =
\begin{cases}
    0 & \text{for } x < x_0 \\[2pt]
    \dfrac{h_0}{l_0}(x-x_0) & \text{for } x_0 \leq x < x_0 + l_0 \\[4pt]
    0 & \text{for } x > x_0 + l_0,
\end{cases}
\end{aligned}
\end{equation*}
where $x_0=0.25$, $l_0=0.25$ and $h_0=2$. The solution to this problem on an infinite domain is given by \eqref{eq:halfwave_ex_u}. 

We  use the primal boundary condition $u_l(t) = 0$ and apply the Dirichlet b.cs $\lambda_T(x) = \lambda_r(t) = 0$.  The results for this setup are shown in Fig.~\ref{fig:bur_p4}.

The geometry in this example is more intricate compared to previous examples, featuring both a fan at $x=0.25$ and a traveling shock. As explained in Appendix.~\ref{app:halfwave}, conservation of the area under the right-angled triangle dictates the speed and the height of the shock, which vary non-linearly in time. The dual formulation successfully captures this phenomenon.

\subsubsection{N-wave} \label{sec:burgers_N-wave}
The objective of this problem is to evaluate the dual formulation in the context of  two expansion fans in opposite directions converging to form a standing shock whose magnitude dissipates over time (Fig.~\ref{fig:N-wave_char}). Consider the following initial condition:
\begin{equation*}
\begin{aligned}
u_0(x) =
\begin{cases}
    0 & \text{for } x < 0.25 \\
    -4h_0(x-0.5) & \text{for } 0.25 \leq x < 0.75 \\
    0 & \text{for } x > 0.75;
\end{cases} 
\end{aligned}
\end{equation*}
\begin{equation*}
    h_0 = 2.
\end{equation*}
The solution to this problem on an infinite domain is given by \eqref{eq:Nwave_ex_u}.

We  use the primal boundary condition $u_l(t) = 0$ and apply the Dirichlet b.cs $\lambda_T(x) = \lambda_r(t) = 0$.

\begin{figure}
\begin{subfigure}[t]{.3\textwidth}
  \centering
  \includegraphics[width=\linewidth]{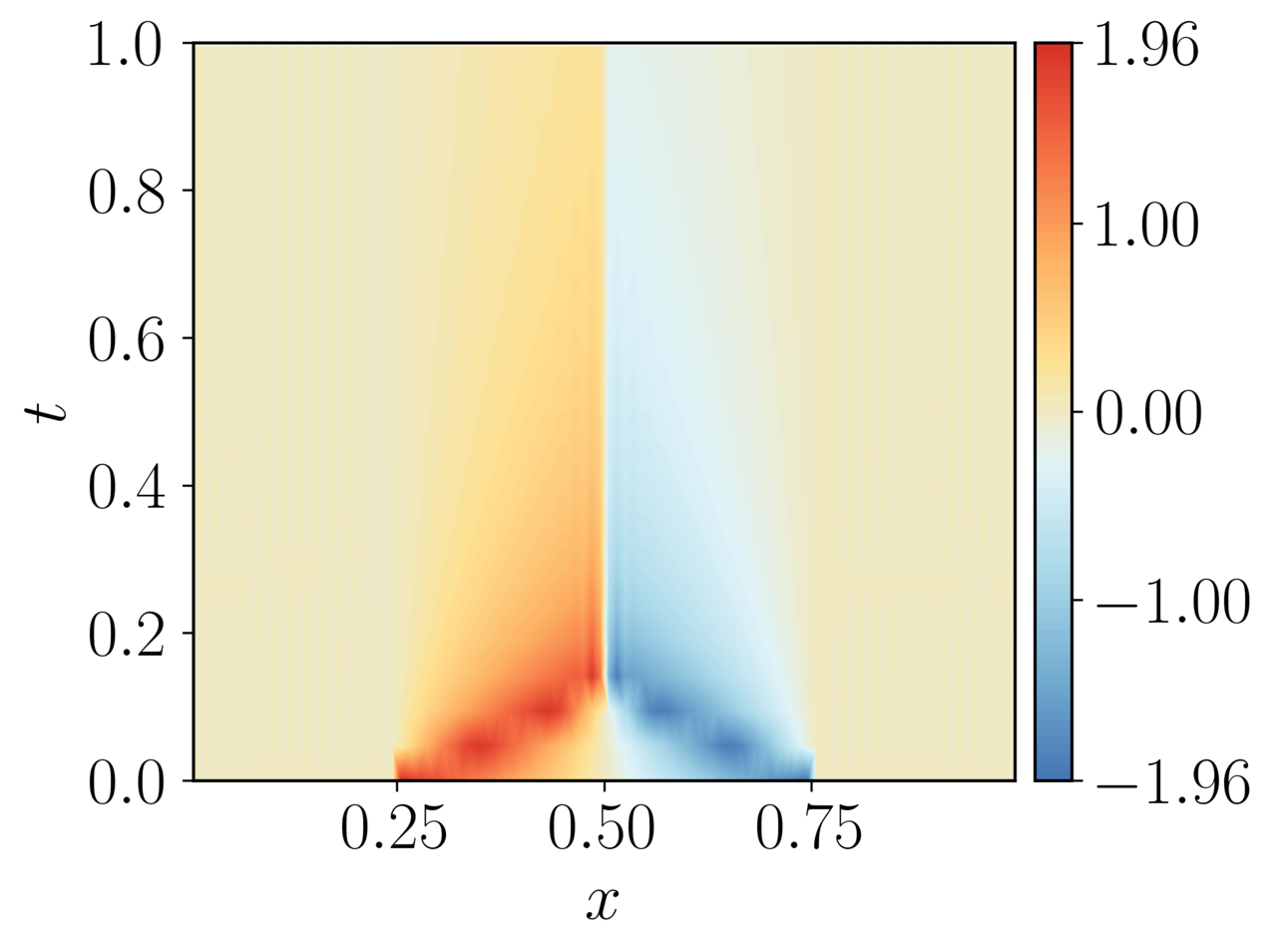}  
  \caption{$u(x,t)$ }
  \label{fig:sub-bur_p5_1}
\end{subfigure}
\begin{subfigure}[t]{.3\textwidth}
  \centering
  \includegraphics[width=0.95\linewidth]{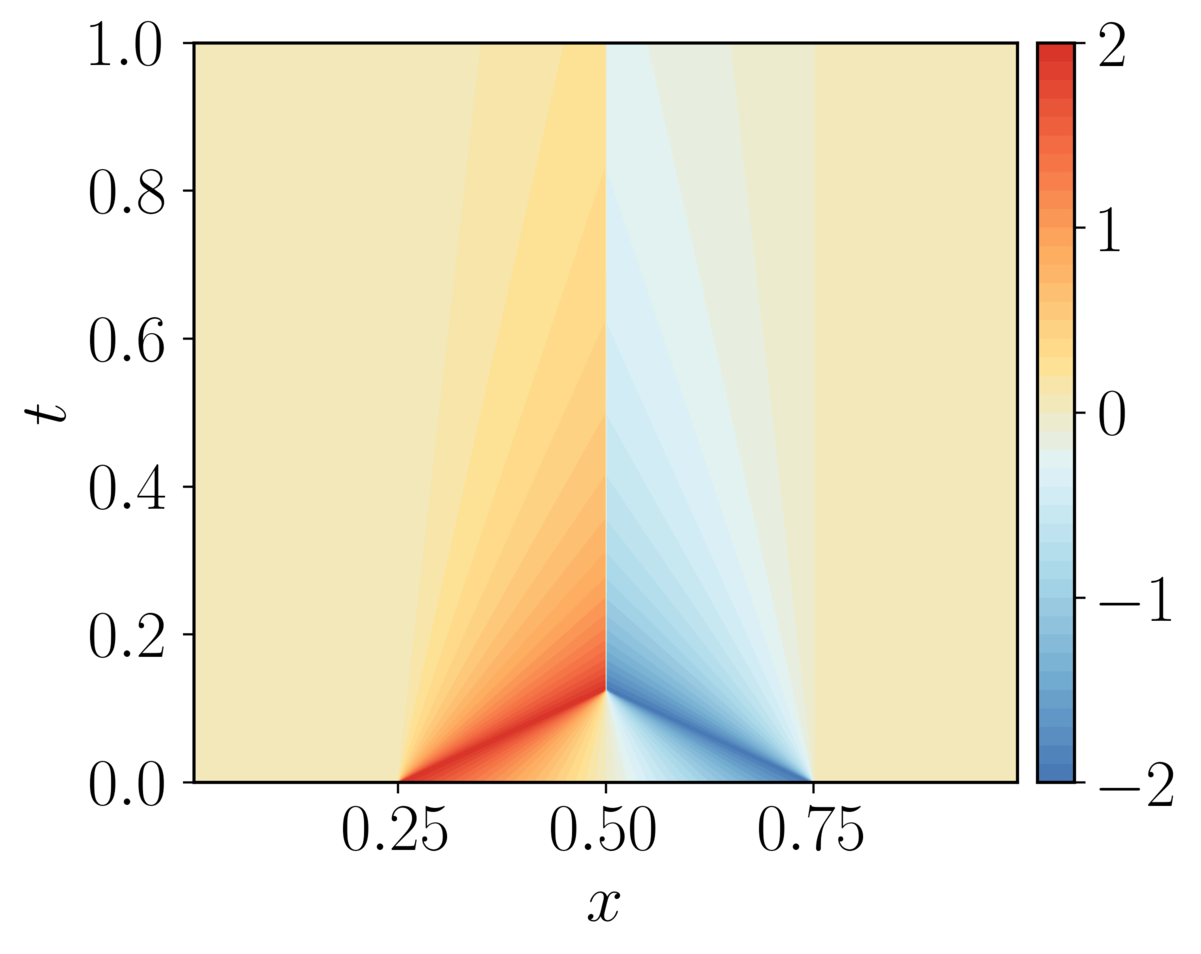}  
  \caption{Exact $u(x,t)$ based on \eqref{eq:Nwave_ex_u}}
  \label{fig:sub-bur_p5_ex_field}
\end{subfigure}
\begin{subfigure}[t]{.39\textwidth}
  \centering
  \includegraphics[width=0.9\linewidth]{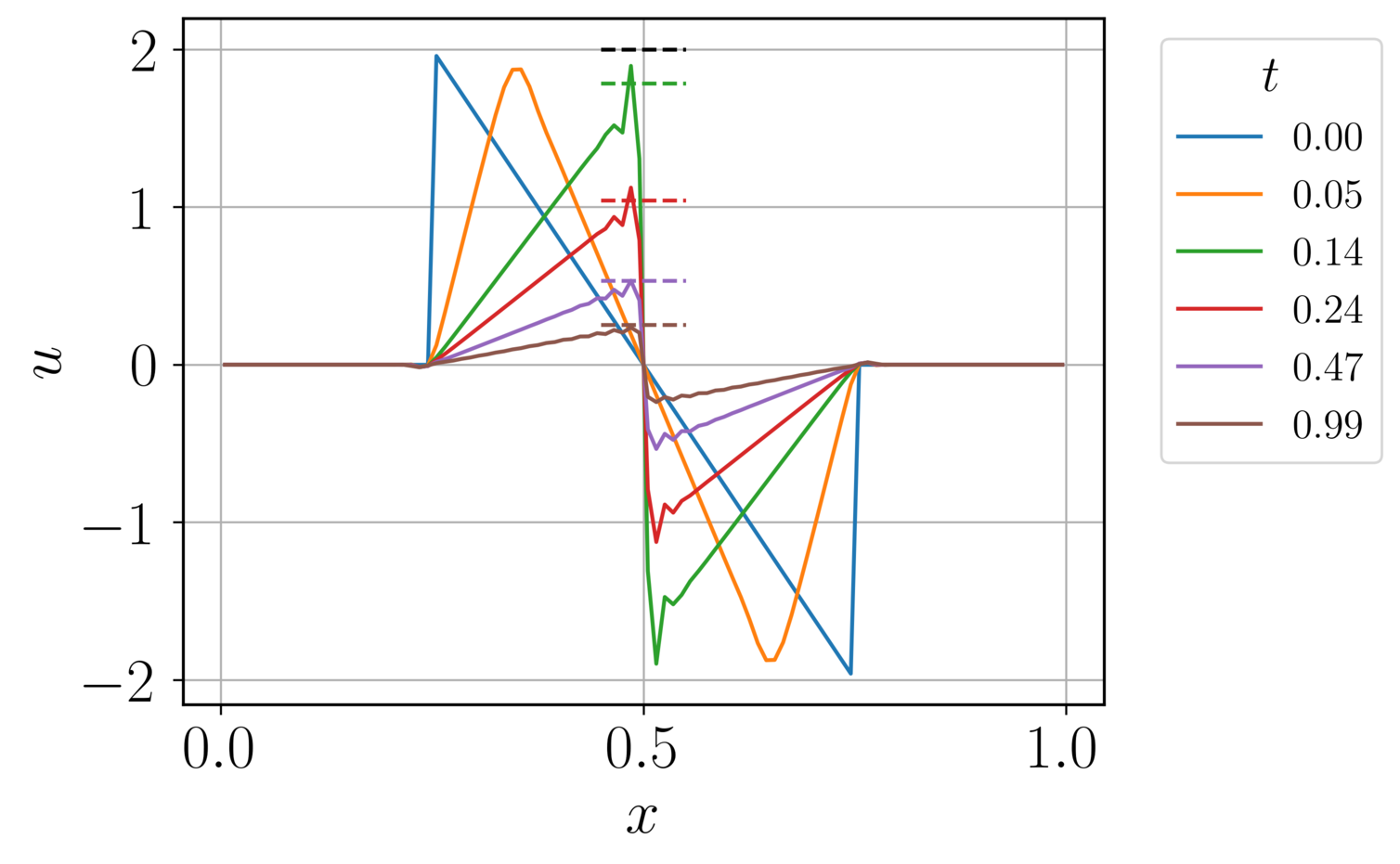}  
  \caption{Line plots for $u(x,t)$}
  \label{fig:sub-bur_p5_2}
\end{subfigure}
\caption{DtP mapping generated primal field $u$ for the N-wave problem.  In Fig.~\ref{fig:sub-bur_p5_2}, the dashed lines use colors matching those of the line plots at specific times to represent the true height at those times according to the exact entropy solution. The black dashed line corresponds to the height of the tip before the shock formation occurs.} 
\label{fig:bur_p5}
\end{figure}

Based on the discussion found in \ref{app_nwave} related to Fig.~\ref{fig:sub-bur_p5_ex}, each of the two tips at $x=0.25$ and $x=0.75$ travel in opposite directions while maintaining a constant speed and height. 
The information originating from the range $x\in(0.25,0.75)$ at time $t=0$ ultimately converges to form a standing shock at $x=0.5$ and $t_m =0.125$. Subsequently, the standing shock gradually self-annihilates. The results based on the dual scheme for this setup, shown in Fig.~\ref{fig:bur_p5}, manages to capture this whole sequence of events quite well.

\subsection{Inviscid Burgers-HJ}\label{sec:HJ_examples}
For the examples presented in this section, the following parameters were set for each stage in the algorithm:  $\beta_Y=\beta_u=10^6$, $N _c = 5$ and $tol=10^{-16}$. The parameters defining the mesh were $N_x=50, N_t^{(s)}=10$ with $T^{(s)}=5\times10^{-5}$. This is again a fine discretization, and our concern here is simply to evaluate the dual formulation of the problem.

\subsubsection{Expansion Fan} \label{sec:HJ_fan}
\begin{figure}
\begin{subfigure}[t]{.49\textwidth}
  \centering
  \includegraphics[width=0.95\linewidth]{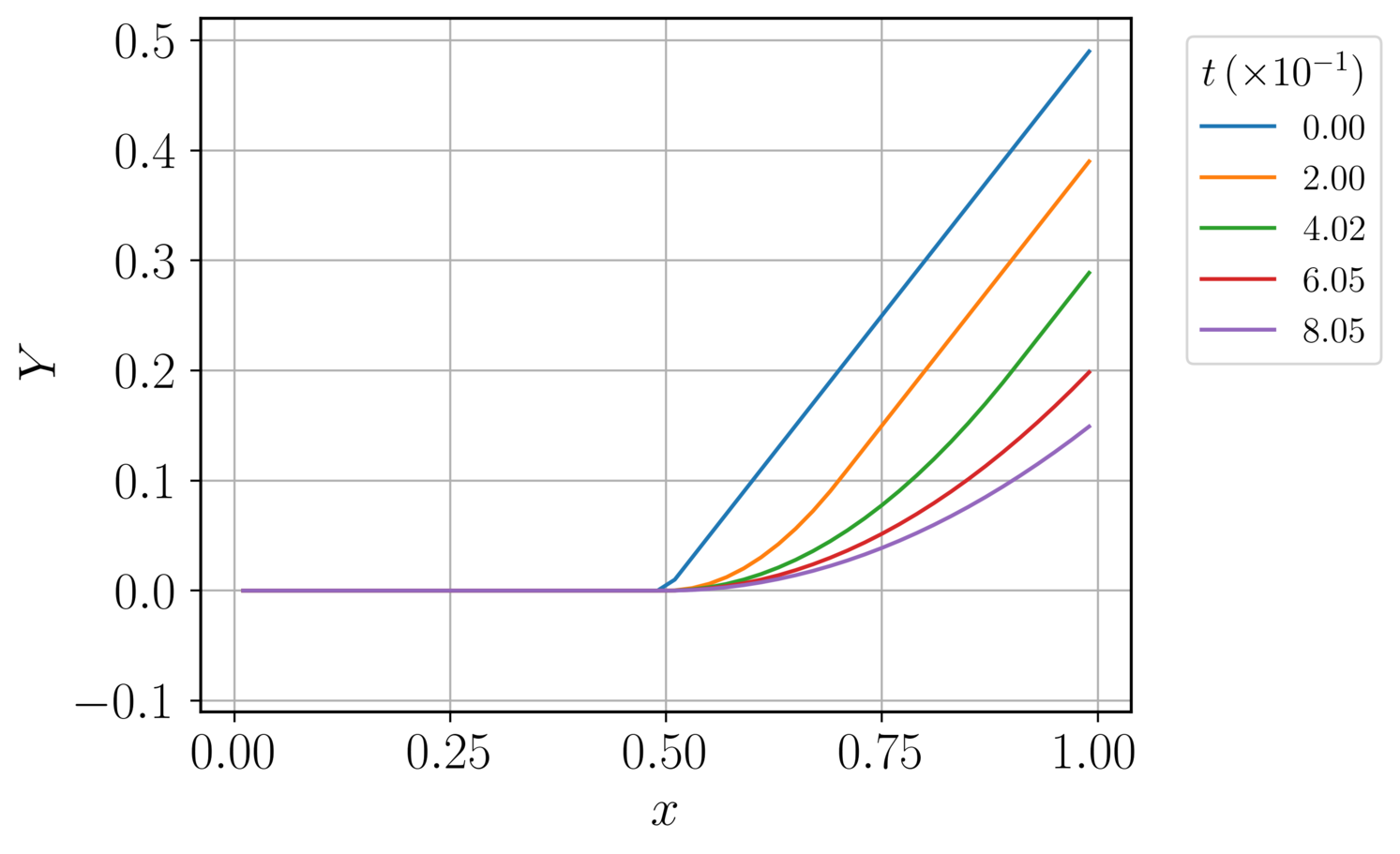}  
  \caption{Exact time plots for $Y(x,t)$}
  \label{fig:sub-burHJ_p1_1}
\end{subfigure}
\begin{subfigure}[t]{.49\textwidth}
  \centering
  \includegraphics[width=0.95\linewidth]{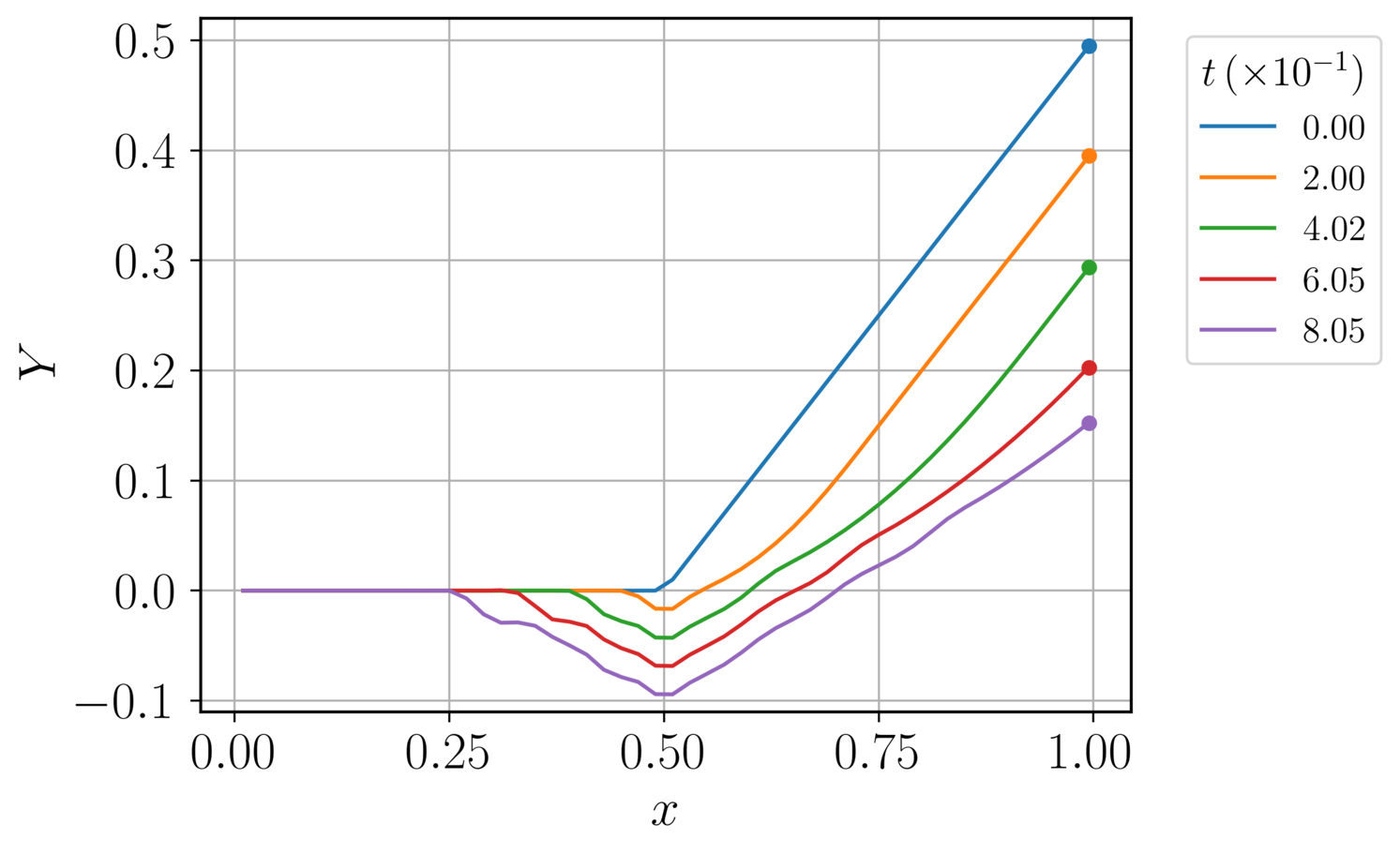}  
  \caption{Line plots for $Y(x,t)$ obtained using the dual scheme}
  \label{fig:sub-burHJ_p1_2}
\end{subfigure}
\caption{The two plots have been drawn on same scale. Fig.(b) shows DtP mapping generated primal field $Y$ for the Fan problem. The dots in various colors at $x=1$ correspond to the exact value of $Y$ derived from the exact entropy solution \eqref{eq:fan_ex_Y}, corresponding to the time indicated by the color of each respective dot.} 
\label{fig:burHJ_p1}
\end{figure}

Consider the following initial condition:
\begin{equation*}
    Y_0(x) = \begin{cases}
0 &  \mbox{for } x< 0.5 \\ 
x &  \mbox{for } x>0.5.
\end{cases}
\end{equation*}
The exact entropy solution for this setup on an infinite domain is given by \eqref{eq:fan_ex_Y}. To evaluate the dual solutions numerically on $\Omega$, we impose $Y_l(t)=0$, naturally. The results obtained for this setup are show in Fig.~\ref{fig:burHJ_p1}, which shows an unexpected dip in the $Y$ profiles at $x=0.5$. 

This discrepancy indicates that the dual scheme does not capture the entropy solution for this particular setup. Upon closer inspection one can notice that the dip actually represents a standing shock which can be confirmed by approximating the spatial gradient ($u$) to the immediate left and to the right of $x=0.5$, and using the jump condition such that 
\[ c\,(x=0.5,t) = \frac{\p_x Y(x=0.5^+,t) + \p_x Y(x=0.5^-,t)}{2} \approx \frac{1 + (-1)}{2}=0.\]
Thus the evolution profiles indicates that a fan breaks up into a locally standing shock which then propagates information in both directions. We note that within an expansion fan, a shock is acceptable as a weak solution.

\subsubsection{Shock}
\begin{figure}
\begin{subfigure}[t]{.49\textwidth}
  \centering
  \includegraphics[width=0.75\linewidth]{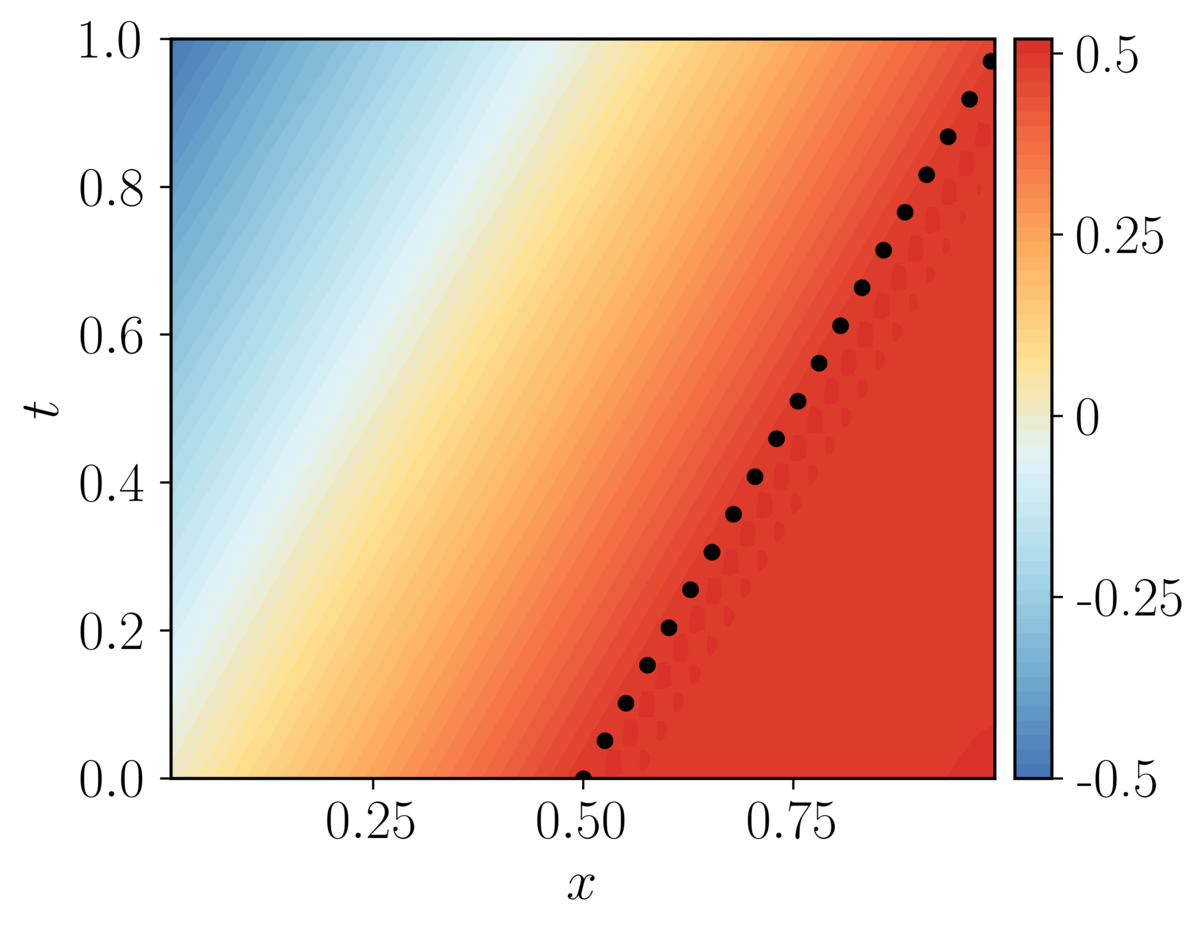}  
  \caption{$Y(x,t)$}
  \label{fig:sub-burHJ_p2_1}
\end{subfigure}
\begin{subfigure}[t]{.49\textwidth}
  \centering
  \includegraphics[width=0.95\linewidth]{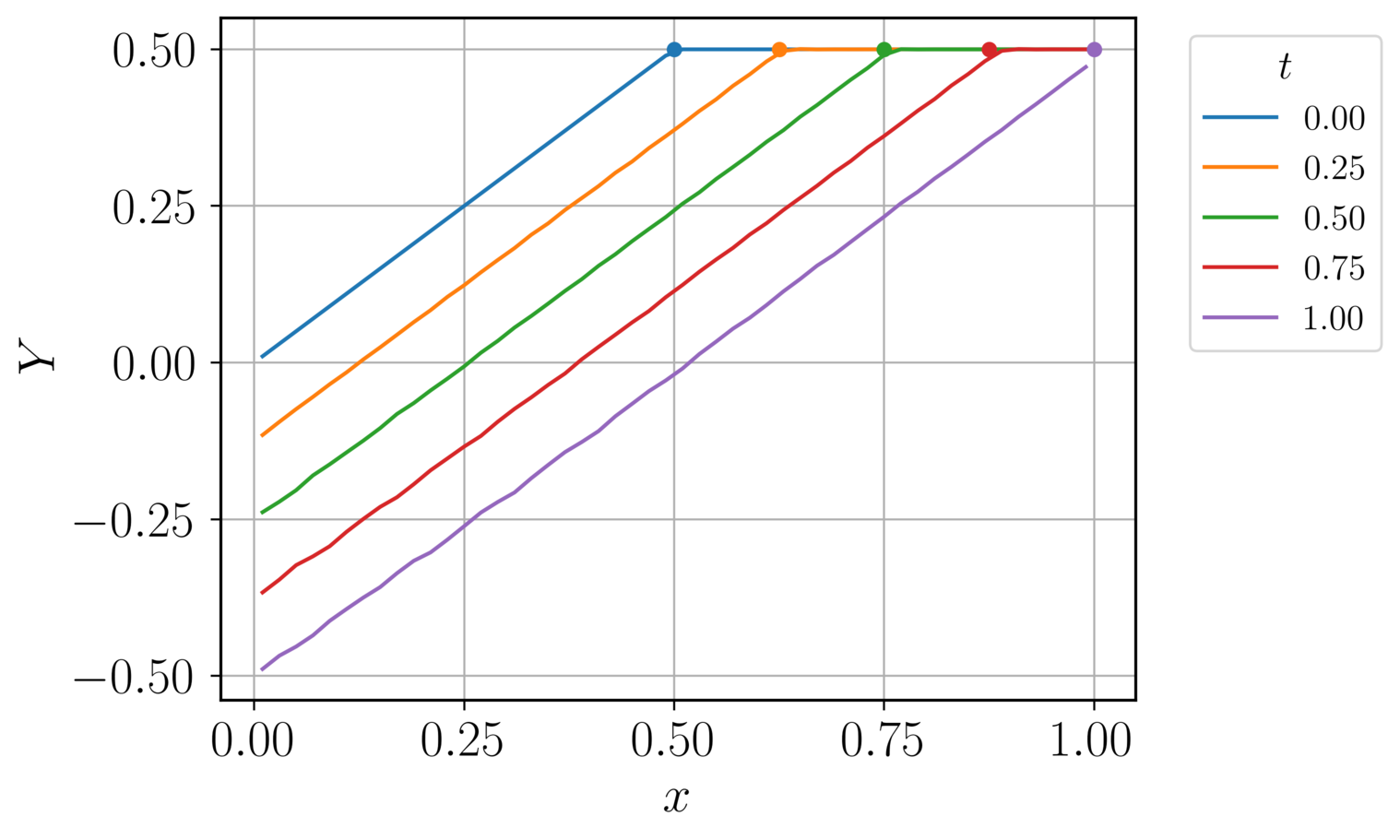}  
  \caption{Line plots for $Y(x,t)$}
  \label{fig:sub-burHJ_p2_2}
\end{subfigure}
\caption{DtP mapping generated primal field $Y$ for the Shock problem. In Figure(a), the black dots represent the trajectory of the kink based on the exact entropy solution \eqref{eq:Shock_ex_Y}. In Figure (b), the dots plotted at $Y=0.5$ with colors matching those of the line plots represent the position of the kink based on the exact $Y$ profile \eqref{eq:Shock_ex_Y} at that specific time.} 
\label{fig:burHJ_p2}
\end{figure}

Consider the initial condition and the boundary condition given by
\begin{equation*}
\begin{gathered}
    Y_0(x) = \begin{cases}
x &  \mbox{for } x< 0.5 \\ 
0.5 &  \mbox{for } x>0.5;
\end{cases} \\
Y_l(t) = -\frac{t}{2},
\end{gathered}
\end{equation*}
respectively. For these conditions, we expect to recover the solution given by \eqref{eq:Shock_ex_Y}, the derivative of which recovers the shock presented in \ref{sec:Burgers_shock}. 
The $Y$ data obtained at different times from the simulation of the dual problem has been shown in Fig.~\ref{fig:burHJ_p2}. 

Referring to Fig.~\ref{fig:sub-burHJ_p2_2}, it is apparent that kink in the $Y$ profiles moves at an approximate speed of 0.5. This kink signifies the presence of a shock in the $u$ profile, where the slope of the $Y$ (equivalent to $u$) profile changes from 1 to 0.  As depicted in Fig. \ref{fig:sub-burHJ_p2_2}, minor oscillation-type instabilities are still noticeable in the $Y$ profiles as time advances.

\subsubsection{Double Shock}
\begin{figure}
\begin{subfigure}[t]{.49\textwidth}
  \centering
  \includegraphics[width=0.75\linewidth]{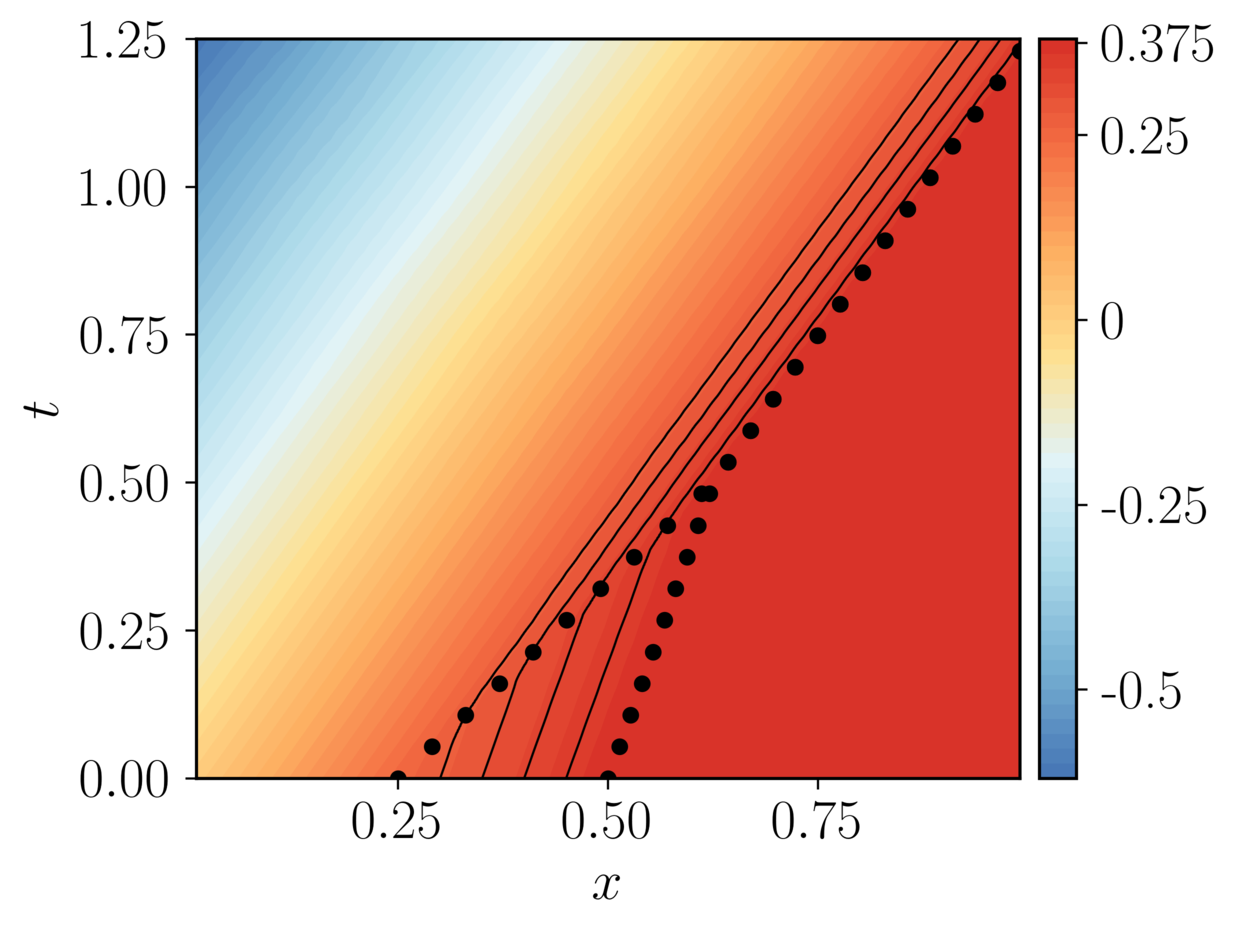}  
  \caption{$Y(x,t)$}
  \label{fig:sub-burHJ_p3_1}
\end{subfigure}
\begin{subfigure}[t]{.49\textwidth}
  \centering
  \includegraphics[width=0.95\linewidth]{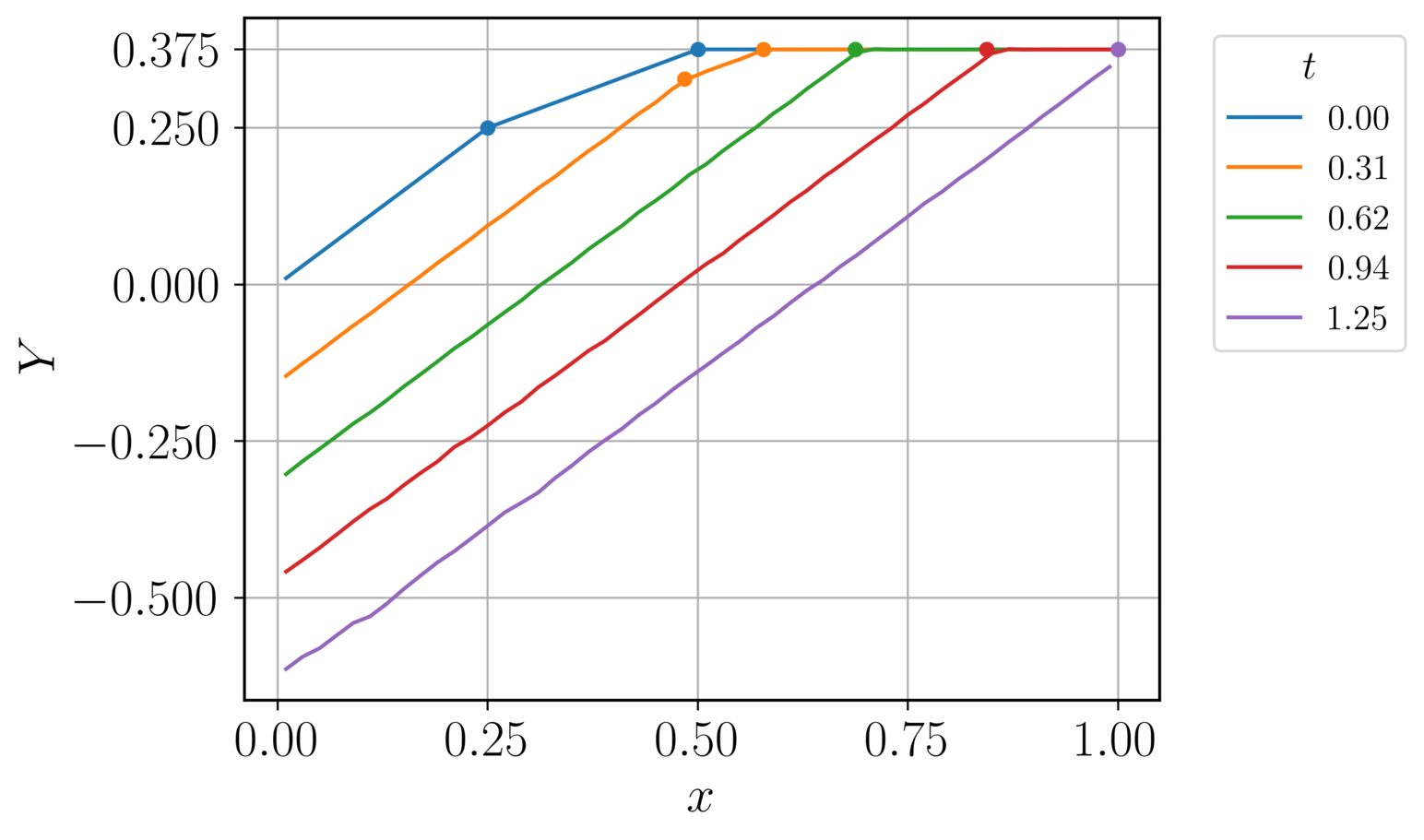}  
  \caption{Line plots for $Y(x,t)$}
  \label{fig:sub-burHJ_p3_2}
\end{subfigure}
\caption{DtP mapping generated primal field $Y$ for the double shock. In Fig.(a), the black dots represent the trajectory of the two kinks in the space-time domain based on the exact expressions \eqref{eq:DShock_ex_Y}. The solid black lines represent contour lines corresponding to the characteristics of $u$ derived through dual scheme in the vicinity of the double shock. In Figure (b), the dots plotted with colors matching those of the line plots represent the position and height of the kinks based on the exact $Y$ profile \eqref{eq:DShock_ex_Y} at that specific time.} 
\label{fig:burHJ_p3}
\end{figure}

Consider the initial condition and the boundary condition given by
\begin{equation*}
\begin{gathered}
Y_0(x) =
\begin{aligned}
\begin{cases}
    
        x & \quad \text{for } x < 0.25 \\
        x/2+0.125 & \quad \text{for } 0.25 < x < 0.5 \\
        0.375 & \quad \text{for } x > 0.5; 
\end{cases} 
\end{aligned} \\
Y_l(t) = -\frac{t}{2},
\end{gathered}
\end{equation*}
respectively. For these conditions, we expect to recover the solution given by \eqref{eq:DShock_ex_Y}, the derivative of which recovers the double shock presented in \ref{sec:Burgers_dshock}. The results for this setup are shown in Fig.~\ref{fig:burHJ_p3}.

Similar to the shock example, this problem has two kinks emnating at $t=0$ corresponding to the two shocks, traveling at different speeds, merging at $t=0.5$. Apparent from  Fig. \ref{fig:sub-burHJ_p3_1}, the characteristics originating between $x=0.25$ and $x=0.5$ move at a slower speed compared to the approaching kink from the left, eventually merging with it. The dual formulation captures this phenomenon, although (very minor) instabilities are again visible as evident in Fig.~\ref{fig:sub-burHJ_p3_2}


\subsubsection{Half N-wave}\label{sec:HJ_halfwave}
\begin{figure}
\begin{subfigure}[t]{.49\textwidth}
  \centering
  \includegraphics[width=0.95\linewidth]{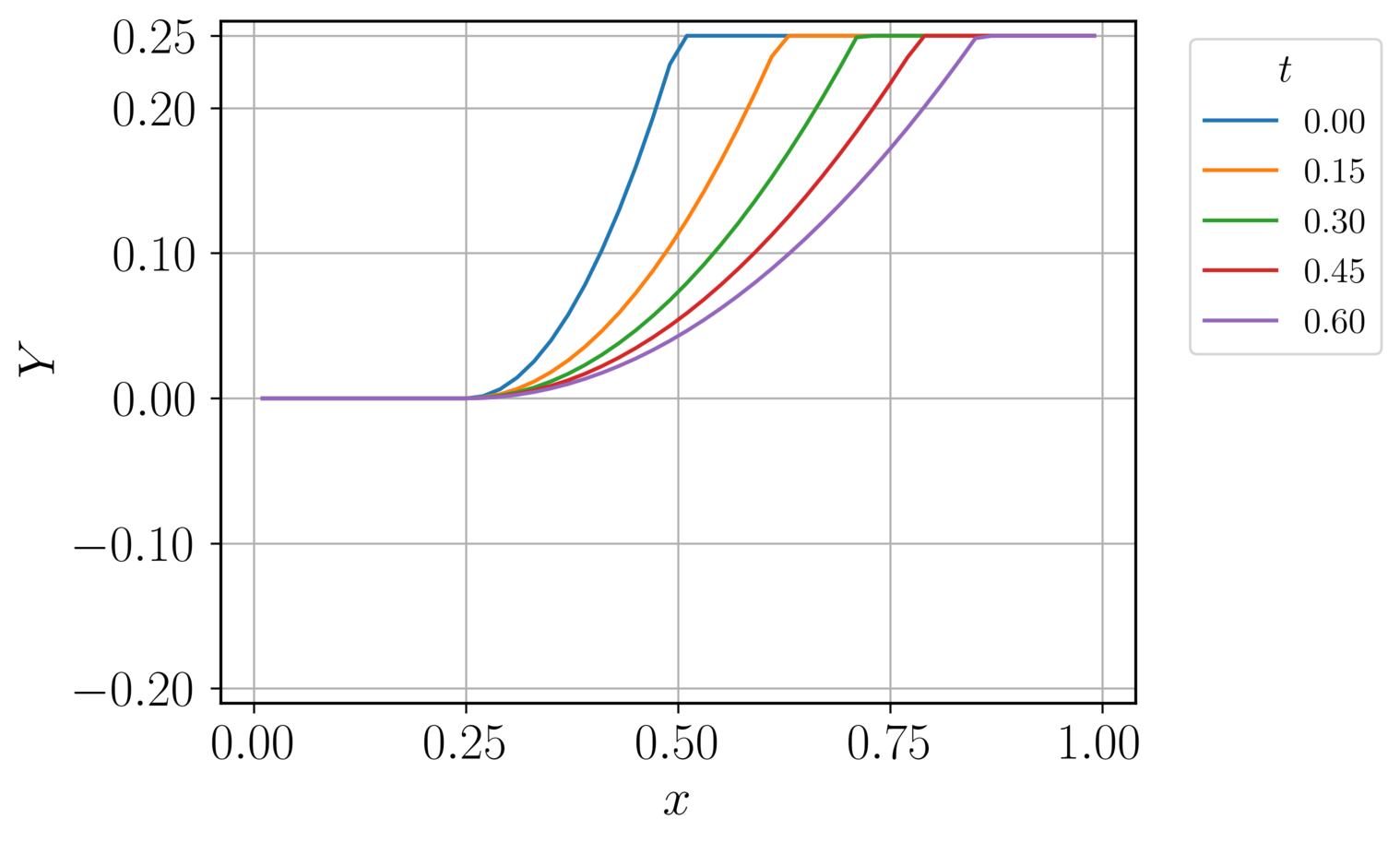}  
  \caption{Exact time plots for $Y(x,t)$}
  \label{fig:sub-burHJ_p4_1}
\end{subfigure}
\begin{subfigure}[t]{.49\textwidth}
  \centering
  \includegraphics[width=0.95\linewidth]{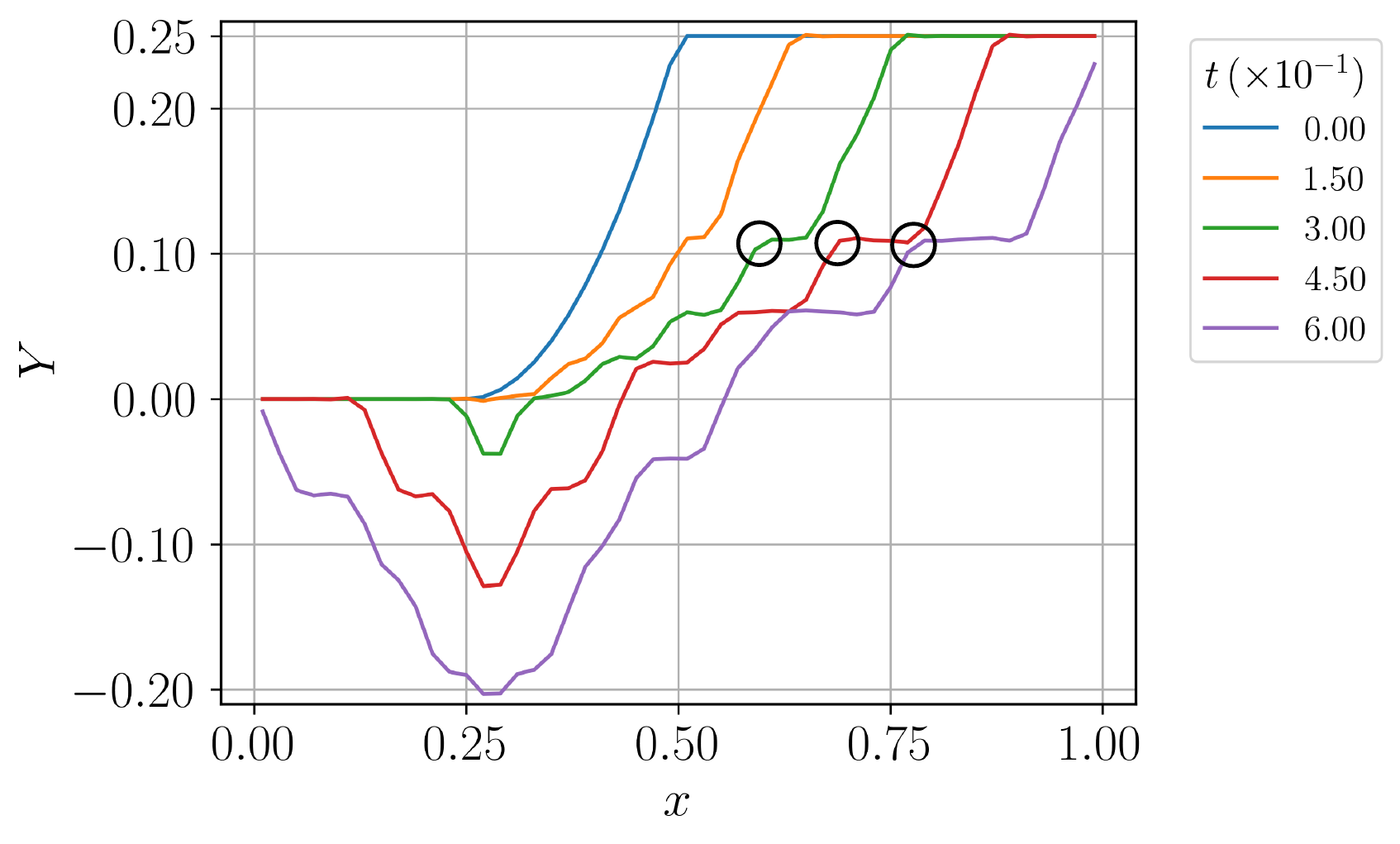}  
  \caption{Line plots for $Y(x,t)$ obtained using the dual scheme}
  \label{fig:sub-burHJ_p4_2}
\end{subfigure}
\caption{The two plots have been drawn on same scale. Fig.(b) represents the DtP mapping generated $Y$ profiles. The black circles along the kinks represent the trajectory of one of the multiple kinks whose spatial derivatives exhibit a shock-like behavior.} 
\label{fig:burHJ_p4}
\end{figure}

Consider the following equation set:
\begin{equation*}
\begin{aligned}
Y_0(x) &=
\begin{cases}
    0 & \text{for } x < x_0 \\[4pt] 
    \dfrac{h_0}{l_0}\left(\dfrac{x^2+x_0^2}{2}-x_0\,x\right) & \text{for } x_0 \leq x < x_0 + l_0 \\[10pt] 
    \dfrac{h_0\,l_0}{2} & \text{for } x \geq x_0 + l_0.
\end{cases} \\
&x_0=0.25;\quad l_0=0.25;\quad h_0=2.
\end{aligned}
\end{equation*}
The exact expression for $Y$ for this setup is given by \eqref{eq:halfwave_ex_Y}, the spatial derivative of which represents the Half N-wave presented in Sec.~\ref{sec:Burgers_halfwave}. The results obtained via our dual scheme is shown in Fig.~\ref{fig:burHJ_p4}. 

Comparing against the profiles for $Y$ based on the exact entropy solution illustrated in \ref{fig:sub-burHJ_p4_1}, we immediately notice the standing shock behavior (as in the case of the fan example) here at $x=0.25$. Additionally, we notice several small shock-like structures nucleating one after other at different points in space as the curved section of the profile between $x = 0.25$ and $x = 0.5$ advances. A `shock' speed analysis similar to the fan example, but now on one such structure characterized by a kink for which the spatial derivative ($u$) approximately goes from 1 to 0 (shown by the circles in Fig.~\ref{fig:sub-burHJ_p4_2}), reveals that the distance travelled by this kink in a given time approximates well the distance suggested by the `shock speed' of the kink within that time, establishing that these structures indeed behave like shocks.

\subsubsection{N-wave}\label{sec:HJ_nwave}
\begin{figure}
\begin{subfigure}[t]{.49\textwidth}
  \centering
  \includegraphics[width=0.95\linewidth]{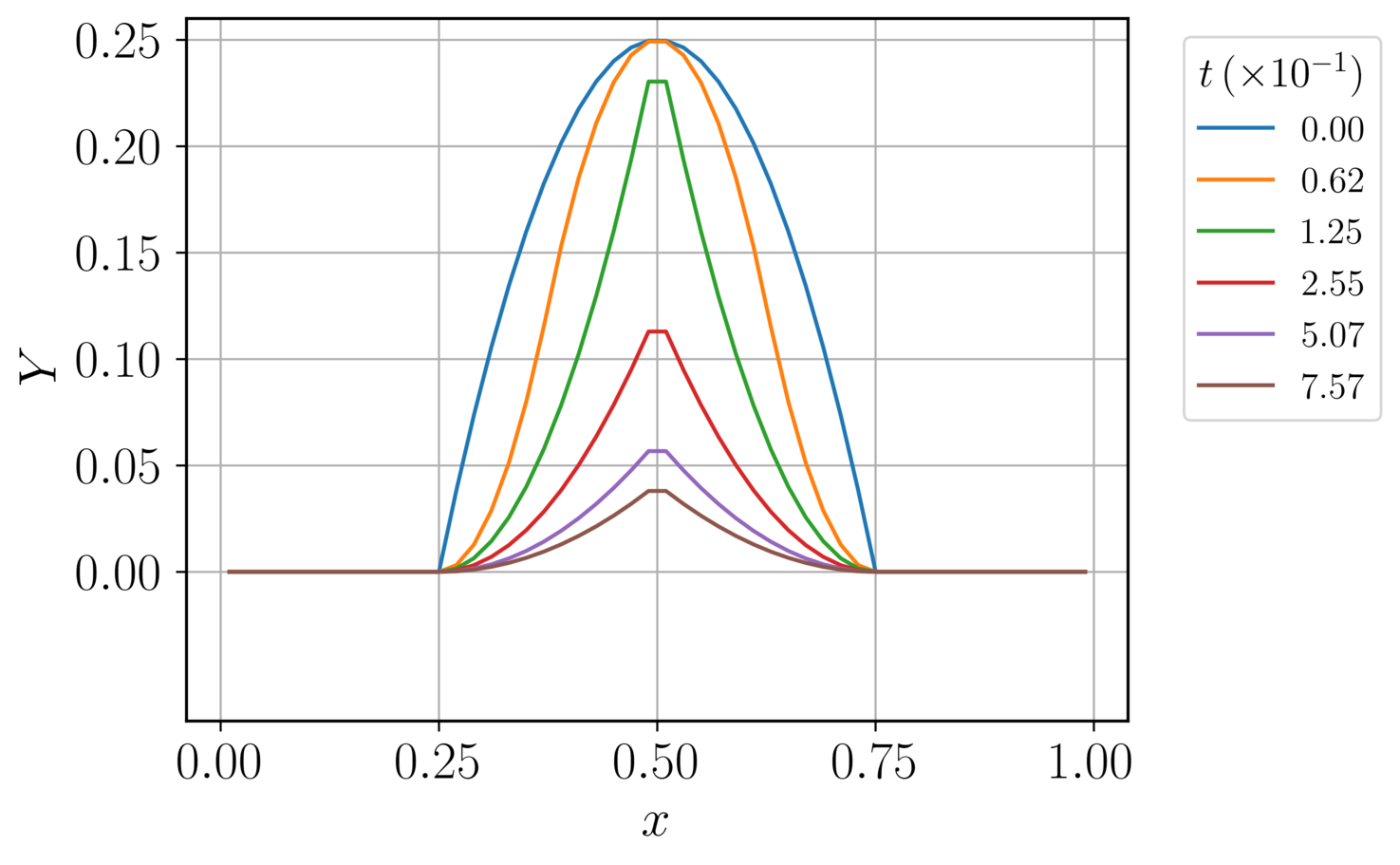}  
  \caption{Exact time plots for $Y(x,t)$}
  \label{fig:sub-burHJ_p5_1}
\end{subfigure}
\begin{subfigure}[t]{.49\textwidth}
  \centering
  \includegraphics[width=0.95\linewidth]{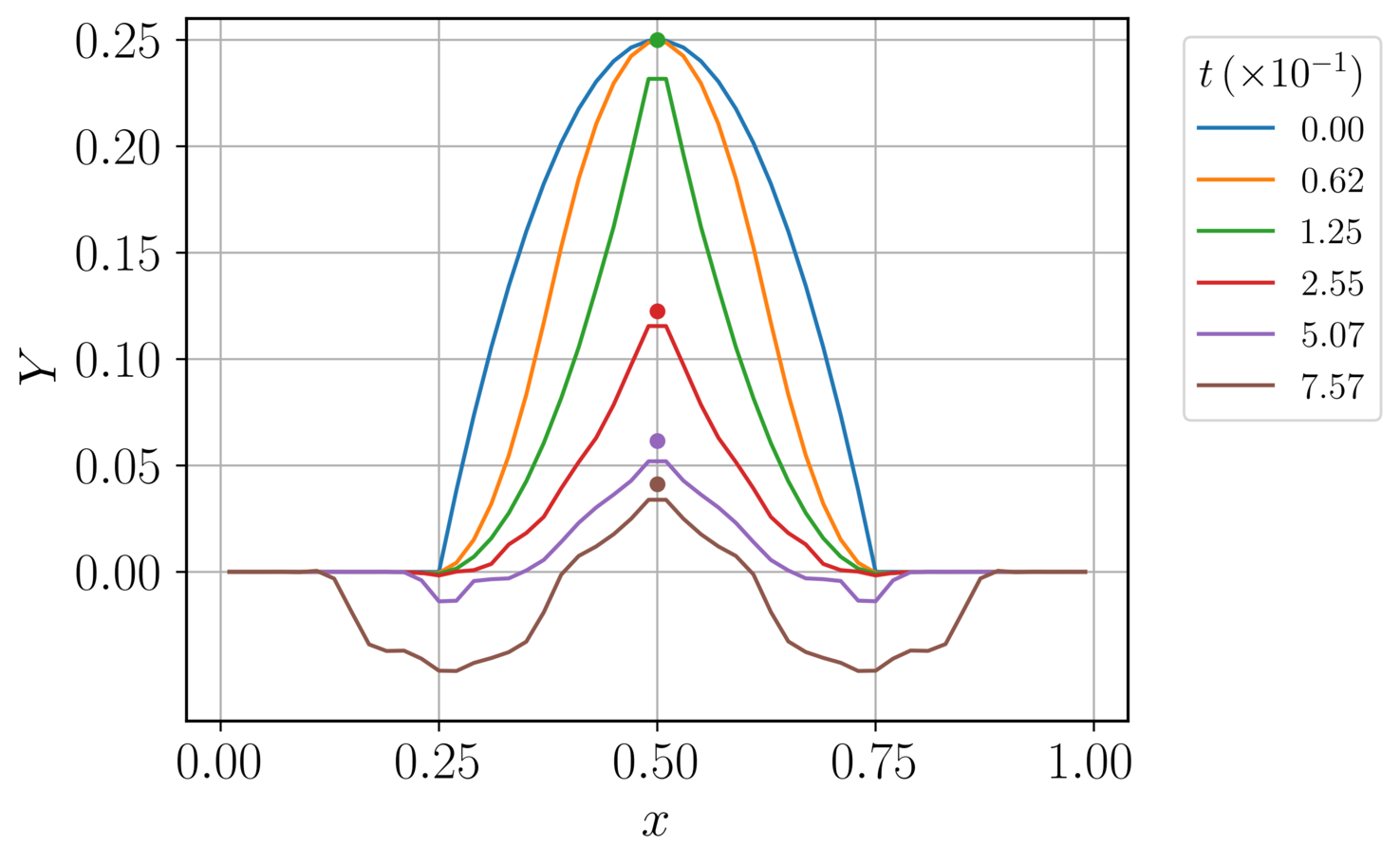}  
  \caption{Line plots for $Y(x,t)$ obtained using the dual scheme}
  \label{fig:sub-burHJ_p5_2}
\end{subfigure}
\caption{The two plots have been drawn on same scale. Fig.(b) represents the DtP mapping generated $Y$ profiles. The dots in various colors at $x=0.5$ correspond to the exact value of $Y$ derived from \eqref{eq:Nwave_ex_Y}, corresponding to the time indicated by the color of each respective dot, after the kink formation takes place at $t=0.125$.} 
\label{fig:burHJ_p5}
\end{figure}

As the last example, the following initial and boundary conditions are considered:
\begin{equation*}
\begin{gathered}
\begin{aligned}
Y_0(x) =
\begin{cases}
    0 & \text{for } x < 0.25 \\
    -2h_0(x^2-x) -3h_0/8  & \text{for } 0.25 \leq x < 0.75 \\
    0 & \text{for } x > 0.75;
\end{cases}
\end{aligned} \\
Y_l(t) = 0; \\
h_0 = 2.
\end{gathered}
\end{equation*}
The exact entropy solution for this setup is given by \eqref{eq:Nwave_ex_Y}, which gets approximated well based on the profiles displayed in Fig.~\ref{fig:burHJ_p5}, although we again notice the standing shock formations at $x=0.25$ and $x=0.75$ after a certain time based on Fig.~\ref{fig:sub-burHJ_p5_2}.

\subsection{Viscous Burgers-HJ}
Evident from the examples presented in the last section, the dual scheme does not recover the entropy solutions in the case of the Fan \ref{sec:HJ_fan}, the Half N-wave \ref{sec:HJ_halfwave} and the N-wave \ref{sec:HJ_nwave}.
As mentioned previously, for any Riemann problem of the form \eqref{eq:ini_Burgers_fan} with $U_R > U_L$, the equation has non-unique weak solutions and the entropy condition, i.e.~information along characteristics travel towards the shock-curve in space-time and not away from it, acts to select a solution, in this case a fan. This ambiguity disappears when viscosity is introduced; indeed, entropy solutions for inviscid Burgers(-HJ) are known to be limits of Burgers solutions as $\nu \to 0$. While our scheme recovers entropy solutions when dealing with the Burgers equation examples in Sec. \ref{sec:Burgers_examples}, this does not happen with the computational scheme for Burgers-HJ, and it appears to recover a greater variety of weak solutions (which, in a sense, is what the scheme, in and of itself, is required to do without any further conditions). The  examples presented in the last section demonstrate this point.

Considering the functional \eqref{eq:HJ_functional_explicit}, it is evident that the derivative $\p_x \lambda$ is not being penalized. This allows for  spatial discontinuities to exist in $\lambda$ (that may arise from even inevitable numerical round-off errors in $\lambda$), and correspondingly in $\hat{u}$ due to the DtP mapping \eqref{eq:DtP_HJ2}. These discontinuous will be propagated in the spatial domain due to the inherent characteristics of the inviscid Burgers equation. To test this hypothesis, we work with a dual formulation of Burgers-HJ \eqref{eq:HJ_original} with $\nu\neq0$.

The first order system associated with \eqref{eq:HJ_original} can be written as
\begin{subequations}
\begin{gather}
    \p_t Y = -\frac{u^2}{2} + \nu \p_x u  \quad \text{in }\Omega; \label{eq:HJ_primal_nu_1}\\
    \p_x Y = u \qquad\text{in }\Omega\label{eq:HJ_primal_nu_2}.
\end{gather}    \label{eq:HJ_primal_nu}%
\end{subequations}
The formulation for this set of equations is analogous to the formulation presented in Sec.~\ref{sec:HJ_formulation}. Consider the following pre-dual functional given by 
\begin{multline}
\widehat{S}_H[\,Y,u,\lambda,\gamma\,] 
    =
    \int_\Omega \,\mathcal{L}(Y,u,\mathcal{D},x,t) \,dt\,dx
     -\left.\int_0^T \,f_i(t)\,\lambda(x,t)\,dt\right|_{x=0}^{x=L} \\ - \int_0^L \,Y_0(x)\,\lambda(x,0)\,dx
    - \int_0^T \,Y_l(t)\,\gamma(0,t)\,dt; \label{eq:HJ_functional_nu}    \end{multline}
    \begin{equation*}
        \mathcal{L}(Y,u,\mathcal{D},x,t) = -Y\,\p_t \lambda \,  + \frac{u\,^2 \lambda}{2} 
    -Y\,\p_x\gamma + \nu \,u\,\p_x\lambda   - u\,\gamma +\frac{\beta_Y}{2}(Y-\bar{Y})^2+\frac{\beta_u}{2}(u-\bar{u})^2,
    \end{equation*}
where $f_i$, $i=1,2$, are prescribed functions on the left and right boundaries defining natural boundary conditions on those boundaries when $\lambda$ is allowed to freely vary there.
Employing the following conditions yields the DtP map
\begin{subequations}
\begin{gather}
\frac{\p \mathcal{L}}{\p Y} = 0:\quad Y=Y^{(H)}(\mathcal{D},x,t) = \bar{Y}+\frac{\p_t \lambda + \p_x \gamma}{\beta_Y};\label{eq:DtP_HJ1_nu}
\\
\frac{\p \mathcal{L}}{\p u} = 0:\quad
u=u^{(H)}(\mathcal{D},x,t) = \bar{u}+\frac{\gamma - \lambda \bar{u} - \nu\, \p_x \lambda}{\beta_u + \lambda},
\label{eq:DtP_HJ2_nu}
\end{gather}
\label{eq:DtP_HJ_nu}%
\end{subequations}
where again $\dee = (D, \nabla D)$.

Substituting the DtP mapping back into  \eqref{eq:HJ_functional_nu} gives us the functional $S_H[\lambda,\mu]$ which can be explicitly written as
\begin{multline*} \label{eq:HJ_nu_functional_explicit}
S_H[\lambda,\gamma] = \int_\Omega \,\left( -\frac{\mathbb{K}_1}{\beta_Y}\,(\p_t \lambda \,+\, \p_x \gamma)^2 -  \frac{\mathbb{K}_2}{\beta_u}\,(\nu\,\p_x \lambda+\bar{u}\lambda - \gamma)^2\right)\,dt\,dx \\
+ \bigintsss_\Omega \,\bigg(-\bar{Y}(\p_t \lambda + \p_x \gamma) + \bar{u}\left(\frac{\lambda\bar{u}}{2} - \gamma + \nu\, \p_x \lambda\right)\bigg)\,dt\,dx \\- \int_0^L \,Y_0(x)\,\lambda(x,0)\,dx
    - \int_0^T \,Y_l(t)\,\gamma(0,t)\,dt,
\end{multline*}
where 
$$\mathbb{K}_1 = \frac{1}{2}; \qquad \mathbb{K}_2 = \frac{1}{2\left(1+\frac{\lambda}{\beta_u}\right)},$$
and for $\beta_u\gg |\lambda|$, $\p_x \lambda$ now gets penalized (\textit{multiplied by $\nu$}). Of course, the E-L equations of $S_H$ is \eqref{eq:HJ_primal_nu} with the standard substitution of $(Y,u) \to (\hat{Y},\hat{u})$.

Considering the ellipticity of E-L equations of $S_H$, i.e.~\eqref{eq:HJ_primal_nu} with the above replacement, the matrix in \eqref{eq:HJ_ellipticity_matrix} now has a modification: $$\mathbb{A}_{1212}=\frac{\nu^2}{\beta_u+ \lambda},$$ with the other terms unchanged.
Correspondingly, we have
\begin{equation*}
    c_i A_{ij} c_j = \frac{(c_1 + c_4)^2}{\beta_Y} + \frac{c_2^2 \, \nu^2}{\beta_u + \lambda}= \frac{(c_1 + c_4)^2}{\beta_Y} + \frac{c_2^2 \, \nu^2}{\beta_u\left(1 + \frac{\lambda}{\beta_u}\right)}\qquad \forall c\in \mathbb{R}^4.
\end{equation*}
 In case when $\lambda=0$, we always have $c_iA_{ij}c_j\geq0$ for non-trivial $c$. \udkk{ Let $\mathcal{N}$ represent a neighborhood around $\mathcal{D} = 0$ given by:
$$\mathcal{N} = \left\{a \in \R^{6}: || a||_{\R^6} < \beta_u \right\}.$$ 
(which is conservative). Thus for $\mathcal{D}\in \mathcal{N}$,
 
  $$c_iA_{ij}c_j\geq0 \quad \forall c\in \mathbb{R}^4.$$ 
  
   This establishes that $\mathbb{A}_{ijk\ell}$ is positive semi-definite in this neighborhood and the equation set \eqref{eq:HJ_primal_nu} is locally degenerate elliptic.} In contrast to the inviscid Burgers-HJ, the only rank-one matrices along which $C : \mathbb{A} C = 0$ are of the form $(0,a_2)\otimes(1,0)$ for any $a_2\in\mathbb{R}$. Based on the discussion in Sec.~\ref{sec:HJ_formulation}, it follows that ellipticity is lost only along the normal $(1,0)$, where the first argument represents the component of the normal along the time direction. This is consistent with solutions of Burgers equation $(\nu \neq 0)$ as spatial discontinuities in $u$ are not allowed. The only discontinuity allowed is in $\p_t \gamma$ across space-like surfaces ($t =$ constant) but the DtP mapping equations \eqref{eq:DtP_HJ_nu} shows that this does not affect $u^{(H)}$ and hence the mapped primal solutions are continuous; a similar conclusion was also obtained in our dual formulation for the heat equation (see-\cite[Sec.~5]{KA1}).

Similar to the inviscid case, we obtain the following side conditions:
\begin{equation}
\label{eq:HJ_side_cond}
    \hat{Y}(x,0) = Y_0(x); \qquad \hat{Y}(0,t) = Y_l(t)
\end{equation}
where $\hat{Y} = Y^{(H)}(\mathcal{D}(x,t),x,t)$
. Additionally, on each of the boundaries, the following operation can be performed which lead to different classes of boundary conditions and hence different problem setups:
\begin{itemize}
    \item  $\lambda$ can be specified as a Dirichlet BC with a fixed arbitrary value. In such a case, $\delta \lambda=0$ on the boundary under consideration and we obtain \eqref{eq:HJ_side_cond} as the only natural b.cs.
    \item If $\lambda$ is not specified and allowed to vary freely, then an extra natural BC for $S_H$,
    $$\nu \,\hat{u} = f_i$$ emerges, corresponding to the boundary under consideration.
\end{itemize}
Our focus is primarily on observing the effect of the $\nu$ term in the bulk of the domain. 
Accordingly, we consider the second BC above, weakly enforcing $u=0$ on both ends. The impact of this term on the boundaries becomes more pronounced for higher values of $\nu$, and we want to operate with as small a value of $\nu$ as possible. In the following, we show calculations for $\nu=10^{-3}$ for the Fan, the Half N-wave and the N-wave problems. 
  
 The residual \eqref{eq:weak_dual_HJ} and Jacobian \eqref{eq:jaco_HJ} and their corresponding discrete versions are now appropriately modified.

\subsubsection{Fan}
\begin{figure}
\begin{subfigure}[t]{.3\textwidth}
  \centering
  \includegraphics[width=\linewidth]{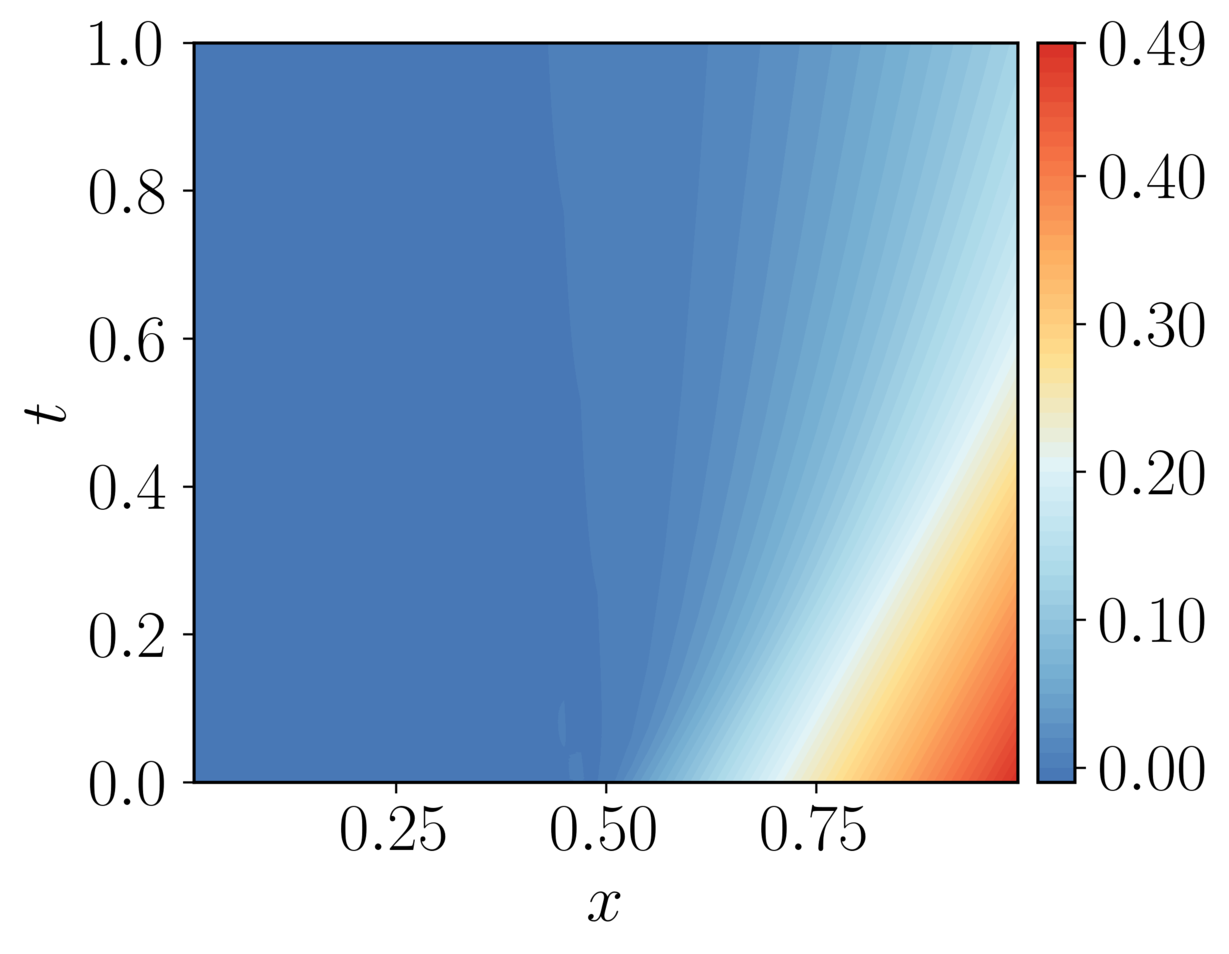}  
  \caption{$Y(x,t)$}
  \label{fig:sub-burHJ_p1_3}
\end{subfigure}
\begin{subfigure}[t]{.3\textwidth}
  \centering
  \includegraphics[width=0.98\linewidth]{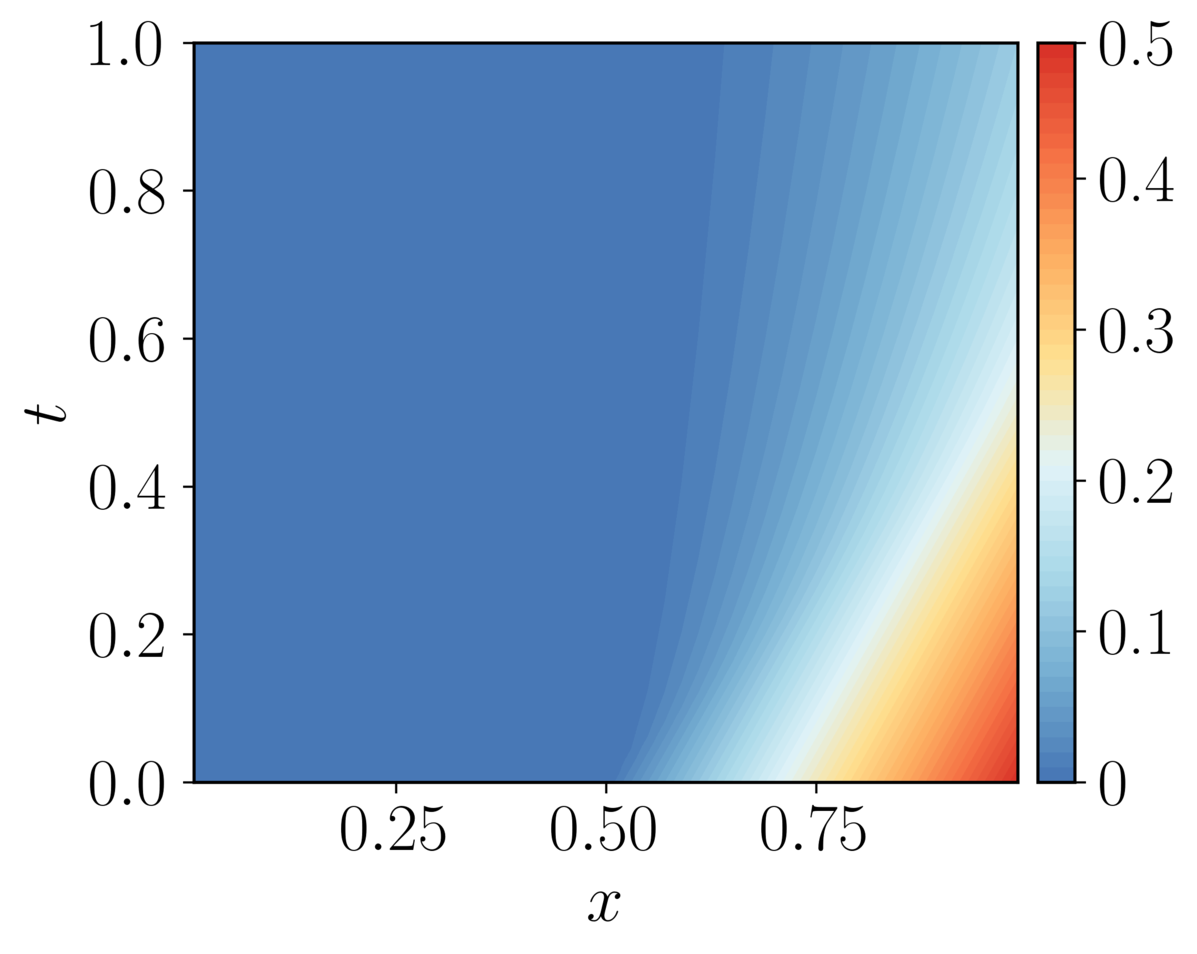}  
  \caption{Exact $Y(x,t)$ based on \eqref{eq:fan_ex_Y}}
  \label{fig:sub-burHJ_p1_6}
\end{subfigure}
\begin{subfigure}[t]{.39\textwidth}
  \centering
\includegraphics[width=0.95\linewidth]{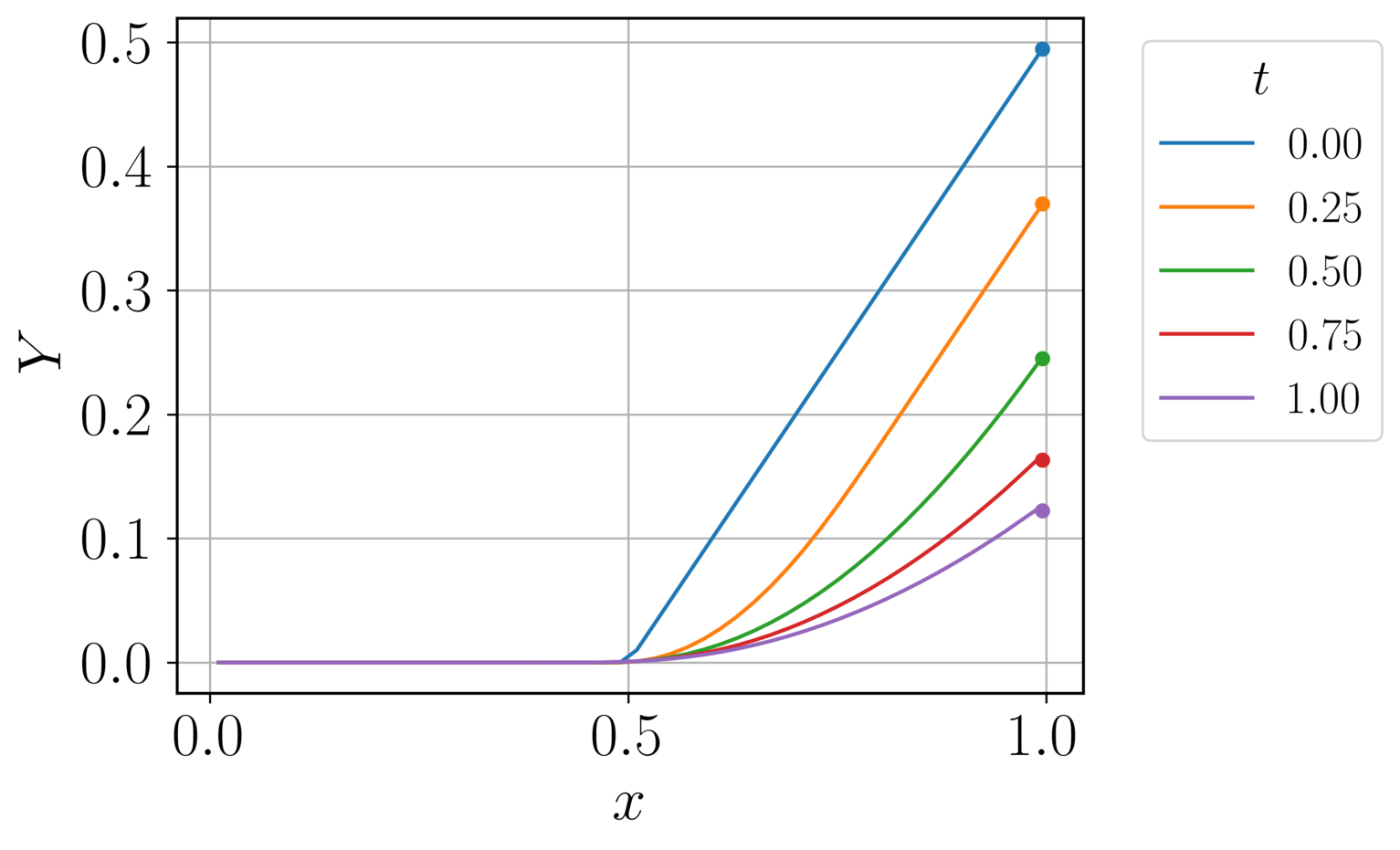}  
  \caption{Line plots for $Y(x,t)$}
  \label{fig:sub-burHJ_p1_4}
\end{subfigure}
\caption{DtP mapping generated primal field $Y$ for the fan example with $\nu=10^{-3}$. In Fig.(c), the dots in various colors at last spatial Gauss point correspond to the exact value of $Y$ derived from \eqref{eq:fan_ex_Y}, corresponding to the time indicated by the color of each respective dot.} 
\label{fig:burHJ_p1_nu}
\end{figure}

Considering \ref{sec:HJ_fan} with an active $\nu$, the outcomes are presented in Fig.~\ref{fig:burHJ_p1_nu}. It is apparent from the figure that the results have exhibited enhancement and closely resemble the exact outcomes of the inviscid case. Upon close inspection, a slight deviation from the exact expression can be noticed close to the right boundary. This discrepancy can be attributed to the influence of an active $\nu$ and more specifically, to the natural imposition of $\nu u = 0$ on that boundary.

\subsubsection{Half N-wave}
\begin{figure}
\begin{subfigure}[t]{.3\textwidth}
  \centering
  \includegraphics[width=\linewidth]{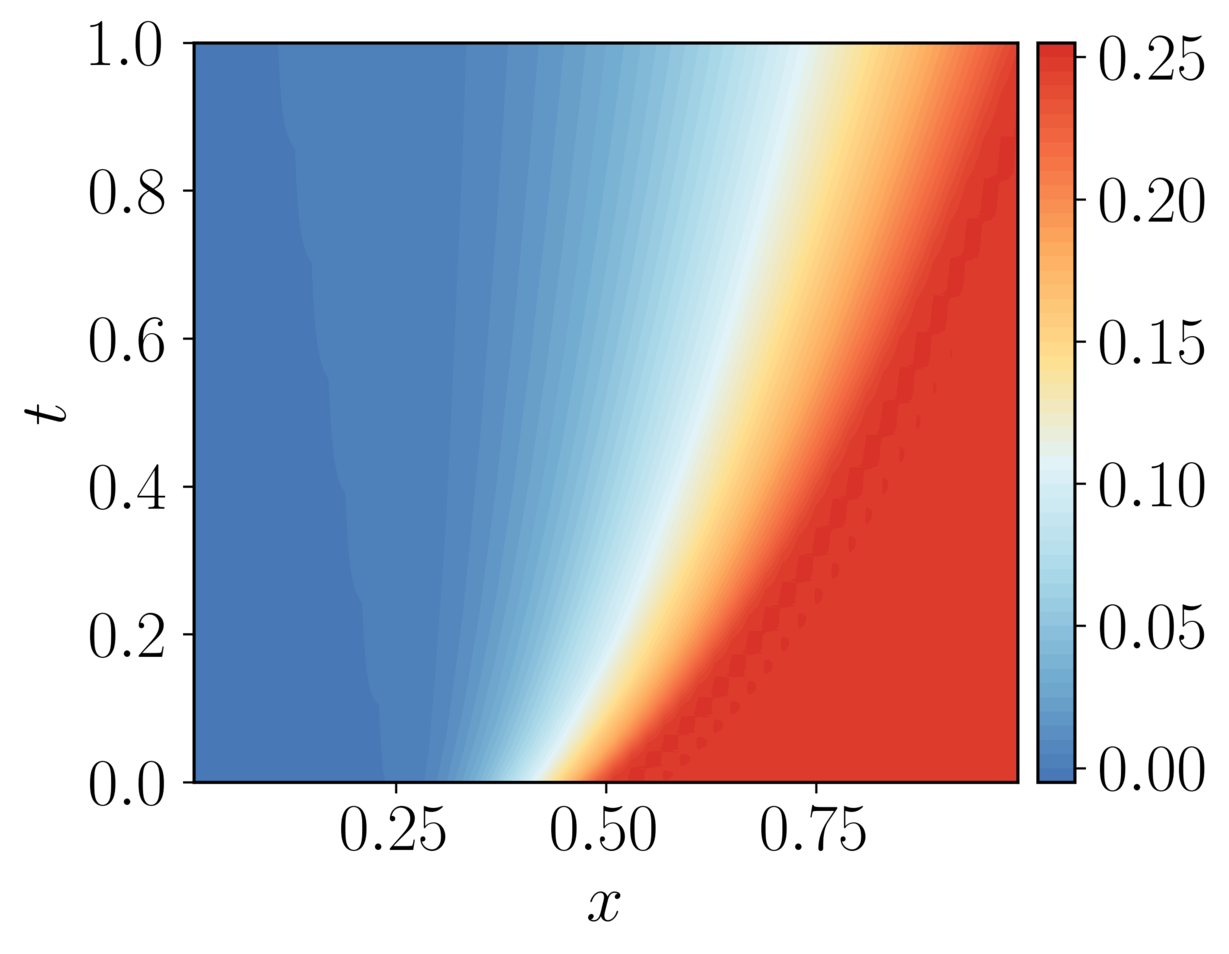}  
  \caption{$Y(x,t)$}
  \label{fig:sub-burHJ_p4_3}
\end{subfigure}
\begin{subfigure}[t]{.3\textwidth}
  \centering
  \includegraphics[width=\linewidth]{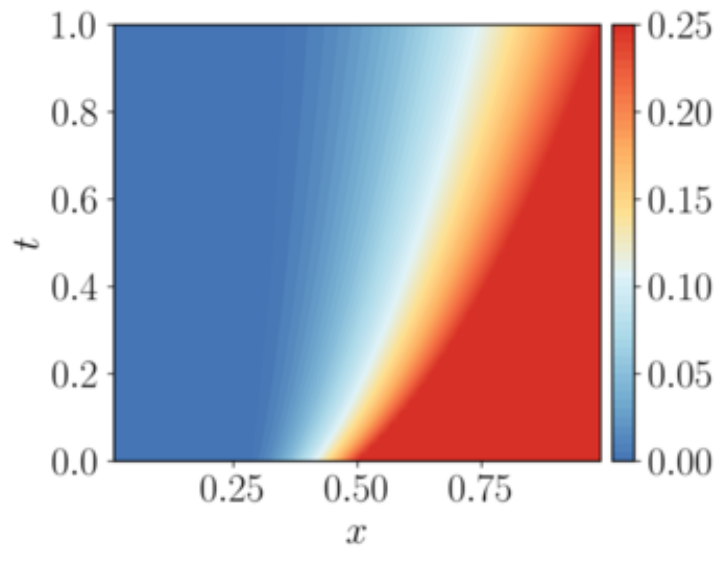}  
  \caption{Exact $Y(x,t)$ based on \eqref{eq:halfwave_ex_Y}}
  \label{fig:sub-burHJ_p4_6}
\end{subfigure}
\begin{subfigure}[t]{.39\textwidth}
  \centering
\includegraphics[width=0.95\linewidth]{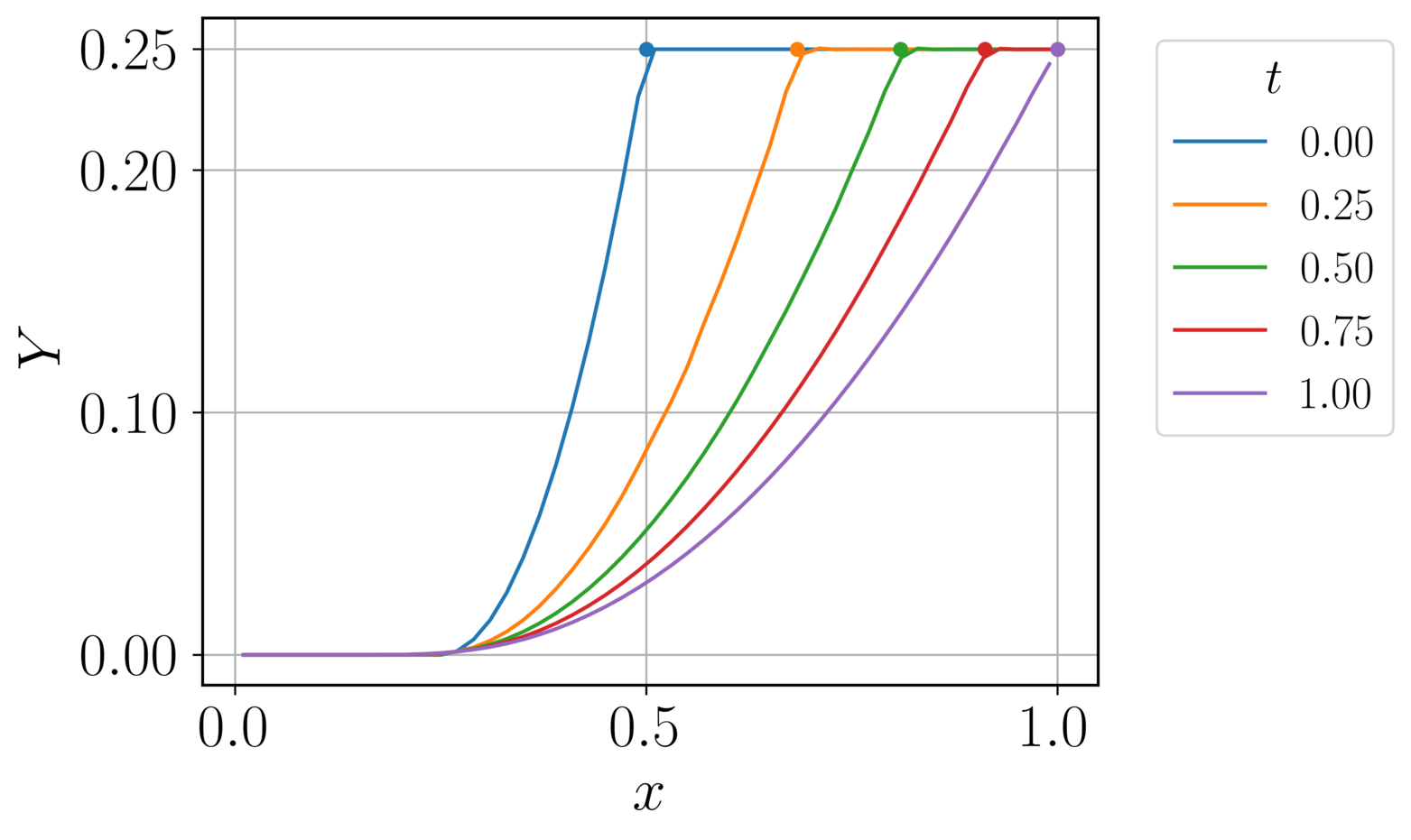}  
  \caption{Line plots for $Y(x,t)$}
  \label{fig:sub-burHJ_p4_4}
\end{subfigure}
\caption{DtP mapping generated primal field $Y$ for the Half N-wave with $\nu=10^{-3}$. In Fig.(c), the coloured dots at $Y=0.25$ represent the trajectory of the kink in the space-time domain based on the exact expressions \eqref{eq:DShock_ex_Y} for $\nu=0$.} 
\label{fig:burHJ_p4_nu}
\end{figure}

Considering \ref{sec:HJ_halfwave} again with an active $\nu$ term, the results as shown in Fig.~\ref{fig:burHJ_p4_nu}. The boundary effects in this problem are not very pronounced because the Half N-wave itself has a value of $u=0$ on the right boundary for a significant duration based on the exact expression. With a small viscosity $\nu$, one would not anticipate the solution for this problem to be significantly different from the case where $\nu=0$, which is also evident from the results presented.

\subsubsection{N-wave}
\begin{figure}[h]
    \begin{minipage}{.5\linewidth}
        \centering
        \subfloat[$Y(x,t)$]{\label{fig:sub-burHJ_nu_p5_1}\includegraphics[width=0.8\linewidth]{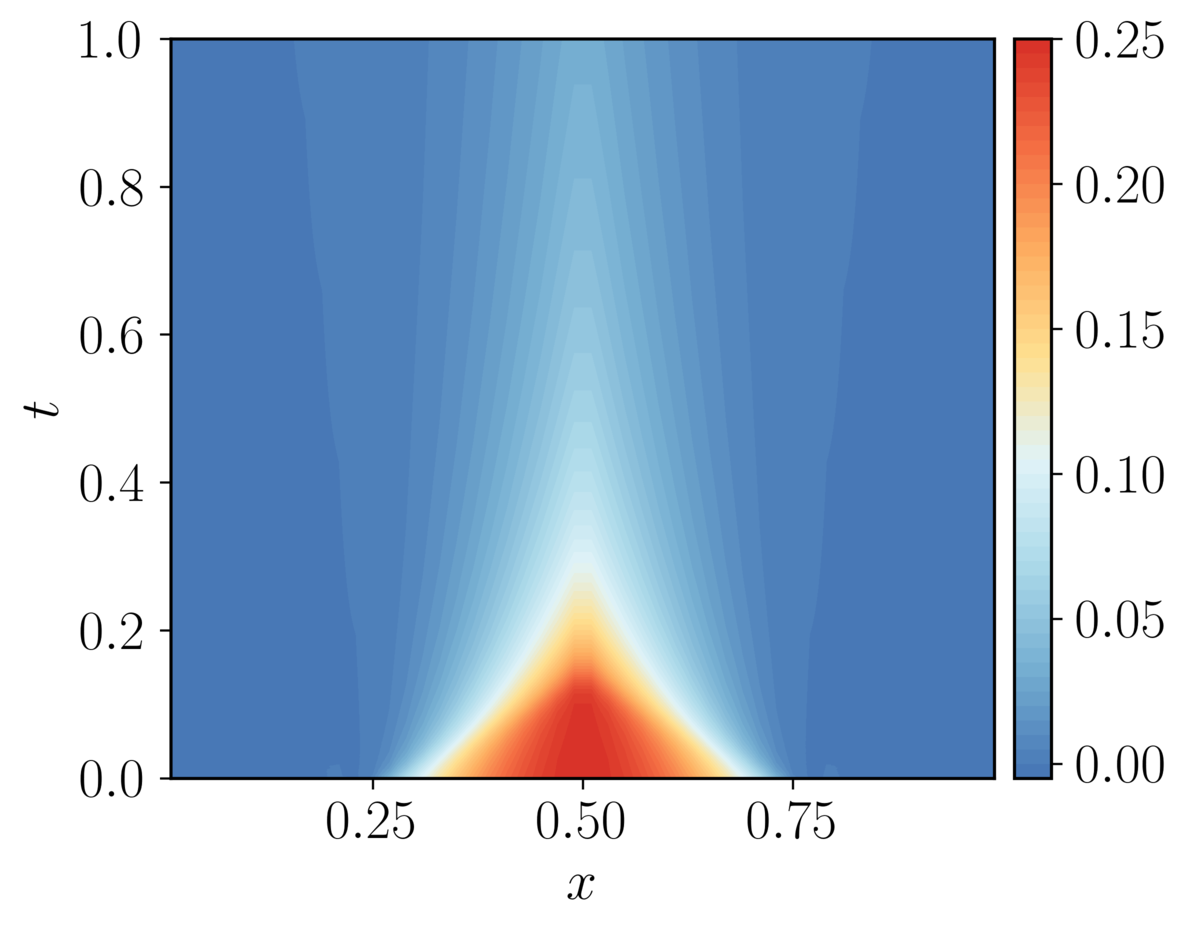}}
    \end{minipage}%
    \begin{minipage}{.5\linewidth}
        \centering
        \subfloat[Exact $Y(x,t)$ based on \eqref{eq:Nwave_ex_Y}]{\label{fig:sub-burHJ_nu_p5_2}\includegraphics[width=0.8\linewidth]{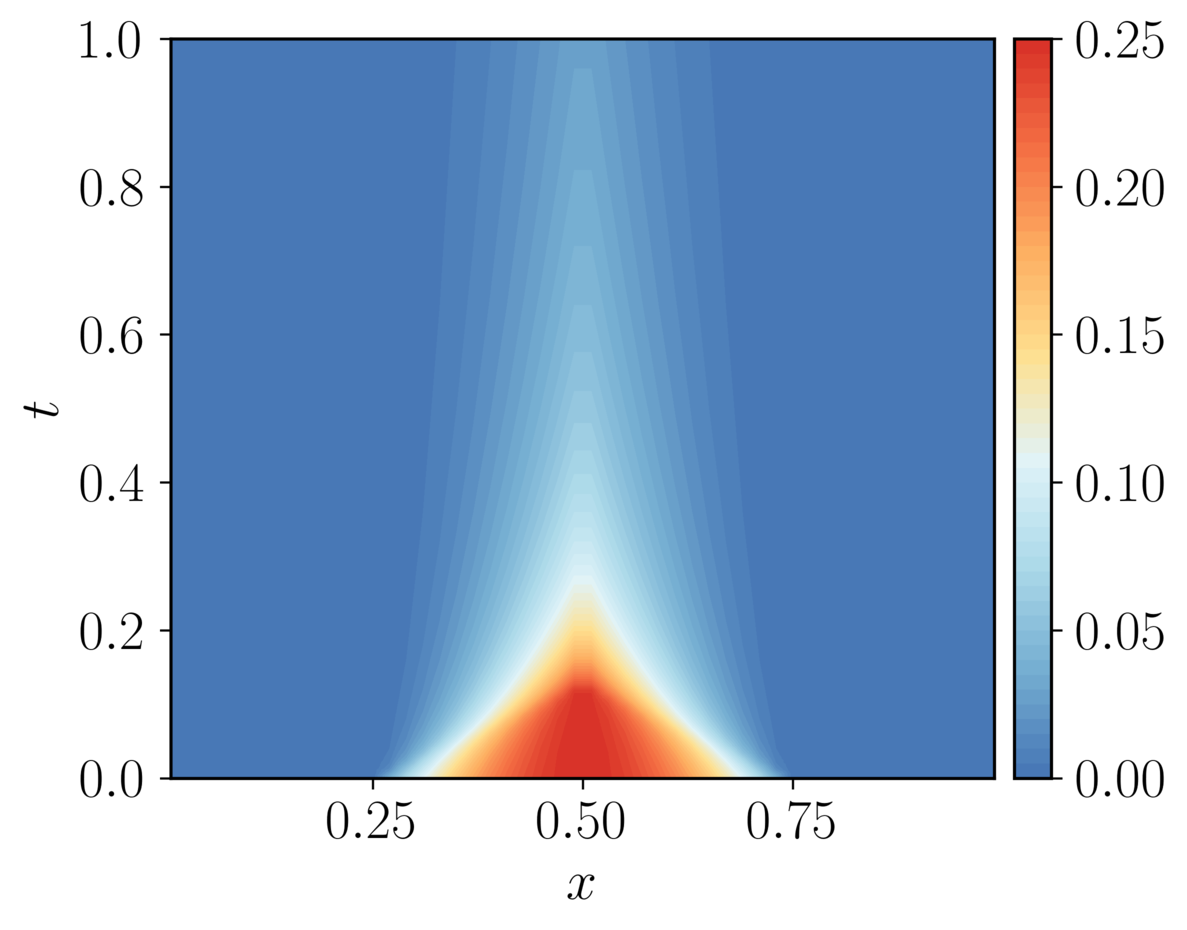}}
    \end{minipage}\par\medskip
    
    \begin{minipage}{.5\linewidth}
        \centering
        \subfloat[Line plots for $Y(x,t)$]{\label{fig:sub-burHJ_nu_p5_3}\includegraphics[width=\linewidth]{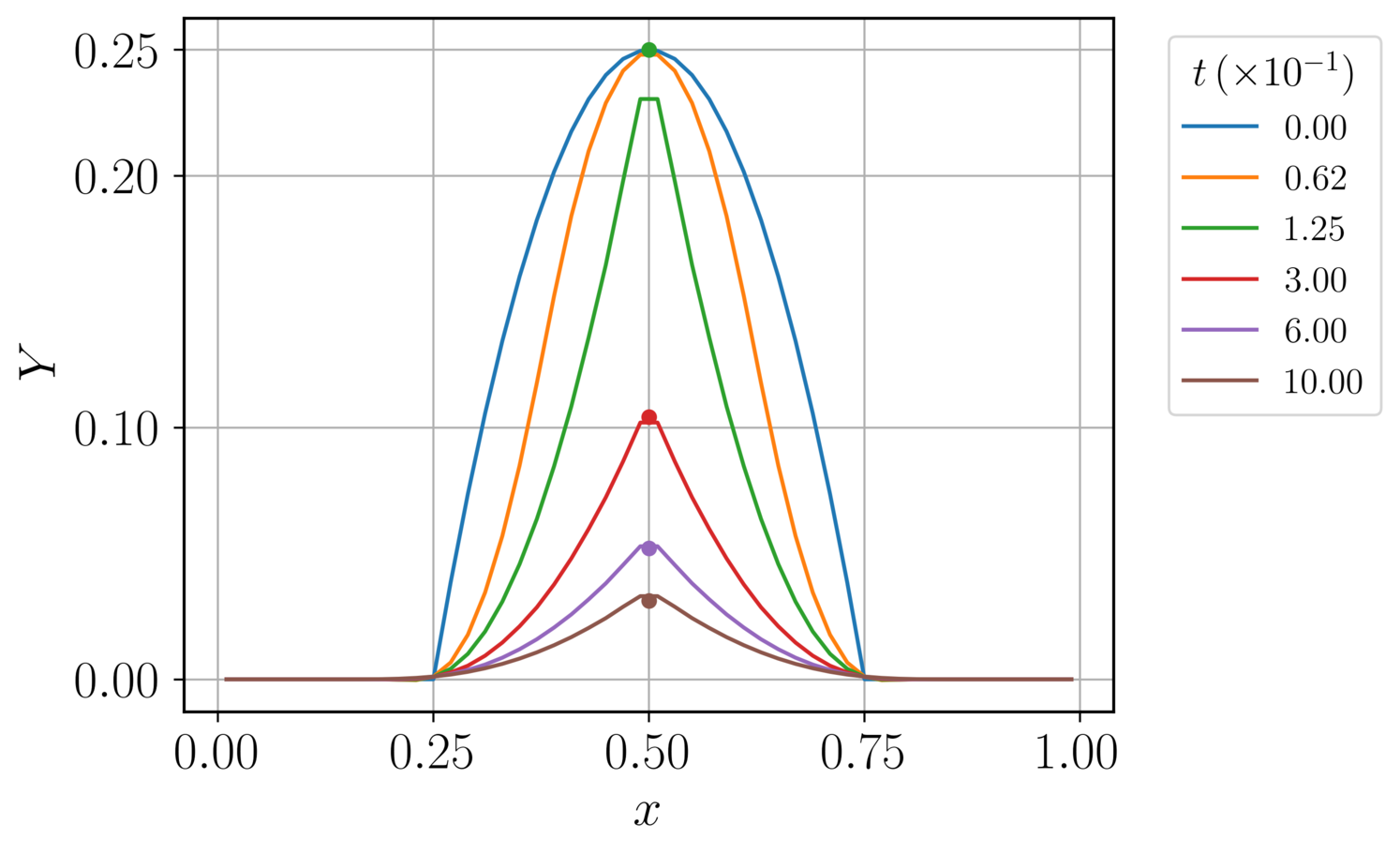}}
    \end{minipage}%
    \begin{minipage}{.5\linewidth}
        \centering
        \subfloat[Line plots for $L^2$ projected $Y(x,t)$]{\label{fig:sub-burHJ_nu_p5_4}\includegraphics[width=\linewidth]{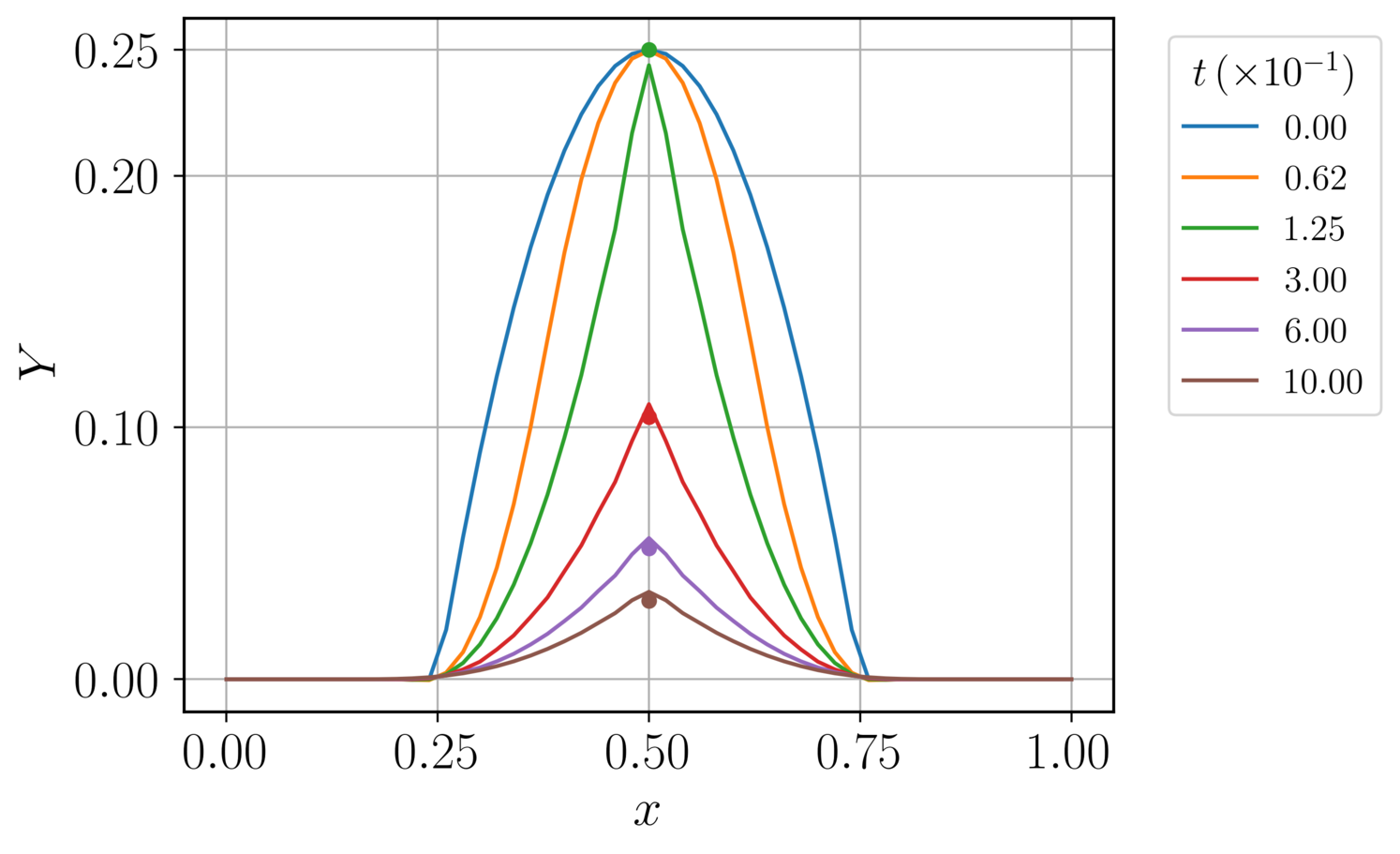}}
    \end{minipage}

    \caption{DtP mapping generated primal field $Y$ for the N-wave with $\nu=10^{-3}$. In Fig.~\ref{fig:sub-burHJ_nu_p5_3} and Fig.~\ref{fig:sub-burHJ_nu_p5_4}, the colored dots at $X=0.5$ represent the trajectory of the kink in the space-time domain based on the exact expressions \eqref{eq:Nwave_ex_Y} for $\nu=0$.}
    \label{fig:burHJ_p5_nu}
\end{figure}

Considering the problem \ref{sec:HJ_nwave} again with an active $\nu$, the results are presented in \ref{fig:burHJ_p5_nu}, from which it is evident that the two dips have been resolved. Furthermore, the line plots for $Y$ derived from the $L^2$ projected data at different nodal times are displayed, revealing the existence of a kink at the node located at $x=0.5$. A slight deviation in the speed at which the kink travels down from the exact expression may again be attributed to the presence of an active $\nu$.

\subsection{Inviscid Burgers-HJ with Viscous base states}

In this section we will explore solutions to the inviscid Burgers-HJ equation generated from the dual scheme with base states provided from exact solutions to the (viscous) Burgers-HJ equation with small viscosity.

Given an initial condition $u_0(x)$ for the Burgers equation, for $x\in\mathbb{R}$, the Hopf-Cole transformation can be used to derive the exact expression for $u(x,t)$ \cite{TP_Liu} 
along with its antiderivative with respect to $x$, $Y(x,t)$, such that  for $t>0$ and $x\in\mathbb{R}$,
\begin{subequations}
\begin{gather}
Y(x,t) = - 2\nu \log\left[\frac{1}{2\sqrt{\pi\nu t}}\int_{-\infty}^{\infty}\,e^{-\frac{Y_0(y)}{2\nu}-\frac{(x-y)^2}{4\nu t}} \,dy\right] ;\label{eq:vis_Y}\\
    u(x,t) = \frac{\int_{-\infty}^{\infty}\frac{x-y}{t}e^{-\frac{(x-y)^2}{4\nu t}-\frac{Y_0(x)}{2\nu}}\,dy}{\int_{-\infty}^{\infty}e^{-\frac{(x-y)^2}{4\nu t}-\frac{Y_0(x)}{2\nu}}dy},\label{eq:vis_u}
\end{gather}
\label{eq:vis_formulas}
\end{subequations}
  where $Y_0$ can be generated via
  \udkk{
  $$Y_0(x) = \int_0^x \,u_0(y)\,dy,$$
  }and we set $Y_0(0)=0$ without loss of generality. These expressions will be denoted as the `viscous' formulae.
  
\begin{figure}[h!]
\begin{subfigure}[t]{.49\textwidth}
  \centering
 \includegraphics[width=0.95\linewidth]{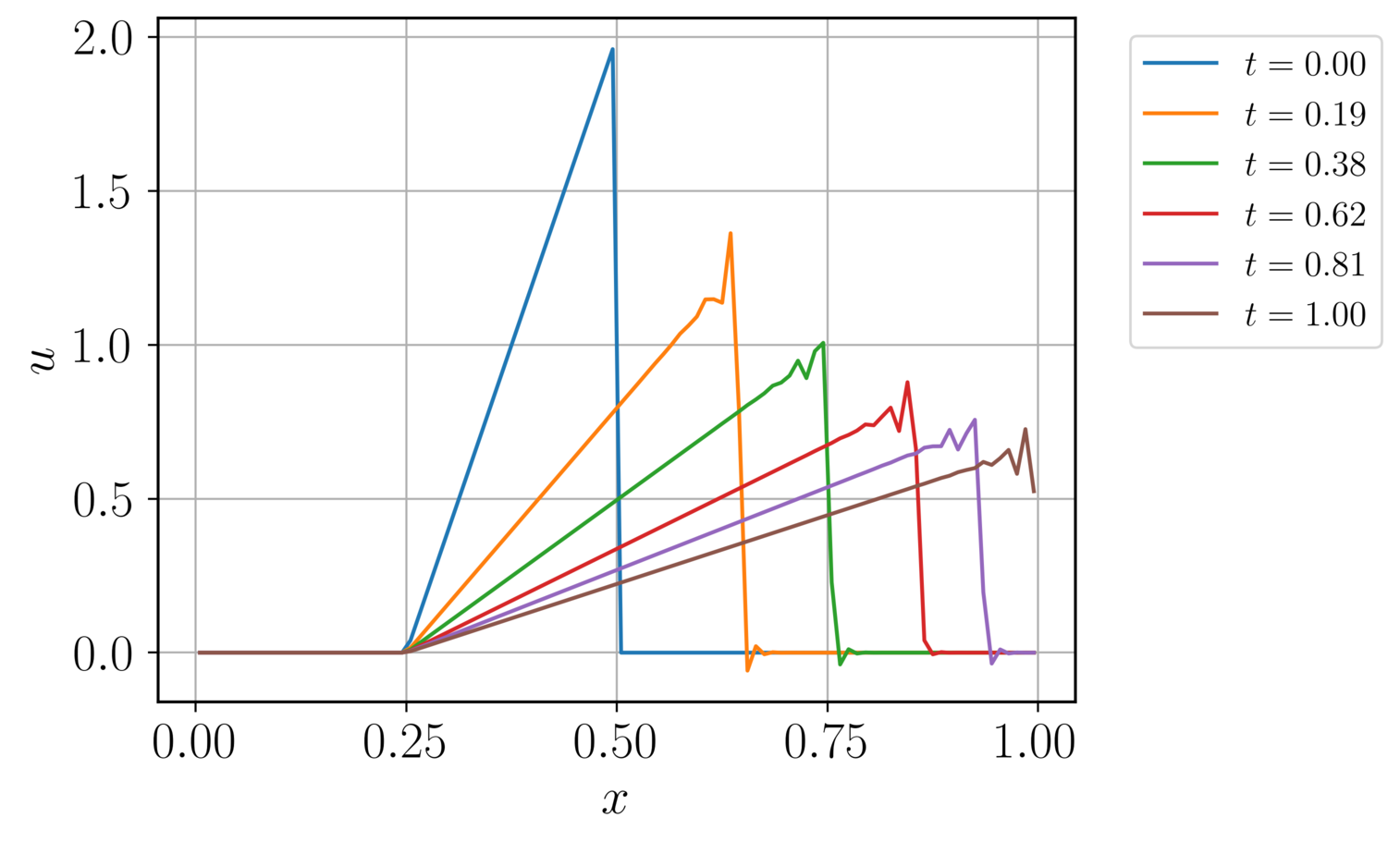}  
  \caption{$u(x,t)$}
  \label{fig:vis_base_1}
\end{subfigure}
\begin{subfigure}[t]{.49\textwidth}
  \centering
 \includegraphics[width=0.95\linewidth]{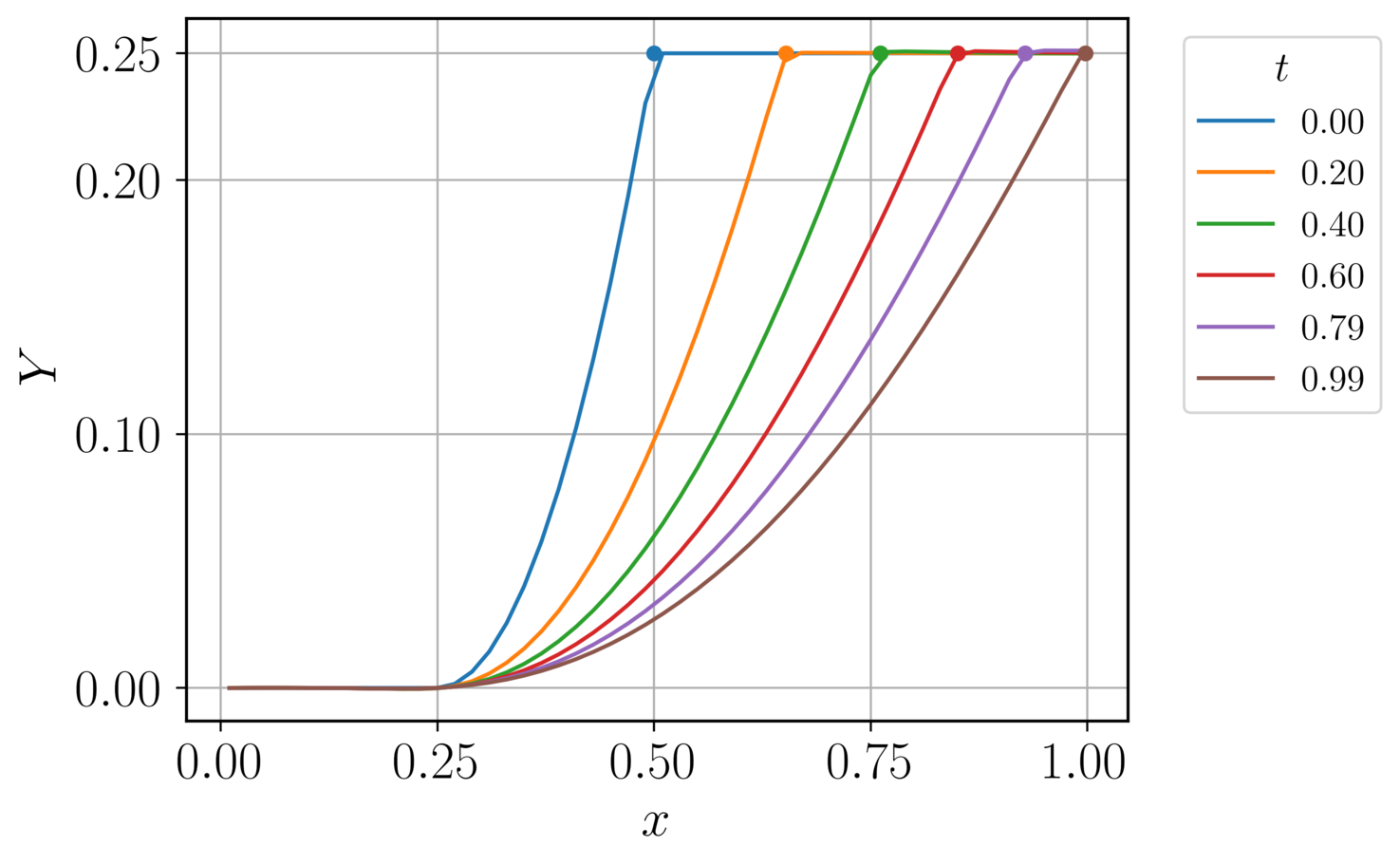}  
  \caption{$Y(x,t)$}
  \label{fig:vis_base_2}
\end{subfigure}
\caption{(a) is obtained from the formulation based on the Burgers equation Sec.~\ref{sec:HJ_examples}. (b) is produced using the formulation based on Burgers-HJ Sec.~\ref{sec:Burgers_examples}. Respective base states used for both the problems are based on the viscous formulas \eqref{eq:vis_formulas}. }
\label{fig:vis_base}
\end{figure}

For $\nu \rightarrow 0$, the formula \eqref{eq:vis_formulas} for $u(x,t)$ recovers the solution to the inviscid Burgers equation. Motivated by this fact, we set a small value $\nu=10^{-3}$ and generate the solutions to the Burgers equation using the expressions \eqref{eq:vis_formulas}. 
We consider the algorithm presented for the Burgers equation (Sec.~\ref{sec:Burgers_examples}) and use the same settings as used for the Half N-wave (Sec.~\ref{sec:Burgers_halfwave}), except that the base state for $u$ is now set based  the expression obtained using the formula  \eqref{eq:vis_u}. The result for such a setup is shown in Fig.~\ref{fig:vis_base_1}.
Similarly for the Burgers-HJ, we utilize the algorithm on the Half N-wave (Sec.~\ref{sec:HJ_halfwave}), with the base states for $Y$ and $u$ obtained from the formulae \ref{eq:vis_Y} and \ref{eq:vis_u}, respectively. We use coarser meshes as compared to \ref{sec:HJ_halfwave} and set $N_x=50$, $N_t=400$ and $T^{(s)} =0.01$ (fine meshes work as well). The result for this setup has been shown in Fig.~\ref{fig:vis_base_2}. 

Based on the results we observe the following:
\begin{itemize}
    \item In the case of the inviscid Burgers equation, the unique entropy solution is again recovered, now with no smoothing operator for base state generation (which is to be expected).
    \item In the case of the dual \textit{inviscid} Burgers-HJ, using the viscous base states is one way of enforcing the entropy condition.
\end{itemize}

\section*{Acknowledgment}
This work was supported by the grant NSF OIA-DMR \#2021019 and was also supported by a Simons Pivot Fellowship to AA. It was completed while he was on sabbatical leave at a) the Max Planck Institute for Mathematics in the Sciences in Leipzig, and b) the Hausdorff Institute for Mathematics at the University of Bonn funded by the Deutsche Forschungsgemeinschaft (DFG, German Research Foundation) under Germany's Excellence Strategy – EXC-2047/1 – 390685813, as part of the Trimester Program on Mathematics for Complex Materials. The support and hospitality of both institutions is acknowledged.

\appendix 
\renewcommand{\thesection}{\Alph{section}}

\section*{Appendices}

\section{Review of the general formalism}\label{app:rev_duality}
The following two subsections are excerpted from \cite{action_2,action_3} to make this paper self-contained.
\subsection{The essential idea: An optimization problem for an algebraic system of equations}\label{sec:fin_dim}
Consider a generally nonlinear system of algebraic equations in the variables $x \in \mathbb{R}^n$ given by
	\begin{equation}\label{eq:alg_sys}
	    A_\alpha (x) = 0,
	\end{equation}
	where $A: \mathbb{R}^n \to \mathbb{R}^N$ is a given function (a simple example would be $A_\alpha (x) = \bar{A}_{\alpha i} \, x^i - b_\alpha$, \udkk{$\alpha = 1,2, \ldots N, i = 1,2, \ldots n$}, where $\bar{A}$ is a constant matrix, \textit{not necessarily symmetric} (when $n = N$), and $b$ is a constant vector). \udkk{For $n>0$ and $N>0$, we allow for all the possibilities such that $n>N$, $n=N$ or $n<N$.}
	
	The goal is to construct an objective function whose critical points solve the system \eqref{eq:alg_sys} (when a solution exists) by defining an appropriate $x^* \in \mathbb{R}^n$ satisfying  $A_\alpha (x^*) = 0$.
	
	For this, consider first the auxiliary function
	\begin{equation*}
	    \widehat{S}_H(x,z) = z^\alpha A_\alpha (x) + H(x)
	\end{equation*}
	(where $H$ belongs to a class of scalar-valued function to be defined shortly) and define
	\begin{equation*}
	    S_H(z) = z^\alpha A_\alpha(x_H (z)) + H(x_H(z))
	\end{equation*}
	with the requirement that the system of equations
	\begin{equation}\label{eq:H_fin_dim}
	    z^\alpha \frac{\p A_\alpha}{\p x^i}(x) + \frac{\p H}{\p x^i}(x) = 0
	\end{equation}
	be solvable for the function $x = x_H(z)$ through the choice of $H$, and \textit{any} function $H$ that facilitates such a solution qualifies for the proposed scheme. 
	
	In other words, given a specific $H$, it should be possible to define a function $x_H(z)$ that satisfies 
	\begin{equation*} 
	z^\alpha \p_{x^i} A_\alpha (x_H(z)) + \p_{x^i} H (x_H(z)) = 0 \quad \forall z \in \mathbb{R}^N
	\end{equation*}
	(the domain of the function $x_H$ may accommodate more intricacies, but for now we stick to the simplest possibility). Note that \eqref{eq:H_fin_dim} is a set of $n$ equations in $n$ unknowns regardless of $N$ ($z$ for this argument is a parameter).
	
	Assuming this is possible, we have
	\begin{equation*}
	    \frac{\p S_H}{\p z^\beta} (z) = A_\beta(x_H(z)) +  \left( z^\alpha \frac{\p A_\alpha}{\p x^i}(x_H(z)) + \frac{\p H}{\p x^i}(x_H(z)) \right) \frac{\p x^i_H}{\p z^\beta}(z) = A_\beta(x_H(z)),
	\end{equation*}
	using \eqref{eq:H_fin_dim}. Thus,
	\begin{itemize}
	    \item if $z_0$ is a critical point of the objective function $S_H$ satisfying $\p_{z^\beta} S_H(z_0) = 0$, then the system $A_\alpha(x) = 0$ has a solution defined by $x_H(z_0)$; 
	    \item if the system $A_\alpha(x) = 0$ has a unique solution, say $y$, and if $z^H_0$ is any critical point of $S_H$, then $x_H\left(z^H_0 \right) = y$, for all admissible $H$.
	    \item If $A_\alpha(x) = 0$ has non-unique solutions, but $\p_{z^\beta} S(z) = 0$ ($N$ equations in $N$ unknowns) has a unique solution for a specific choice of the function $z \mapsto x_H(z)$ related to a choice of $H$, then such a choice of $H$ may be considered a selection criterion for imparting uniqueness to the problem $A_\alpha(x) = 0$.
        \item Finally, to see the difference of this approach with the Least-Squares (LS) Method, we note that the optimality condition for the objective $A_\alpha(x) A_\alpha(x)$ is $A_\alpha(x) \p_{x^i} A_\alpha(x) = 0 \centernot \implies A_\alpha(x) = 0$. 
        
        For a linear system $\bar{A} x = b$, the LS governing equations are given by
        \[
        \bar{A}^T \bar{A} z = \bar{A}^T b,
        \]
        with LS solution defined as $z$ even when the original problem $\bar{A} x = b$ does not have a solution (i.e., when $b$ is not in the column space of $\bar{A}$). The LS problem always has a solution, of course. In contrast, in the present duality-based approach with quadratic $H(x) = \frac{1}{2} x^T x$ the governing equation is
        \[
        - \bar{A}\bar{A}^T z = b
        \]
        with solution to $\bar{A} x = b$ given by $x = - \bar{A}^T z$, and the problem has a solution only when $\bar{A} x = b$ has a solution, since the column spaces of the matrices $\bar{A}$ and $\bar{A}\bar{A}^T$ are identical.
        \end{itemize}
 \subsection{The idea behind the general formalism}
The proposed scheme for generating variational principles for nonlinear PDE systems may be abstracted as follows: We first pose the given system of PDE as a \textit{first-order} system (introducing extra fields representing (higher-order) space and time derivatives of the fields of the given system); as before let us denote this collection of primal fields by $U$. `Multiplying' the primal equations by dual Lagrange multiplier fields, the collection denoted by $D$, adding a function $H(U)$, solely in the variables $U$ (the purpose of which, and associated requirements, will be clear shortly), and integrating by parts over the space-time domain, we form a `mixed' functional in the primal and dual fields given by
\begin{equation*}
    \widehat{S}_H [U,D] = \int_{[0,T]\times \Omega} \ \scl_H (\dee,U) \, dt dx 
\end{equation*}
where \udkk{$\dee := (D, \nabla D)$ and we define 
\[
\delta \dee [D; \delta D] \big|_{(x,t)} := \frac{d}{d \veps} \left( (D + \veps \delta D)\big|_{(x,t)}, \ \nabla (D + \veps \delta D)\big|_{(x,t)} \right) \bigg|_{\veps = 0} = \left(\delta D|_{(x,t)}, \nabla (\delta D)|_{(x,t)} \right),
\]
(and, $A|_{(x,t)} := A(x,t)$ represents the evaluation for any function $(x,t) \mapsto A(x,t)$). For the class of functions on which $\delta \dee$ is evaluated in specific examples, it may not be well-defined (due to $\nabla \delta D$) on sets of vanishing measure in space-time; this typically is not a problem as $\delta \dee$ appears within integrals in products with other functions for which the integrals can be made sense of.} 

We then require that the family of functions $H$ be such that it allows the definition of a function $U_H(\dee)$ such that
\begin{equation*}
    \frac{\p \scl_H}{\p U} (\dee, U_H(\dee)) = 0
\end{equation*}
so that the \textit{dual} functional, defined solely on the space of the dual fields $D$, given by
\begin{equation*}
    S_H[D] = \int_{[0,T]\times \Omega} \scl_H(\dee,U_H(\dee)) \, dt dx 
\end{equation*}
has the first variation
\begin{equation*}
    \delta S_H = \int_{[0,T]\times \Omega} \frac{\p \scl_H}{\p \dee} \delta \dee \, dt dx.
\end{equation*}
By the process of formation of the functional $\widehat{S}_H$, it can then be seen that the (formal) E-L equations arising from $\delta S_H$ have to be the original first-order primal system, with $U$ substituted by $U_H(\dee)$, regardless of the $H$ employed.

Thus, the proposed scheme may be summarized as follows: we wish to pursue the following (local-global) critical point problem
\begin{equation*}
   \begin{smallmatrix} \mbox{extremize}\\ D\end{smallmatrix} \int_{[0,T]\times \Omega} \begin{smallmatrix} \mbox{extremize}\\ U\end{smallmatrix} \  \scl_H (\dee(t,x),U) \, dt dx,
\end{equation*}
where the pointwise extremization of $\scl_H$ over $U$, for fixed $\dee$, is made possible by the choice of $H$.

Furthermore, assume the Lagrangian $\scl_H$ can be expressed in the form
\begin{equation*}
    \scl_H(\dee, U) := - P(\dee)\cdot U + f(U,D) + H(U)
\end{equation*}
for some function $P$ defined by the structure of the primal first-order system ((linear terms in) first derivatives of $U$ after multiplication by the dual fields and integration by parts always produce such terms), and for some function $f$ which, when non-zero, does not contain any linear dependence in $U$. Our scheme requires the existence of a function $U_H$ defined from `solving $\frac{\p \scl}{\p U} (\dee, U) = 0$ for $U$,' i.e.~$\exists \  U_H(P(\dee),\dee)$ s.t. the equation
\begin{equation*}
    - P(\dee) + \frac{\p f}{\p U}(U_H(P(\dee),\dee), \dee) + \frac{\p H}{\p U}\left(U_H(P(\dee),\dee)\right) = 0
\end{equation*}
is satisfied. This requirement may be understood as follows: define
\begin{equation*}
    f(U, \dee) + H(U) =: M(U, \dee)
\end{equation*}
and assume that it is possible, through the choice of $H$, to make the function $\frac{\p M}{\p U}(U, \dee)$ \textit{monotone} in $U$ so that a function $U_H(P,\dee)$ can be defined that satisfies
\begin{equation*}
    \frac{\p M}{\p U}(U_H(P,\dee), \dee) = P, \quad \forall P.
\end{equation*}
Then the Lagrangian is
\begin{equation*}
    \scl(\dee, U_H(P(\dee),\dee)) = - P(\dee) \cdot U_H(P(\dee),\dee) + M(U_H(P(\dee),\dee), \dee) =: - M^*(P(\dee), \dee)
\end{equation*}
where $M^*(P,\dee)$ is the Legendre transform of the function $M$ w.r.t $U$, with $\dee$ considered as a parameter.

Thus, our scheme may also be interpreted as designing a concrete realization of abstract saddle point problems in optimization theory \cite{rockafellar1974conjugate}, where we exploit the fact that, in the context of `solving' PDE viewed as constraints implemented by Lagrange multipliers to generate an unconstrained problem, there is a good deal of freedom in choosing an objective function(al) to be minimized. We exploit this freedom in choosing the function $H$ to develop dual variational principles corresponding to general systems of PDE.

\section{Application of base states: a simple example}\label{app:base_state_example}
\begin{figure}
    \centering
    \includegraphics[width=0.3\textwidth]{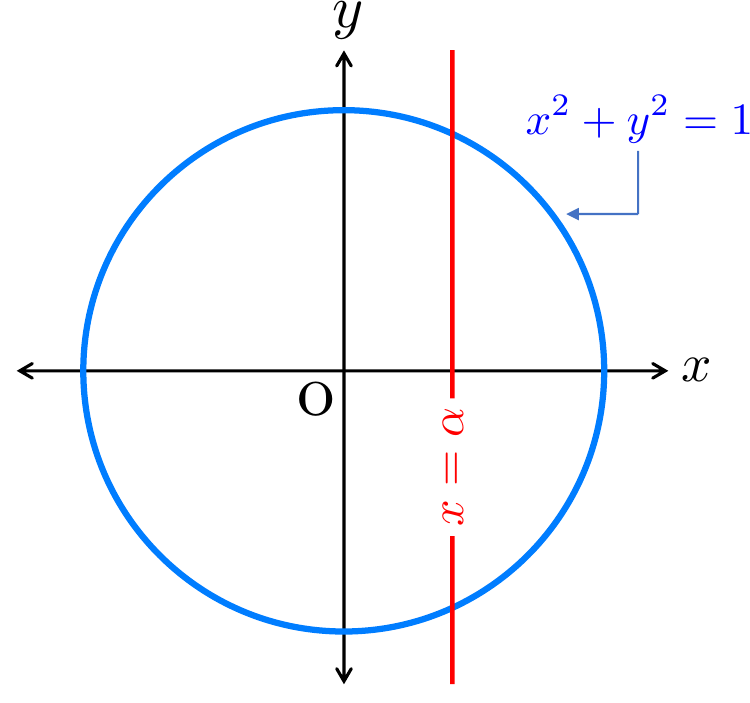}
    \caption{Schematic of Circle-line intersection.}
    \label{fig:circle_line}
\end{figure}
Consider the following algebraic system of equations for $(x,y)\in\mathbb{R}^2, \alpha \in \mathbb{R}$:
\begin{equation}
\begin{aligned}
    x^2 + y^2 & = 1\\
    x & = \alpha.
\end{aligned}
\label{eq:circle_line_primal}%
\end{equation}
A schematic is shown in Fig.~\ref{fig:circle_line}, and  solutions are given by $\{(\alpha, \pm \sqrt{1-\alpha^2}): |\alpha| \leq 1\}$.
We use the logic of Appendix \ref{sec:fin_dim} to multiply each of the above equations by a dual multiplier and add a quadratic auxiliary potential $H$:
\begin{equation*}
    \widehat{S}(x,y,\lambda,\gamma) = \lambda (x^2 + y^2 -1) + \gamma (x - \alpha) + \frac{1}{2}(x - \bar{x})^2 + \frac{1}{2}(y - \bar{y})^2,
\end{equation*}
where $\bar{x}, \bar{y}\in \mathbb{R}$ are constants (in this algebraic problem), and we will refer to the pair as a `base state.' Next we need to generate the analog of the mapping function $x_H$ of Appendix \ref{sec:fin_dim}, also referred to as the DtP mapping in the text:
\begin{equation} \label{eq:line_circle_dtp}
    \begin{aligned}
        \frac{\p \widehat{S}}{\p x} = 0 : \quad 2 \lambda x + \gamma + (x-\bar{x}) = 0 \quad & \Rightarrow \quad \mbox{for} \ \lambda \neq - \half, \quad x_H(\lambda,\gamma) = \frac{\bar{x} - \gamma}{2 \lambda + 1}; \\
        \frac{\p \widehat{S}}{\p y} = 0 : \quad 2 \lambda y + (y-\bar{y}) = 0 \quad & \Rightarrow \quad \mbox{for} \ \lambda \neq - \half, \quad  y_H(\lambda) = \frac{\bar{y}}{2 \lambda + 1}.
    \end{aligned} 
\end{equation}
{\boldmath$\lambda \neq - \half$}: Considering only the case $\lambda \neq - \half$ for the moment, the dual objective function is now obtained by substituting the DtP mapping into $\widehat{S}$:
\begin{equation*}
    S(\lambda,\gamma) = \lambda \left(x_H^2(\lambda,\gamma) + y_H^2(\lambda,\gamma) -1 \right) + \gamma \left(x_H(\lambda,\gamma) - \alpha \right) + \frac{1}{2}\left(x_H(\lambda,\gamma) - \bar{x} \right)^2 + \frac{1}{2}\left(y_H(\lambda,\gamma) - \bar{y} \right)^2.
\end{equation*}

The critical point equations for this objective, by design, are the equations \eqref{eq:circle_line_primal} with the substitution $(x \to x_H, y \to y_H)$:
\begin{equation}
\begin{aligned}
    \frac{(\bar{x}-\gamma)^2}{(2\lambda+1)^2}+\frac{\bar{y}^2}{(2\lambda+1)^2} & = 1 \\ 
    \frac{\bar{x} - \gamma}{2\lambda+1} & = \alpha.
\end{aligned}
\label{eq:circle_line_dual}
\end{equation}
A necessary condition for solutions is
\[
(2\lambda + 1)^2 (1 - \alpha^2) - \bar{y}^2  = 0
\] 
which implies that dual solutions exist only for $|\alpha| \leq 1$ and when $|\alpha| = 1$ only if $\bar{y} = 0$. 

Thus, for $|\alpha| < 1$, $\bar{y} \neq 0$,
\begin{equation*}
    \lambda = \frac{1}{2}\left(\pm\frac{|\bar{y}|}{\sqrt{1-\alpha^2}}-1\right); \qquad \gamma = \bar{x} - \alpha(2\lambda + 1).
\end{equation*}
are the extrema of $S$.

Since only $\lambda \neq - \half$ is being considered, we do not consider the case $|\alpha| < 1$, $\bar{y} = 0$ here.

For $|\alpha| = 1$, $\bar{y} = 0$, the pairs
\begin{equation*}
    - \half \neq \lambda \in \R \ \mbox{arbitrary}; \qquad \gamma = \bar{x} - \alpha(2\lambda + 1).
\end{equation*}
are the extrema of $S$. 

Putting these dual solutions back into the DtP mapping \eqref{eq:line_circle_dtp}, we recover the correct primal solutions as expected.

\noindent {\boldmath$\lambda = - \half$}: We note that the function $S$ is not unambiguously defined at $\lambda = - \half$ and therefore its gradient at any $\left(\lambda,\gamma \right) = \left(-\half,\gamma \right)$ cannot be defined. Hence, it does not make strict sense to talk about critical points of $S$ at such locations.

Nevertheless, it is noted that for $\bar{y} = 0$, the DtP mapping equation \eqref{eq:line_circle_dtp} and the critical point equation \eqref{eq:circle_line_dual} for $\lambda = -\half$ can be made sense of, and are satisfied, at $(\lambda = -\half, \gamma = \bar{x})$. However, even with such a relaxed interpretation of a `critical point' for the dual problem, this dual solution does not define solutions to the primal problem in a non-vacuous manner (i.e.~any primal state $(x,y)$ is admitted as a solution, and the procedure does not provide specific guidance for generating primal solutions).

The following conclusion can be drawn from this simple, yet non-trivial, example:
\begin{itemize}
    \item A `good' choice of the auxiliary function can be crucial for the success of the dual scheme in generating solutions to the primal problem. For example, if `no base states' are invoked within this class of quadratic auxiliary functions, i.e.~$(\bar{x}, \bar{y}) = (0,0)$, while dual extrema exist, the scheme essentially fails to define primal solutions, except for the case $|\alpha| = 1$. And when the latter is the case for the primal problem, then a base state with non-zero $\bar{y}$ is not feasible for defining dual and primal solutions to the problem.

    Of course, the choice of base states and the auxiliary function $H$ defined with (or without) their aid is largely arbitrary within the scheme, so that the flexibility in making such  choices is not a limitation of the approach.
\end{itemize}

\section{Weak form and conservation law for the Burgers equation} \label{app:new}
\subsection{Burgers equation and jump condition from the Conservation Law} \label{sec:conservation_law}

Burgers equation was designed to study the Navier-Stokes equations in the simplest possible setting. It represents an idealization of the equation of balance of linear momentum for a Newtonian viscous fluid, written in an Eulerian (as opposed to a Lagrangian) setting. The role of the mass density is ambiguous in this idealization since the mass density does not appear in the equation, but the fluid is not incompressible either. The one-dimensional conservation law governing the fluid velocity $u$ is defined as follows: For any arbitrarily fixed spatial domain $(l,r) \subset (x_l,x_r)$,
\
\begin{equation}
    \frac{d}{dt}\int_{l}^{r}  u \, dx + \Big( \tilde{F}(u(r,t)) - \tilde{F}(u(l,t)) \Big) = 0
    \label{eq:conservation_law}
\end{equation}
where $$\tilde{F}(u) = \frac{u^2}{2} - \nu\, \p_x u$$  represents the flux of (idealized) linear momentum. The first term is the flux arising from advection, and the second from the diffusive viscous stress. The case $\nu=0$ represents the inviscid approximation of Burgers equation. We will often write  $\tilde{F}(u(x,t))=: F(x,t)$ with a slight abuse of notation. 

In the inviscid case, we think of solutions to be in the class of bounded, piecewise-continuously-differentiable functions on the space-time domain 
\begin{equation}
\Omega:=(x_l,x_r) \times (0,T)
\label{eq:domain_appendix}
\end{equation}
except, possibly, on a finite number of curves. In the viscous case, we think of solutions $u$ belonging to the space $C^{0}$ (and possibly with higher regularity).

 Considering a generalization of the integral  on the left of the equality given by 
 \begin{equation*}
     I_{(l,r)} := \int_{l}^{r} u\,dx
 \end{equation*}
in \eqref{eq:conservation_law}(the subscript ``$(l,r)$'' denotes the spatial limits of the integral), we examine
\begin{equation*}
    \int_{a(t)}^{b(t)}  u \, dx,
\end{equation*}
where $a,b: [0,T] \to \R$ are prescribed continuous functions of time.

Assuming that $u$ is continuously differentiable for $x\in(a(t),b(t)) \subset \R$ and $t\in[0,T]$, we define a  smooth map 
\[\hat{x}:[0,1]\times [0,T]\rightarrow [a(t), b(t)] \subset \R, \qquad (\xi,t) \mapsto \hat{x}(\xi,t)
\]
which is injective for each fixed $t$ and satisfies $\hat{x}(0,t) = a(t)$ and $\hat{x}(1,t) = b(t)$. Then, at any fixed time $t$,
$$\int_{a(t)}^{b(t)}  u(x,t) \, dx = \int_0^1\,u(\hat{x}\big(\xi,t),t\big)\,\frac{\p \hat{x}}{\p \xi}(\xi,t)\, d\xi.$$
  Consider 
\begin{align}
\frac{d}{dt}\int_{a(t)}^{b(t)}  u \, dx &= \frac{d}{dt}\int_0^1\,u\,\frac{\p \hat{x}}{\p \xi}\, d\xi =  \int_0^1\,\frac{\p}{\p t}\left(u\,\frac{\p \hat{x}}{\p \xi}\,\right) d\xi \nonumber \\
&= \int_0^1 \frac{\p \hat{x}}{\p \xi}\left( \frac{\p u}{\p x}\frac{\p \hat{x}}{\p t} + \frac{\p u}{\p t} \right)d\xi \,+ \int_0^1 u\,\frac{\p^2 \hat{x}}{\p \xi \p t} d\xi
\nonumber \\ &= \int_0^1 \frac{\p}{\p \xi}\left( u \frac{\p \hat{x}}{\p t} \right)d\xi + \int_0^1\,\frac{\p u}{\p t}\frac{\p \hat{x}}{\p \xi} \, d\xi, \label{eq:reduction}
\end{align}
where we have used
$$\frac{\p u}{\p x} \frac{\p \hat{x}}{\p \xi} = \frac{\p u}{\p \xi}.  $$
The first integral in \eqref{eq:reduction} becomes
\begin{equation*}
\int_0^1 \frac{\p}{\p \xi}\left( u \frac{\p \hat{x}}{\p t} \right)\,d\xi = \left. u(\hat{x}(\xi,t),t) \,\frac{\p \hat{x}}{\p t}(\xi,t) \right|_{\xi=0}^{\xi=1} = u\big(b(t),t\big) \frac{d b}{dt}(t) - u\big(a(t),t\big)\frac{d a}{dt}(t).
\end{equation*}
The second integral in \eqref{eq:reduction} becomes
\begin{equation*}
\int_0^1\,\frac{\p u}{\p t}(\hat{x}(\xi,t),t)\frac{\p \hat{x}}{\p \xi}(\xi,t) \,d\xi = \int_{a(t)}^{b(t)} \frac{\p u}{\p t}(x,t)\, dx.
\end{equation*}
Thus
\begin{equation}
    \frac{d}{dt}\int_{a(t)}^{b(t)}  u \, dx = \int_{a(t)}^{b(t)} \frac{\p u}{\p t} \,dx + u\big(b(t),t\big) \frac{d b}{dt}(t) - u\big(a(t),t\big)\frac{d a}{dt}(t).
    \label{eq:leibniz}
\end{equation}
For the current section, we perform the following analysis by focusing on local space-time subdomains 
\begin{equation}
\Omega\supset\Omega_s:=(l,r)\times(t_i, t_f) \label{eq:subdomain}
\end{equation}
such that at most one interface on which $u$ is not continuously differentiable traverses through $\Omega_s$. Let $\mathcal{I}$ denote such an interface and be defined as follows:  \begin{equation}
\mathcal{I}= \{(x,t):x=s(t), x\in(l,r) \mbox{ and } t\in(t_i,t_f)\}.\label{eq:curveSt}
\end{equation}
where $s: (t_i,t_f) \to (l, r) \in \R$ is a smooth function of time. We additionally define the following limits at any time $t$:
\begin{align*}
  & u^-(s(t),t) := lim_{x\to s(t)}\, u(x,t), \quad x<s(t); \\
  & u^+(s(t),t) := lim_{x\to s(t)}\, u(x,t), \quad x>s(t);\\
  & F^-(s(t),t) := lim_{x\to s(t)}\, F(x,t), \quad x<s(t); \\
  & F^+(s(t),t) := lim_{x\to s(t)}\, F(x,t), \quad x>s(t).
\end{align*}
Utilizing the result \eqref{eq:leibniz}, we obtain the following within $\Omega_s$:
\begin{align}
     \frac{dI_{(l,r)}}{dt} &= \frac{d}{dt}\left(\int_{l}^{r}  u \, dx\right) = \frac{d}{dt}\left(\int_{l}^{s(t)}  u \, dx + \int_{s(t)}^{r}  u \, dx\right) \nonumber \\
    & = \left(\int_{l}^{s(t)} \frac{\partial u}{\partial t} dx + u^- \frac{ds}{dt} \right) +\left(\int_{s(t)}^{r} \frac{\partial u}{\partial t} dx -u^+ \frac{ds}{dt}\right). \label{eq:partition}
\end{align}
Substituting the last expression in \eqref{eq:conservation_law}, we obtain:
\begin{align}
&\frac{dI_{(l,r)}}{dt} + \Big( F(r,t) - F(l,t) \Big) = 0. \nonumber \\
\Rightarrow & \begin{aligned}[t]
&\left(\int_{l}^{s(t)} \frac{\partial u}{\partial t} dx + u^- \frac{ds}{dt}\right) + \left(\int_{s(t)}^{r} \frac{\partial u}{\partial t} dx - u^+ \frac{ds}{dt}\right) + \Big( F(r,t) - F(l,t) \Big)  = 0. \nonumber
\end{aligned} \\
\Rightarrow & \int_{l}^{s(t)} \frac{\partial u}{\partial t} dx  + \int_{s(t)}^{r} \frac{\partial u}{\partial t} dx   
 + \big( F(r,t) - F(l,t) \big) - \frac{ds}{dt} \big( u^+ - u^- \big) = 0.\label{eq:eliminate}
\end{align}
By letting $l\rightarrow s(t)$ and $r\rightarrow s(t)$, the first two integrals in the last expression vanish and we are left with
\begin{equation}
\frac{ds}{dt}(t) =\frac{  F^+(s(t),t) - F^-(s(t),t)} { u^+(s(t),t) - u^-(s(t),t)},
\label{eq:jump_conservation}
\end{equation}
which is the jump condition (commonly referred to as the Rankine-Hugoniot condition). In case no such interface exists ($u^+ = u^-$), \eqref{eq:eliminate} reduces to:
\begin{align}
&\int_{l}^{r} \frac{\partial u}{\partial t} dx + \Big( F(r,t) - F(l,t) \Big) = 0. \nonumber \\
\Rightarrow & \int_{l}^{r} \left(\frac{\partial u}{\partial t} + \frac{\p F}{\p x} \right) dx = 0. \nonumber\\
\Rightarrow & \,\frac{\partial u}{\partial t} +  \frac{\p F}{\p x} = 0 \quad \forall x\in(l,r) \mbox{ and }t\in(t_i,t_f). \label{eq:burger_conservation}
\end{align}

\subsection{Burgers equation and the jump condition from the weak form}\label{sec:weak_form_implications}
We start by considering the weak form of the inviscid Burgers equation. Let 
\begin{equation}
    \mathcal{V}=\{\delta u: \delta u \in C^1(\bar{\Omega}) \text{ and } \delta u(x,t) = 0 \, \forall \,(x,t)\in\p \Omega\},
    \label{eq:subspace}
\end{equation}
where $\Omega$ was defined in \eqref{eq:domain_appendix} and $\p \Omega$ represents the boundary of $\Omega$. We consider $u$ to be piecewise continuously differentiable in $\Omega$ with discontinuities in $u$ and/or $\nabla u$ concentrated on the union of at most a finite number of curves in $\Omega$. Additionally, we require $u$ to satisfy, for any $\delta u \in \calV$, the following weak form of Burgers equation:
\begin{equation}
    \int_\Omega \big(u\,\p_t \delta u + F\, \p_x \delta u \big)\, dt\,dx  = 0,
\label{eq:burger_distributional}
\end{equation} 
where $F \circ u$ is $C^1$ at any point of $\Omega$ where $u$ is $C^1$.
\begin{figure}
    \centering
    \begin{subfigure}[b]{0.4\textwidth}
        \centering
        \includegraphics[scale=0.35]{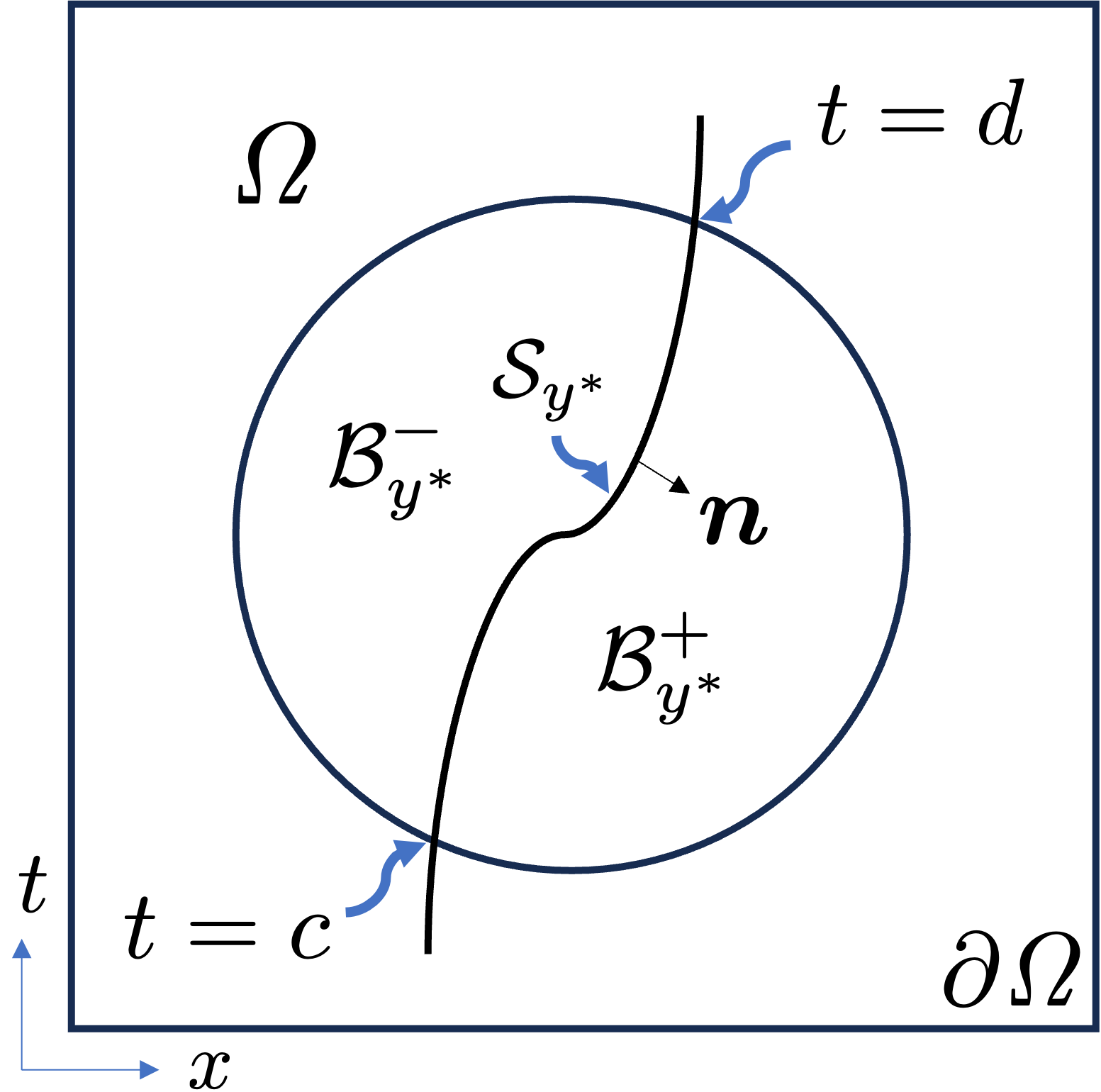}
    \caption{ }
    \label{fig:discontinuity}
    \end{subfigure}
    \hfill
    \begin{subfigure}[b]{0.4\textwidth}
\centering        \includegraphics[scale=0.5]{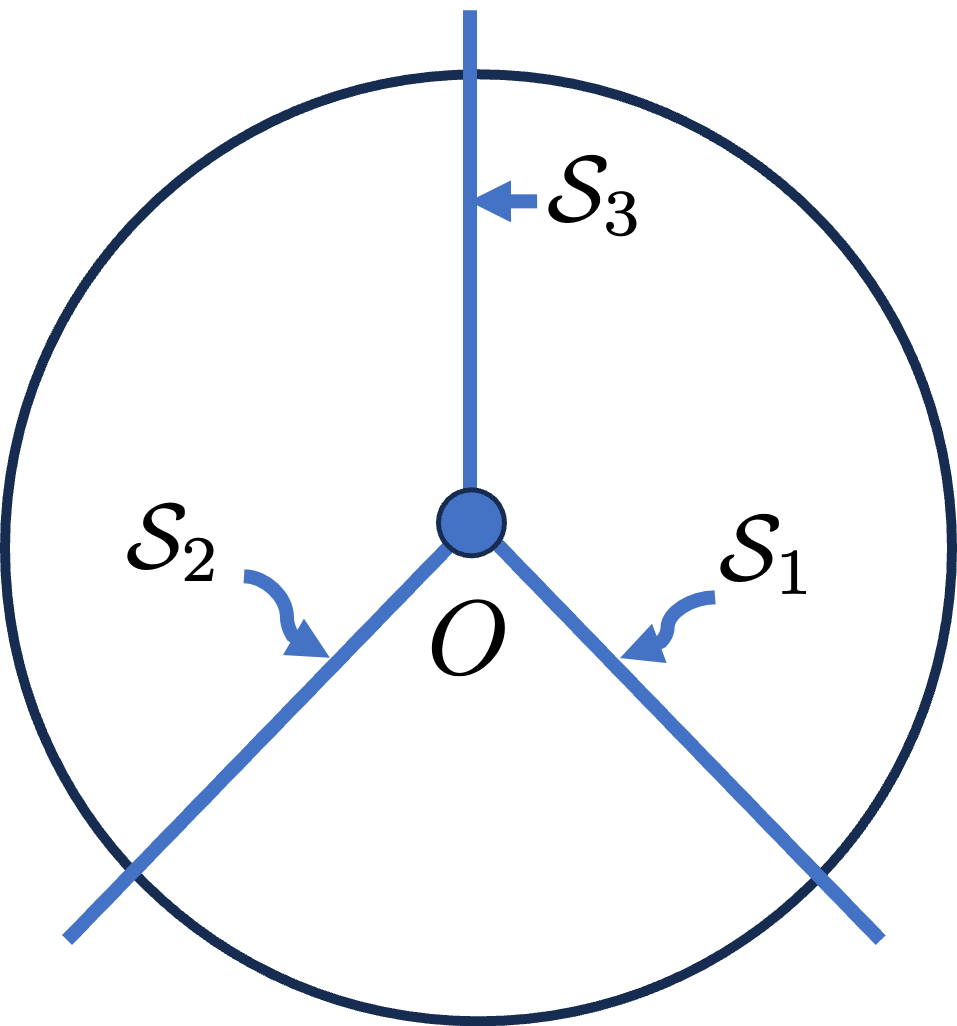}
    \caption{}
    \label{fig:discon2}
    \end{subfigure}
    \caption{(a) $\boldsymbol{n} = (n_1,n_2)$ represents the outward unit normal to the region $\mathcal{B}_{y^*}^-$ along the curve $\mathcal{S}_{y^*}$. $x$ and $t$ directions will be represented by indices 1 and 2, respectively. (b) $\mathcal{S}_1,\mathcal{S}_2$ and $\mathcal{S}_3$ represent three different curves on which $u$ is not continuously differentiable, with $O$ representing their junction point.}
    \label{fig:mainfigure}
\end{figure}

 Let $\calS = \overline{\cup_i^n \calS_i}$ where $\{\calS_i: i = 1,\ldots,n\}$ is a collection of disjoint curves in $\Omega$ across which $u$ is not continuously differentiable. For $y^* := (x^*,t^*)\in\Omega\backslash\mathcal{S}$, consider an open ball $\mathcal{B}_{y*}$ centered at $y^*$ in which $u$ is continuously differentiable. For $\delta u \in \calV$ with $\delta u(x,t)=0 \, \forall \, (x,t) \notin \mathcal{B}_{y^*}$, eq.~\eqref{eq:burger_distributional} reduces to
\begin{equation*}
    \int_{\mathcal{B}_{y^*}} \big(u\,\p_t \delta u + F\, \p_x \delta u \big)\, dt\,dx  = 0 \quad \Longrightarrow \quad \int_{\mathcal{B}_{y^*}} -\delta u\big(\p_t u +  \,\p_x F\big)\, dt\,dx=0.
\end{equation*}
Noting that the coefficient of $\delta u$ in the above integral is continuous on $\calB_{y^*}$, it has to vanish for the class of  $\delta u\in \mathcal{V}$ with support in $\calB_{y^*}$ being considered for the integral to vanish. Hence,
\begin{equation}
    \p_t u +  \,\p_x F = 0 \mbox{ on } \Omega\backslash \mathcal{S},\label{eq:weak_1}  
\end{equation}
since the point $y^*$ was chosen arbitrarily from $\Omega\backslash \mathcal{S}$. 

Next, consider a point $y^* \in \calS_i$ for some $i \in \{1,2,\ldots,n\}$, and an open ball $\mathcal{B}_{y^*}$ around it which does not intersect any other $\calS_j, j \in \{1,2,\ldots,n\}, j \neq i$. Let $\calS_{y^*} = \calS_i \cap \calB_{y^*}$ be parametrized as
\begin{equation}
\mathcal{S}_{y^*}= \{(x,t):x=s(t), t\in(c,d)\, \subset (0,T)\},\label{eq:curveSt2}
\end{equation} 
where $s:(c,d)\rightarrow(x_l,x_r)$ is a smooth function representing $\mathcal{S}_{y^*}$, parameterized by $t$ ($s$ from \eqref{eq:curveSt} is being redefined for the current section). The schematic for such a setup has been shown in \ref{fig:discontinuity}, where
\begin{align*}&\mathcal{B}_{y^*} \supset \mathcal{B}_{y^*}^-:= \{(x,t):x<s(t)\}; \\&\mathcal{B}_{y^*} \supset\mathcal{B}_{y^*}^+:= \{(x,t):x>s(t)\}; 
\\ u(x,t) \to & \,u^+(s(\hat{t}\,),\hat{t}\,) \text { for } \mathcal{B}_{y^*}^+ \ni (x,t) \to (s(\hat{t}\,),\hat{t}\,); \\ 
 u(x,t) \to &\,u^-(s(\hat{t}\,),\hat{t}\,) \text { for } \mathcal{B}_{y^*}^- \ni (x,t) \to (s(\hat{t}\,),\hat{t}\,); \\
F(x,t) \to &\,F^+(s(\hat{t}\,),\hat{t}\,) \text { for } \mathcal{B}_{y^*}^+ \ni (x,t) \to (s(\hat{t}\,),\hat{t}\,); \\ 
 F(x,t) \to &\,F^-(s(\hat{t}\,),\hat{t}\,) \text { for } \mathcal{B}_{y^*}^- \ni (x,t) \to (s(\hat{t}\,),\hat{t}\,),
\end{align*}
and $\left(s(\hat{t}\,),\hat{t}\,\right)\in \mathcal{S}_{y^*}$. We assume that the various limits $()^{\pm}$ exist and are continuous on $\calS_{y^*}$. Let $\delta u(x,t)=0 \, \forall \, (x,t) \notin \mathcal{B}_{y^*}$, $\delta u \in \calV$.
Thus eq.~\eqref{eq:burger_distributional} reduces to
\begin{equation*}
    \int_{\mathcal{B}_{y^*}} \big(u\,\p_t \delta u + F\, \p_x \delta u \big)\, dt\,dx  = 0 \Rightarrow \int_{\mathcal{B}_{y^*}^-} \big(u\,\p_t \delta u + F\, \p_x \delta u \big)\, dt\,dx + \int_{\mathcal{B}_{y^*}^+} \big(u\,\p_t \delta u + F\, \p_x \delta u \big)\, dt\,dx=0,
\end{equation*}
which upon utilizing the divergence theorem and \eqref{eq:weak_1} leads to the following line integral on $\mathcal{S}_{y^*}$:
\begin{equation}\label{eq:full}
\int_{\mathcal{S}_{y^*}} \delta u(s(t),t)\,\Big(\big(F^- - F^+\big)\, n_1 + \big(u^- - u^+\big)\, n_2 \Big)d\mathcal{S}=0.
\end{equation} 
Here, $\boldsymbol{n}:=(n_1,n_2)$ represents the unit normal to  $\mathcal{S}_{y^*}$ pointing outwards from $\mathcal{B}_{y^*}^-$ as shown in Fig.~\ref{fig:discontinuity} and $d\mathcal{S}$ represents the infinitesimal arc length.
Using \eqref{eq:curveSt2}, we obtain the following:
$$\frac{d\mathcal{S}}{dt} = \sqrt{1 + \left(\frac{ds}{dt}\right)^2}; \quad \boldsymbol{n} =  \frac{1}{\sqrt{1 + \left(\frac{ds}{dt}\right)^2}} \left(1,-\frac{ds}{dt}\right).$$
Substituting the above expressions into \eqref{eq:full} leads to:
\begin{equation}
     \int_a^b \delta u(s(t),t)\,\left(\big(F^- - F^+\big)\, - \big(u^- - u^+\big)\,\frac{ds}{dt} \right)dt=0.
    \label{eq:full2}
\end{equation}
Since this holds for all smooth functions $\delta u$ considered, the coefficient of $\delta u(s(t),t)$ being continuous along the curve ($u^+,u^-,F^+$ and $F^-$ are continuous on $\mathcal{S}_{y^*}$) implies:
$$\big(F^- - F^+\big)\, - \big(u^- - u^+\big)\,\frac{ds}{dt} = 0.$$
The above statement can be rewritten as
\begin{equation}
 \frac{ds}{dt}(t) = \frac{  F^+(s(t),t) - F^-(s(t),t) } { u^+(s(t),t) - u^-(s(t),t)} \quad\forall t\,\in(c,d).
 \label{eq:weak_jump}
 \end{equation}
The ball $\mathcal{B}_{y^*}$ can then be centered around any point on any of the curves $\mathcal{S}_i$. The last statement is the jump condition.

At points where multiple curves on which $u$ is not continuously differentiable intersect, e.~\!g., the point $O$ in Fig.~\ref{fig:discon2}, the functions $u^+, u^-, F^+, F^-$ need not necessarily be continuous, and there is no analog of the jump condition implied by the weak formulation of the inviscid Burgers equation. This situation arises in the N-wave problem (Sec.~\ref{sec:burgers_N-wave}, App.~\ref{app_nwave}) where two lines of weak discontinuity  ($\p_x u$ is discontinuous across the two inclined red lines in Fig.~\ref{fig:app_bur_p5}) join to give rise to the origin of a strong/shock discontinuity curve (the junction of the two aforementioned lines and the vertical $x$-patterned line). Such a condition also arises in the double-shock problem (Sec.~\ref{sec:Burgers_dshock}, App.~\ref{app_dshock}), where two lines of strong discontinuity merge to result in one strong discontinuity.

It is also straightforward to show that the Burgers equation and jump conditions together imply the weak form of the equation.

\subsection{Burgers equation and the jump condition imply that $u$ is conserved in the Half N-wave problem (Sec.~\ref{sec:Burgers_halfwave}, App.~\ref{app:halfwave})} 
Based on the Burgers equation \eqref{eq:weak_1} (or \eqref{eq:burger_conservation}) and the jump condition \eqref{eq:weak_jump} (or \eqref{eq:jump_conservation}), the conservation law can be achieved as follows:
Consider the construction presented in the Appendix \ref{sec:conservation_law}. Then
\begin{equation}
    \frac{ds}{dt} \Big(u^+-u^-\Big) =\Big(F^+-F^-\Big) \quad \forall (x,t)\in\mathcal{I},
\end{equation}
which on substituting into \eqref{eq:partition} leads to the following: 
$$\frac{dI_{(l,r)}}{dt} 
 = \left(\int_{l}^{s(t)} \frac{\p u}{\p t} dx   
      +\int_{s(t)}^{r} \frac{\partial u}{\partial t} dx \right) 
 - \Big(F^+-F^-\Big).$$
Adding and subtracting $F(l,t)$ and $F(r,t)$, the last equation upon rearrangement gives:
\begin{align*}
\frac{dI_{(l,r)}}{dt} 
 =& \left(\int_{l}^{s(t)} \frac{\p u}{\p t} dx +F^- - F(l,t)\right)   
      +\left(\int_{s(t)}^{r} \frac{\partial u}{\partial t} dx +F(r,t) - F^+\right) \\
&  \qquad- \Big(F(r,t)-F(l,t)\Big) \\
=& \int_{l}^{s(t)} \left(\frac{\partial u}{\partial t} + \frac{\partial F}{\partial x}\right) dx 
      +\int_{s(t)}^{r} \left(\frac{\partial u}{\partial t} + \frac{\partial F}{\partial x}\right) dx -\Big(F(r,t) - F(l,t)\Big).
\end{align*} 
The domain integrals in the last expression can be eliminated using \eqref{eq:burger_conservation} (or using \eqref{eq:weak_1}) since $u$ is continuously differentiable within the limits of the integrals. The remainder implies the conservation law for the arbitrary interval $(l,r)$. 

The above analysis can also be extended to the case when multiple non-intersecting interfaces are allowed to pass through $\Omega_s$ by partitioning, for each fixed time $t$, the spatial domain $(l,r)$ into smaller intervals defined by the intersection of the interfaces with the set $(l,r) \times \{t\}$.

  \noindent   \underline{Half N-wave problem (Sec.~\ref{sec:Burgers_halfwave}, App.~\ref{app:halfwave}):} In this case, $\Omega$ contains a weak stationary discontinuity (at $x=x_0$) and a strong traveling discontinuity (at $x = x_0+l(t)$) that do not intersect. Hence, the above analysis applies. Also, when the flux at the domain boundaries in this problem vanishes for the time under analysis, i.e.~$F(x_l,t) = F(x_r,t) = 0$, the integral of $u$ across $(x_l,x_r)$ remains conserved:
\begin{equation*}
   \frac{d}{dt} \int_{x_l}^{x_r} u\,dx = 0.
\end{equation*}
 \underline{A sketch of the argument for the general case}: In cases when multiple interfaces of discontinuity (weak or strong) interact directly (for e.g.~merge together at a junction), we consider the following argument. Let there be a finite number of time instants $t^*_i, i = 1, \ldots, n$, for which the number of intersections of the jump interfaces in $\Omega$ with the set $(l,r) \times \{t^*_i\}$ given by $\calN_i$ is such that $\calN_i \neq \calN_{i+1}$. Then, the conservation law will hold for all times in the intervening time partitions $(t^*_i, t^*_{i+1})$ based on the construct presented at the beginning of this section: 
\begin{equation*}
    \frac{dI_{(l,r)}}{dt} + \Big(F(r,t) - F(l,t) \Big) = 0 \quad \forall\, t \in (t^*_i, t^*_{i+1}), i = 1,\ldots, n-1, \mbox{and} \ t \in (t^*_n, T).
\end{equation*}
Furthermore,

\begin{itemize}
    \item For the class of solutions of the inviscid Burgers equation considered, $F(r,t) - F(l,t)$ can at most be discontinuous as a function of time (since $F(x,t) = \tilde{F}(u(x,t))$) (where any $t$ at which such a discontinuity arises need not necessarily be an element of $\{t^*_i: i = 1,\ldots, n\}$).
    \item For $t\in(0,T)$, we assume:
\begin{equation*}
    I_{(l,r)}(t) = \lim_{s \to t,\,s<t} I_{(l,r)}(s).
\end{equation*}
\end{itemize}
Based on the last two statements, $I_{(l,r)}$ can be expressed as a continuous function in $t$ given by:
\begin{equation}\label{eq:int_conv_u}
I_{(l,r)}(t) = I_{(l,r)}(0) + \int_{0}^t \Big(F(r,z)-F(l,z) \Big)\,dz \quad \forall t\in(0,T).
\end{equation}
Recalling that the interval $(l,r)$ was chosen arbitrarily, the satisfaction of Burgers equation and the jump condition on appropriate subsets of $\Omega$ implies that the conservation statement \eqref{eq:int_conv_u} holds for any interval $(l,r) \subset (x_l,x_r)$ for all $t \in (0,T)$. 

In case when $F(r,z)-F(l,z)=0$ for almost all times $z \in (0,T)$, $I_{(l,r)}$ remains conserved (as observed in the case of the N-wave problem, Sec.~5.1.5).

\section{$L^2$ projection}\label{app:L2}
 Let $u$ represent the primal field under consideration, which may depend on the derivatives of the dual fields (approximated using linear FE shape functions). To achieve a $C^0$ continuous approximation of $u$ in the domain $\Omega\subset\mathbb{R}\mbox{ or }\mathbb{R}^2 $, we employ the following method: Let $u_h$ represent the projection of $u$ onto a space $V_h$ formed by the linear span of globally continuous, piecewise smooth finite element shape functions corresponding to a FE mesh for $\Omega$.
We enforce the following condition upon $u_h$:
\begin{equation*}
     u_h = \operatorname*{arg\,min}_{v\in V_h}  \int_\Omega \frac{1}{2} |u-v|^2 \, d\Omega.\label{eq:weak_L2}
\end{equation*} 
Let $N^A$ represent the basis functions associated with any node with an index $A$. The discrete version of the optimality condition of the above statement is
\begin{equation}\label{eq:L2_project}
  \sum_{A=1}^N\sum_{B=1}^N\delta u^A_h\biggl(\int_\Omega N^A\, N^B \,d\Omega \biggl) u_h^B=  \sum_{A=1}^N\delta u^A_h \int_\Omega N^A\, u\,d\Omega,
\end{equation}
where $u_h^B$ denotes the nodal value at node $B$ of the sought continuous projection of $u$, and  $\delta u_h := \delta u^A_h N^A$ is a test function. The primal data for the integral in the right hand side of \eqref{eq:L2_project} is obtained from the Gauss points of 1-D/2-D element. In solving \eqref{eq:L2_project}, we impose any known data, e.g.~initial and boundary conditions of the primal problem, as known function values $u^A_h$, at corresponding nodes $A$ where the data is known, with $\delta u^A_h = 0$. Using the arbitrariness of the remaining `free' $\delta u^A_h$, we obtain a system of linear equations to solve for the unconstrained nodal values $u_h^B, B \in \{1, \ldots, N\}$.

If the $u$ profile of the primal field is required only at a specific time, then conducting a 1-D $L^2$ projection along that timeline is advisable since it is computationally inexpensive. As a consequence of $C^0$ type FE shape functions being employed, the time derivatives of the dual fields experience discontinuities across nodal timelines. 
In such a case the value of the time derivative on any nodal timeline can be approximated as an average of the value obtained from the element above and below the spatial point of interest (these values remain constant in the direction of time within any element).   With the inputs to the DtP mapping generated in this manner, the value of primal fields at the Gauss points (of this nodal timeline) are obtained.

\section{The base state smoothing operator}\label{app:Smoothing}
The smoothing operator $\mathcal{S}$ serves to mitigate the high-wave number contributions in any function. Let $u:=\mathcal{S}[f]$ represent the smoothed output produced for an input function $f$. Then we require $u$ to satisfy:
\begin{equation}\label{eq:smooth_strong}
    u- \eta\, \p_{xx} u = f, \qquad \eta > 0.
\end{equation}
Here, $\eta$ represents a diffusion-like control parameter.

Let $t_s$ represent the time at which the smoothing operation occurs. When boundary conditions (BC) are available at $x=0$ $(u_l(t_s))$ or (and) $x=L$ $(u_r(t_s))$ from the primal problem, we employ them as the BCs for $u$. However, when the boundary data is not available, we calculate the average of the Gauss data values  for the first element ($f_1^{g_1}$ and $f_1^{g_2}$) and set it as the left BC, and the average of the Gauss data values for the last element ($f_l^{g_1}$ and $f_l^{g_2}$) set it as the right BC on $u$. 
Correspondingly, we generate a weak form for equation
\eqref{eq:smooth_strong}. Let $\delta u$ be the space of continuous test functions on $(0,L)$ with $\delta u(0) = \delta u(L)=0$. The we intend to find $u$ which satisfies, for such test functions,
\begin{equation*} 
    \int_0^L dx\bigl( u \,\delta u + \eta\, \p_x u \, \p_x (\delta u) - f\, \delta u \bigl) = 0.
\end{equation*}
We use a standard Galerkin FEM discretization to solve the above with an $\eta=10^{-4}$.

\section{Exact entropy solutions for the provided examples}\label{app:exact_sol}
The exact entropy solutions presented in the following examples are based on the concepts and examples outlined in \cite{GILBERT_STRANG}. Given an initial condition $x \mapsto u_0(x) \in \R$, entropy solutions to the Burgers equation ($u$) are obtained using two important concepts, recalled below from \cite{GILBERT_STRANG}: 
\begin{enumerate}
    \item \label{jump_cond} Rankine-Hugoniot/jump condition: This is simply a consequence of the underlying conservation statement for Burgers equation as well as looking for weak solutions (in the sense of distributions) to it. If the solution to $\p_t u + \p_x\left(\frac{u^2}{2}\right)=0$ has different values across a shock curve $x = X(t)$, then the shock speed $s = \frac{dX}{dt}$ must satisfy $$s = \frac{1}{2}\left(u(X^+(t)) + u(X^-(t))\right),$$ where $X^-$ and $X^+$ denote the regions immediately to the left and right of the shock, respectively.
    \item \label{item:entropy}The Entropy condition: The characteristic lines must go into the shock as $t$ increases, so those on the left go faster (measured by $\frac{dX}{dt}$) than the shock and those on the right go slower. With respect to the Burgers equation this implies $$u(X^-(t))>s>u(X^+(t)).$$
\end{enumerate}

\subsection{Expansion fan}
Consider the initial condition 
\begin{equation*}
    u_0(x) = \begin{cases}
0 &  \mbox{for } x< 0.5 \\ 
1 &  \mbox{for } x>0.5.
\end{cases}
\end{equation*}
The initial signals for $x>0.5$ travel along characteristics with slope $1$, while the signals at $x<0.5$ travel along vertical characteristics (slope $0$, i.e. do not move in space). This leads to an expansion fan filling the space-time in between $x=0.5$ and $x=0.5+t$, which is a rarefaction wave traveling along $x$. The corresponding expression for $u(x,t)$ can then be given as
\begin{equation}\label{eq:Fan_ex_u}
u(x,t) = \begin{cases}
\quad \, 0 &  \mbox{for } \quad \, x< 0.5 \\[2pt]  
\dfrac{x-0.5}{t} &  \mbox{for } \,0.5+t>x\geq0.5 \\[4pt] 
\quad \, 1 & \mbox{for } \quad \,x\geq0.5+t.
\end{cases}
\end{equation}
The corresponding $Y$ profile for $Y_0(0) = 0$ is given by
\begin{equation}
Y(x,t) = \begin{cases}
\qquad\quad\, 0 &  \mbox{for } \qquad  \, x< 0.5 \\[2pt]  
\dfrac{x^2}{2t} - \dfrac{0.5\,x}{t} + \dfrac{1}{8t} &  \mbox{for } \,0.5+t>x\geq0.5 \\[6pt] 
\quad \, x - 0.5 - \dfrac{t}{2} & \mbox{for } \quad \,x\geq0.5+t.
\end{cases}
\label{eq:fan_ex_Y}
\end{equation}

\subsection{Shock} \label{sec:shock_exact}
Consider the initial condition
\begin{equation*}
    u_0(x) = \begin{cases}
1 &  \mbox{for } x< 0.5 \\ 
0 &  \mbox{for } x>0.5.
\end{cases}
\end{equation*}
Given the presence of a jump in the initial condition from 1 to 0 at $x=0.5$, the characteristics originating from the left will merge with the characteristics from the right, resulting in the traveling shock wave with velocity $s = \frac{1+0}{2} = \frac{1}{2}$. The exact entropy solution is given by
\begin{equation}
    u(x,t) = \begin{cases}
1 &  \mbox{for } x< 0.5+t/2 \\ 
0 &  \mbox{for } x>0.5+t/2.
\end{cases}
\label{eq:Shock_ex_u}
\end{equation}
The corresponding $Y$ profile for $Y_0(0) = 0$ is given by
\begin{equation} \label{eq:Shock_ex_Y}
    Y(x,t) = \begin{cases}
x-\frac{t}{2} &  \mbox{for } x< 0.5+t/2 \\ 
0.5 &  \mbox{for } x>0.5+t/2.
\end{cases}
\end{equation}

\subsection{Double shock} \label{app_dshock}
Consider the initial condition
\begin{equation*}
u_0(x) =
\begin{cases}
    \begin{aligned}
        1 & \quad \text{for } x < 0.25 \\
        0.5 & \quad \text{for } 0.25 < x < 0.5 \\
        0 & \quad \text{for } x > 0.5.
    \end{aligned}
\end{cases}
\end{equation*}
Following the same ideas from example \ref{sec:shock_exact}, the shock starting at $x=0.25$ will travel with a velocity of $s_1 = \frac{1+0.5}{2} = 0.75$, while the shock starting at $x=0.5$ will travel at with a velocity $s_2 = \frac{0.5+0}{2} = 0.25$. Since, $s_1>s_2$, the first shock will merge with the second shock at a time $t_m$ and travel with a velocity $s=\frac{1+0}{2} = 0.5$. The time $t_m$ for the faster shock to catch up with the slower shock is given by
\begin{equation*}
    s_1 \,t_m = s_2\,t_m +0.25 \quad \Rightarrow \quad  t_m = 0.5.
\end{equation*}
Consequently the shocks meet at $x=0.5+s_2 t_m=0.625$. The solution $u(x,t)$ is then be given by

\begin{equation}\label{eq:DShock_ex_u}
u(x,t) =
\begin{cases}
    \begin{rcases}
        1 & \,\text{for }\, \qquad x < 0.25 + 0.75t \\
        0.5 & \,\text{for }\, 0.25 + 0.75t < x < 0.5 + 0.25t \\
        0 & \,\text{for } \,\qquad x > 0.5 + 0.25t
    \end{rcases}
     t < t_m \\
     \\
    \begin{rcases}
        1 & \quad\text{for } \, \qquad x < 0.625 + 0.5(t-t_m) \\
        0 & \quad \text{for }\, \qquad x > 0.625 + 0.5(t-t_m)
    \end{rcases}
     t > t_m.
\end{cases}
\end{equation}
The corresponding $Y$ profile for $Y_0(0) = 0$ is given by
\begin{equation}\label{eq:DShock_ex_Y}
Y(x,t) =
\begin{cases}
    \begin{rcases}
          x - t/2&\,\text{for }\, \qquad x < 0.25 + 0.75t \\
        x/2 +0.125  -t/8 & \,\text{for }\, 0.25 + 0.75t < x < 0.5 + 0.25t \\
        0.375 & \,\text{for } \,\qquad x > 0.5 + 0.25t
    \end{rcases}
     t < t_m \\
     \\
    \begin{rcases}
        x-t/2 & \qquad\text{for } \, \qquad x < 0.625 + (t-t_m) \\
        0.375 & \qquad \text{for }\, \qquad x > 0.625 + (t-t_m)
    \end{rcases}
     t > t_m.
\end{cases}
\end{equation}
\begin{figure}\label{fig:halfwave_1}
    \centering
    \includegraphics[width=0.4\textwidth]{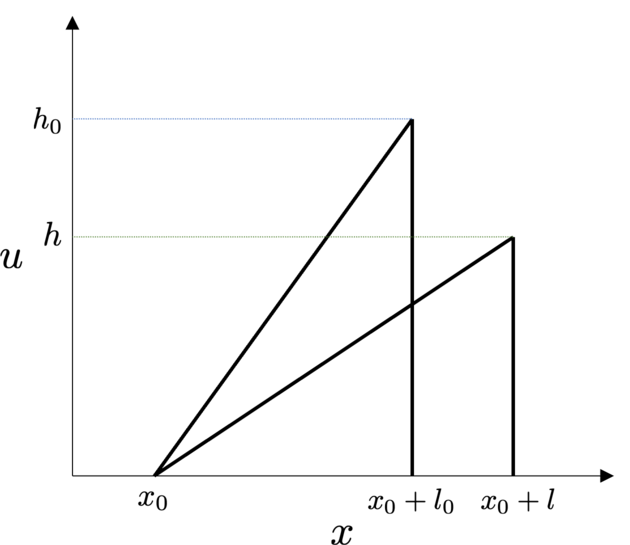}
    \caption{The right-angled triangle with a height $h_0$ indicates the initial profile for the Half N-wave. The right-angled triangle with the height $h$ indicates the anticipated profile at some later time.} 
    \label{fig:app_halfwave}
\end{figure}

\subsection{Half N-Wave}\label{app:halfwave}
Consider the initial condition
\begin{equation*}
\begin{aligned}
u_0(x) =
\begin{cases}
    0 & \text{for } x < x_0 \\[2pt]
    \dfrac{h_0}{l_0}(x-x_0) & \text{for } x_0 \leq x < x_0 + l_0 \\[4pt]
    0 & \text{for } x > x_0 + l_0.
\end{cases}
\end{aligned}
\end{equation*}
This profile has been depicted in Fig.~\ref{fig:app_halfwave}, and it is apparent that the point at $x=x_0+l_0$ corresponds to a shock, as $u_0((x_0+l_0)^-)>u_0((x_0+l_0)^+)$. Furthermore, the point at $x_0$ serves as a source for an expansion fan due to the linear increase observed from $x=x_0$ to $x=x_0 + l_0$. For any time $t$, let the shock be present at a point $X(t) = x_0 + l(t)$, where $l(0)=l_0$ and the height of shock be given by $h(t)$, where $h(0)=h_0$. To obtain the behavior of the shock front, we utilize the following information:
\begin{itemize}
    \item  Burgers equation on $\R$ is a conservation law for $u$ which is conserved:
  \begin{equation*}
       \frac{d}{dt}\int_{x_l}^{x_r}  u \, dx = 0  \quad\Rightarrow
      \int_{x_l}^{x_r}  u \, dx = \mbox{Const.},
  \end{equation*}  
  where $x_l<x_0$ and $x_l>x_0 + l_0$ such that $x_l$ is sufficiently large so that for the time under analysis, the triangular profile stays within $[x_l, x_r]$.
   This implies
  \begin{equation*}
   \half h_0\, l_0 = \half h(t) \,l(t).   
  \end{equation*}
  \item Using the jump condition and the last obtained expression across the shock, we establish that
  \begin{equation*}
      s(t) = \frac{dX(t)}{dt}=h(t)/2 \quad\Rightarrow \frac{dl(t)}{dt}=\half  \frac{h_0\,l_0}{l(t)}
  \end{equation*}
\end{itemize}
The last expression leads to 
\begin{equation*}
    l(t) =  \sqrt{h_0\,l_0\,t + l_0^2} ; \qquad X(t) = x_0 +l(t); \quad h(t) = \frac{h_0\,l_0}{l(t)}.
\end{equation*}
The solution $u(x,t)$ is given by
\begin{equation}
\label{eq:halfwave_ex_u}
\begin{aligned}
u(x,t) =
\begin{cases}
    0 & \text{for } x < x_0 \\[2pt]
    \dfrac{h(t)}{l(t)}(x-x_0) & \text{for } x_0 \leq x < x_0 + l(t) \\[4pt]
    0 & \text{for } x > x_0 + l(t).
\end{cases}
\end{aligned}
\end{equation}
The corresponding $Y$ profile for $Y_0(0) = 0$ is given by
\begin{equation}
\label{eq:halfwave_ex_Y}
\begin{aligned}
Y(x,t) =
\begin{cases}
    0 & \text{for } x < x_0 \\[4pt] 
    \dfrac{h(t)}{l(t)}\left(\dfrac{x^2+x_0^2}{2}-x_0\,x\right) & \text{for } x_0 \leq x < x_0 + l(t) \\[10pt] 
    \dfrac{h_0\,l_0}{2} & \text{for } x > x_0 + l(t).
\end{cases}
\end{aligned}
\end{equation}
For the example considered in this paper, we have used the values $x_0=0.25$, $l_0=0.25$ and $h_0=2$.
\subsection{N-wave}\label{app_nwave}
\begin{figure}
\begin{subfigure}[t]{.55\textwidth}
  \centering
  \includegraphics[width=0.9\linewidth]{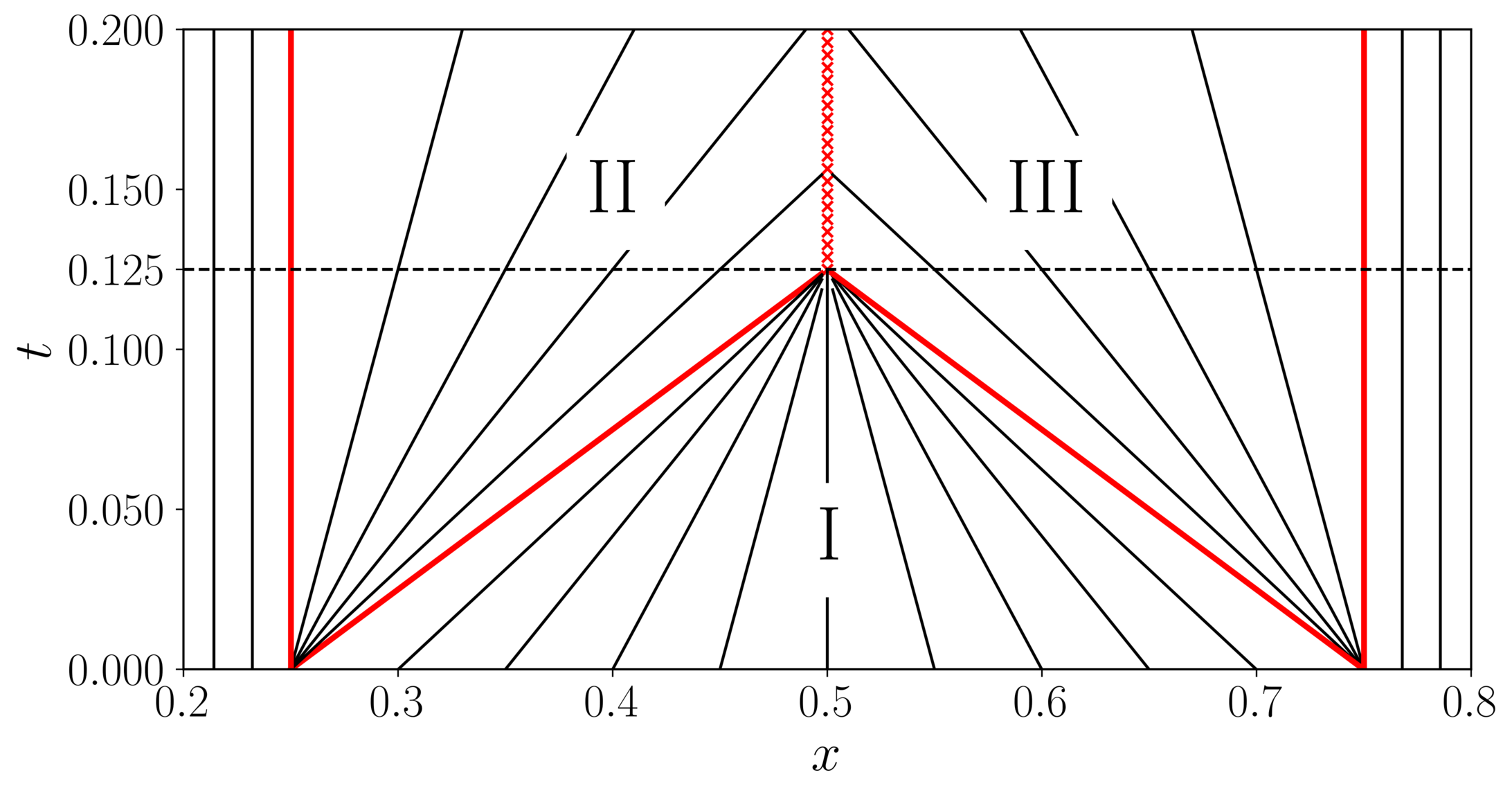}  
  \caption{Characteristics}
  \label{fig:N-wave_char}
\end{subfigure}
\begin{subfigure}[t]{.44\textwidth}
  \centering
  \includegraphics[width=\linewidth]{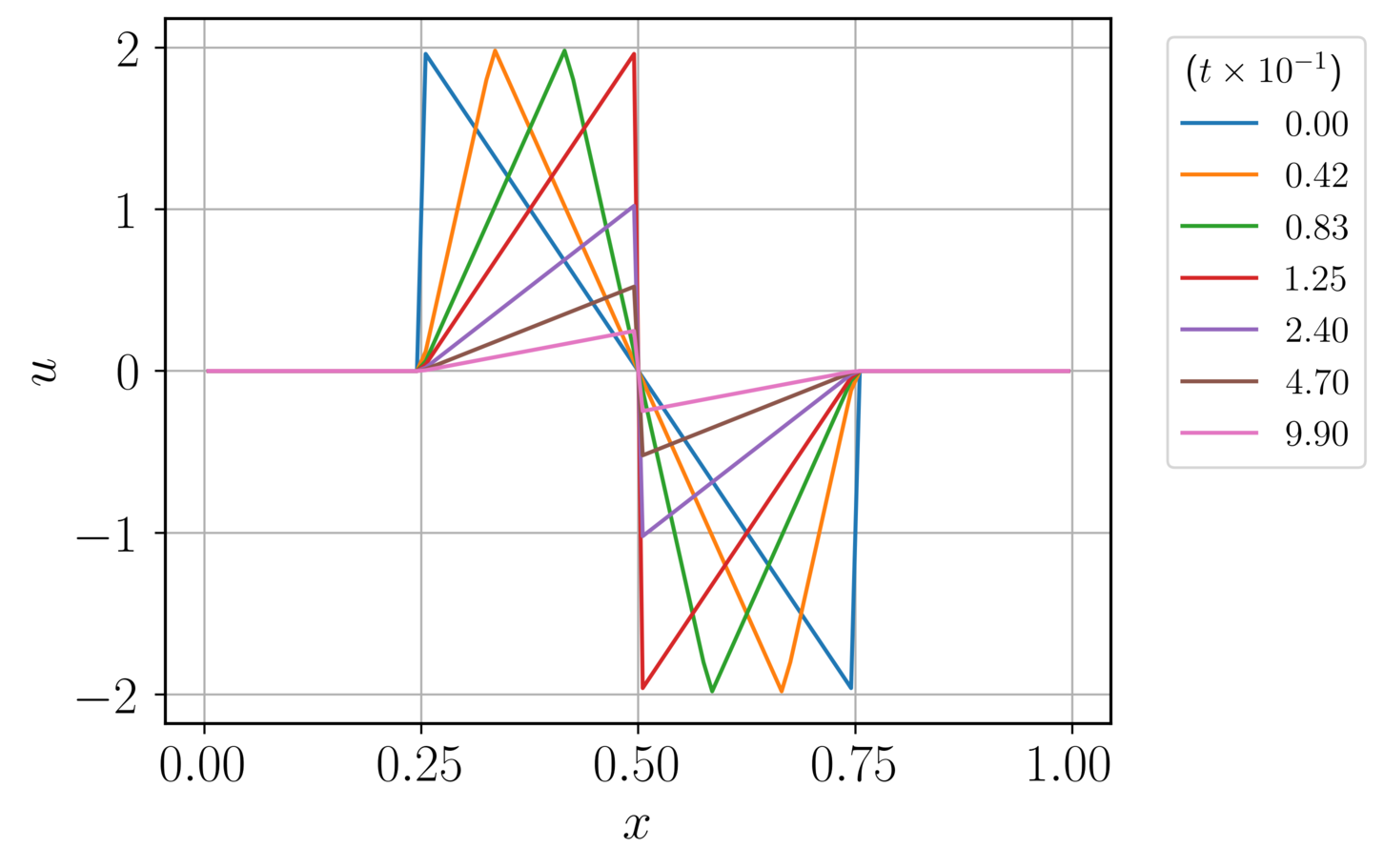}  
  \caption{$u(x,t)$}
  \label{fig:sub-bur_p5_ex}
\end{subfigure}
\caption{In Fig.~\ref{fig:N-wave_char}, slope of any characteristic in the $x-t$ plane is given by $1/u(x,t)$. Three distinct regions are delineated by the red-colored lines. The red coloured, vertical $x$-patterned line beyond $t_m$ at $x=0.5$ represents the shock.} 
\label{fig:app_bur_p5}
\end{figure}

Consider the profile given by
\begin{equation}
\begin{aligned}
u_0(x) =
\begin{cases}
    0 & \text{for } x < 0.25 \\
    -4h_0(x-0.5) & \text{for } 0.25 < x < 0.75 \\
    0 & \text{for } x > 0.75.
\end{cases}
\end{aligned}
\label{eq:app_N-wave_ini}
\end{equation} 
For the example considered in this paper ($h_0=2$), the characteristic plane has been shown in Fig.~\ref{fig:N-wave_char}. Each characteristic ray from $x^*$ at time $t=0$ is a straight line with slope $u(x^*)$  (Strang, pg.~591 \cite{GILBERT_STRANG}). Accordingly, for any $x_1,x_2\in(0.25,0.75)$, the characteristics from each of these points will be given by 
\begin{equation} \label{eq:x1_x2_form}t = \frac{1}{u_0(x_1)}(x - x_1); \quad t = \frac{1}{u_0(x_2)}(x - x_2).
\end{equation}
Since $u_0(x)$ is a linearly decreasing function in the considered domain, if $(x_m,t_m)$ represents the point where two such characteristics meet, substituting \eqref{eq:app_N-wave_ini} into \eqref{eq:x1_x2_form} and eliminating $t_m$ gives
$$-\frac{x_m-x_1}{4h_0(x_1-0.5)}=-\frac{x_m-x_2}{4h_0(x_2-0.5)} \quad \Rightarrow x_m=0.5.$$
Using \eqref{eq:x1_x2_form}, $t_m = \frac{1}{4h_0}$. Hence, for linearly decreasing functions $u_0(x)$, we find that $(x_m,t_m)$ is independent of $x_1$ and $x_2$ and all such characteristics converge at the single point in space-time. For $h_0=2$, $(x_m, t_m) = (0.5, 0.125)$.

Beyond $t=t_m$, shock formation takes place as a consequence of characteristics traveling in the opposite directions. In the regions $x<0.25$ and $x>0.75$, the rays travel parallel to time-axis. 
Furthermore, the jumps at $x=0.25$ and $x=0.75$ turn into expansion fans as shown in Fig.~\ref{fig:N-wave_char}. 
The solution for any point $\{(x,t): t>t_m, x\neq 0.5\}$ can be found based on the characteristics coming from the expansion fans. 
We define the following regions:
\begin{align*}
    &\text{I}:  \{(x,t): -2t+x - 0.25\geq0 \} \quad \cap \quad  \{ (x,t):2t+x - 0.75\leq 0\};\\
    &\text{II}: \{(x,t): \quad 0.25\leq x < 0.5\} \quad \cap \quad \{(x,t):-2t+x - 0.25\leq0\};\\
    &\text{III}: \{(x,t): \quad 0.5< x \leq 0.75 \} \quad \cap \quad \{(x,t) : 2t+x - 0.75\geq0\}.  
\end{align*}

The complete $u$ profile can thus be given by
\begin{equation}
u(x,t) =
\begin{cases}
        \dfrac{8(x-0.5)}{8\,t-1}  & \,\text{for Region I} \\[8pt] \dfrac{x-0.25}{t}
        & \,\text{for Region II} \\[8pt]
        \dfrac{x-0.75}{t} & \,\text{for Region III} \\
        0&\,\text{Otherwise. }\, \qquad  
     
\end{cases}
\label{eq:Nwave_ex_u}
\end{equation} 
The continuous $Y$ profile corresponding to $Y(0,t) = 0$ can be obtained as
\begin{equation}
Y(x,t) =
\begin{cases}
\begin{rcases}
           \dfrac{2 x^2 - x}{4t} + \dfrac{1}{32t} & \,\text{for }\,  0.25 \leq x < 0.25 + 2t \\[8pt] 
         \dfrac{4 (x^2 - x)}{8t - 1} + \dfrac{8t+3}{4(8t-1)}
        & \text{for } \,  0.25 + h_0\,t \leq x < 0.75 - h_0\,t \\[8pt]
        \dfrac{2x^2-3x}{4t} + \dfrac{9}{32t}  & \text{for } \,  0.75 - h_0\,t \leq x < 0.75 \\[8pt]
        0&\text{Otherwise }\\
        \end{rcases}
     0<t < t_m \\[10pt]
     \\
     \begin{rcases}
         \dfrac{2 x^2 - x}{4t} + \dfrac{1}{32t}
        & \,\,\text{for }\,  0.25 \leq x < 0.5 \\[8pt]
        \dfrac{2x^2-3x}{4t} + \dfrac{9}{32t}  & \,\,\text{for }\,  0.5 \leq x < 0.75 \\[8pt]
        0&\,\text{Otherwise }
        \end{rcases}
     t \geq t_m. 
\end{cases}
\label{eq:Nwave_ex_Y}
\end{equation}

\printbibliography
\end{document}